\documentclass{amsart}
\title{Minimal surfaces in pseudohermitian geometry}
\author{ Jih-Hsin Cheng }
\address[Cheng]
{Institute of Mathematics,
Academia Sinica, \newline%
\indent Nankang,
Taipei, Taiwan, 11529,
R.O.C.}
\email[]{cheng@math.sinica.edu.tw}%
\author{Jenn-Fang Hwang}
\address[Hwang]
{Institute of Mathematics,
Academia Sinica, \newline%
\indent Nankang,
Taipei, Taiwan, 11529,
R.O.C.}
\email[]{majfh@math.sinica.edu.tw}%
\author{\\ Andrea Malchiodi}
\address[Malchiodi]
{School of Mathematics,  Institute for Advanced Study,\newline%
\indent Princeton, NJ 08540, U.S.A.}
\email[]{malchiod@ias.edu}%
\author{Paul Yang}
\address[Yang]
{Department of Mathematics, Princeton University,\newline%
\indent Princeton, NJ 08544, U.S.A.}
\email[]{yang@Math.Princeton.EDU}%
\date{}

\begin{document}
\maketitle
\centerline{\bf Abstract}

We consider surfaces immersed in three-dimensional pseudohermitian
manifolds. We define the notion of (p-)mean curvature and of the
associated (p-)minimal surfaces, extending some concepts
previously given for the (flat) Heisenberg group. We interpret the
p-mean curvature not only as the tangential sublaplacian of a
defining function, but also as the curvature of a characteristic
curve, and as a quantity in terms of calibration geometry.

As a differential equation, the p-minimal surface equation is
degenerate (hyperbolic and elliptic). To analyze the singular set,
we formulate some {\em extension} theorems, which describe how the
characteristic curves meet the singular set. This allows us to
classify the entire solutions to this equation and to solve a
Bernstein-type problem (for graphs over the $xy$-plane) in the
Heisenberg group $H_1$. In $H_{1}$, identified with the Euclidean
space $R^{3}$, the p-minimal surfaces are classical ruled surfaces
with the rulings generated by Legendrian lines. We also prove a
uniqueness theorem for the Dirichlet problem under a condition on
the size of the singular set in two dimensions, and generalize to
higher dimensions without any size control condition.

We also show that there are no closed, connected, $C^{2}$ smoothly immersed
constant p-mean curvature or p-minimal surfaces of genus greater than one in the standard $S^{3}.$
This fact continues to hold when $S^{3}$
is replaced by a general spherical pseudohermitian 3-manifold.

\begin{center}
{\bf Contents}
\end{center}

\begin{enumerate}
  \item Introduction and statement of the results
  \item Surfaces in a 3-dimensional pseudohermitian manifold
  \item The singular set-proof of Theorem B
  \item A Bernstein-type theorem and properly embedded
p-minimal surfaces
  \item Comparison principle and uniqueness for the
Dirichlet problem
  \item  Second variation formula and area-minimizing
property
  \item Closed p-minimal surfaces in the standard $S^3$ and proof of Theorem E
\end{enumerate}
\ \ \ \ \ \ Appendix. Basic facts in pseudohermitian geometry

\bigskip

\section{{\bf Introduction and statement of the results}}

Minimal surfaces in a Riemannian manifold play an important role
in the study of topology and geometry of the ambient manifold. For
instance, the positive mass theorem originally was proved with the
aid of minimal surface theory ([SY]). In order to study the mass
in an analogous manner
or to formulate a boundary value problem for prescribing the
Webster scalar curvature on a domain with boundary in
pseudohermitian geometry, we find it necessary to formulate a
notion of mean curvature for a surface in a pseudohermitian
manifold.

Let $\Sigma$ be a surface in a 3-dimensional pseudohermitian
manifold (see the Appendix for some basic facts in pseudohermitian
geometry, and Section 2 for rigorous definitions). To every
non-singular point of $\Sigma$, we associate a vector field {\em
normal} to $\Sigma$, called the Legendrian normal. If for example
$\Sigma$ bounds some smooth set $\Omega$, then the p-area, in
analogy with the Riemannian case, is obtained as the variation of
the volume of $\Omega$ in the direction of the Legendrian normal,
see formulas (2.4) and (2.5).

The p-mean curvature $H$ of $\Sigma$ is in turn given by the first
variation of the p-area.  It is easy to see that $H$ equals the
negative subdivergence of the Legendrian normal $e_{2}.$ Suppose
$\psi $ is a defining function of $\Sigma $ such that
$e_{2}=\nabla _{b}\psi /|\nabla _{b}\psi |_{G}$ ($G$ is the Levi
metric. See the definition in Section 2). Then the p-mean
curvature equation and the p-minimal surface equation ($H \equiv
0$) read

\begin{equation}
\tag{pMCE}-div_{b}(\nabla _{b}\psi /|\nabla _{b}\psi |_{G})=H
\end{equation}
and
\begin{align}
 div_{b}(\nabla _{b}\psi /|\nabla _{b}\psi |_{G})=0,\tag{pMSE}
\end{align}

\noindent respectively.

Alternatively, having in mind the Gauss map, the mean curvature
can be defined in terms of the covariant derivative (with respect
to the pseudohermitian connection) along $\Sigma$ of the {\em
Legendrian tangent} $e_1$ to $\Sigma$, see (2.1).

Of course the case of the Heisenberg group as a pseudohermitian
manifold is one of the most important. Indeed it is the simplest
model example, and represents a blow-up limit of general
pseudo-hermitian manifolds. In the case of a smooth surface in the
Heisenberg group, our definitions coincide with those given in
([CDG]), ([DGN]) and ([Pau]). In particular these notions,
especially in the framework of geometrci measure theory, have been
used to study existence or regularity properties of minimizers for
the relative perimeter or extremizers of isoperimetric
inequalities, see (DGN), ([GN]), ([LM]), ([LR]), ([Pan]).  The
p-area can also be identified with the 3-dimensional spherical
Hausdorff measure of $\Sigma $ (see, e.g., [B], [FSS]).

In this paper, we study the subject mainly from the viewpoint of
partial differential equations and that of differential geometry.
Our basic results are the analysis of the singular set (see
Section 3). As consequences, we can prove a Bernstein-type theorem
(see Section 4 and Theorem A) and the nonexistence of closed
hyperbolic p-minimal surfaces (see Section 7 and Theorem E). We
also establish a comparison theorem (see Section 5) which is a
substitute for the maximum principle and may become a useful tool
in the subject.

For a p-minimal graph $(x,y,u(x,y))$ in the Heisenberg group
$H_{1}$, the above equation $(pMSE)$ reduces to

\begin{equation}
\tag{pMGE}
(u_{y}+x)^{2}u_{xx}-2(u_{y}+x)(u_{x}-y)u_{xy}+(u_{x}-y)^{2}u_{yy}=0
\end{equation}

\noindent by taking $\psi =z-u(x,y)$ on the nonsingular domain$.$
This is a degenerate (hyperbolic and elliptic) partial
differential equation. It is degenerate hyperbolic (on the
nonsingular domain) having only one characteristic direction (note
that a 2-dimensional hyperbolic equation has two characteristic
directions [Jo]). We call the integral curves of this
characteristic direction the characteristic curves. We show that
the p-mean curvature is the line curvature of a characteristic
curve. Therefore the characteristic curves of $(pMGE)$ are
straight lines. Moreover, the value of $u$ along a characteristic
curve is determined in a simple way (see (2.22), (2.23)).

The analysis of the singular set is necessary to characterize the
solutions. As long as the behavior of $H$ (consider $(pMCE)$ for a
graph in $H_1$) is not too bad (say, bounded), we show that the
singular set consists of only isolated points and smooth curves
(see Theorem B below). Under a quite weak growth condition on $H$,
a characteristic curve $\Gamma$ reaches a singular point $p_0$ in
a finite arc-length parameter and has an approximate tangent. From
the opposite direction, we find another characteristic curve
${\Gamma}^{\prime}$ reaching also $p_0$ with the opposite
approximate tangent. The union of $\Gamma$, $p_0$, and
${\Gamma}^{\prime}$ forms a smooth curve (see Corollary 3.6 and
Theorem 3.10). Making use of such {\em extension} theorems, we can
easily deal with the singular set, in order to study the Bernstein
problem. Namely, we study entire p-minimal graphs (a graph or a
solution is called entire if it is defined on the whole
$xy$-plane). The following are two families of such examples (cf.
[Pau]):

\begin{align}
&u=ax+by+c \, \, \mbox{(a plane with a,b,c being real constants);}& \tag{1.1}\\
&u=-abx^{2}+(a^{2}-b^{2})xy+aby^{2}+g(-bx+ay)& \tag{1.2} \\
&(a,b \mbox{ being real constants such that }a^{2}+b^{2}=1 \mbox{ and }
g\in C^{2}).&\nonumber
\end{align}

We have the following classification result (see Section 4).

\bigskip
\noindent{\bf Theorem A.} {\it (1.1) and (1.2) are the only entire
$C^{2}$ smooth solutions to the p-minimal graph equation
(pMGE).}
\bigskip

To prove Theorem A, we analyze the
characteristic curves and the singular set of a solution for the case $H=0$.
Observe that the
characteristic curves are straight lines which intersect at singular
points. Let $S(u)$ denote the singular set consisting of all points where
$u_{x}-y=0$ and $%
u_{y}+x=0$. Let $N(u)$ denote the $xy-$plane projection of the
negative Legendrian normal $-e_{2}.$ It follows that
$N(u)=(u_{x}-y,u_{y}+x)/D$ where
$D=\sqrt{(u_{x}-y)^{2}+(u_{y}+x)^{2}}$. On the nonsingular domain,
$(pMGE)$ has the form

$${div}N(u)(=H)=0,$$

\noindent where $div$ denotes the $xy$-plane divergence. The following
result gives a
local description of the singular set (see Section 3).

\bigskip
\noindent {\bf Theorem B.} {\it Let $\Omega$ be a domain in the
$xy-$plane. Let $u\in C^{2}(\Omega)$ be such that $div$$N(u)=H$ in
$\Omega\backslash S(u)$. Suppose $|H|\leq C\frac{1}{r}$ near a
singular point $p_{0}\in S(u)$ where $r(p)=|p-p_{0}|$ for $p\in
\Omega $ and $C$ is a positive constant. Then either $p_{0}$ is
isolated in $S(u)$ or there exists a small neighborhood of $p_{0}$
which intersects with $S(u)$ in exactly a $C^{1}$ smooth curve
past $p_{0}.$}
\bigskip

We show that the restriction on $H$ is necessary by giving a $%
C^{\infty }$ smooth counterexample. In additon the blow-up rate
$H=C\frac{1}{r}$ is realized by some natural examples (see Section
3). Theorem B follows from a characterization for a singular point
to be non-isolated (see Theorem 3.3).

When two characteristic lines meet at a point of a singular curve,
they must form a straight line (see Lemma 4.4). So we can describe all
possible configurations of characteristic lines as if singular curves are
not there. It turns out that there are only two possible configurations of
characteristic lines. Either all characteristic lines intersect at one
singular point or they are all parallel. In the former case, we are led to
the solution $(1.1)$ while for the latter case, $(1.2)$ is the only possible
solution.

The characteristic curves on a p-minimal surface are the
Legendrian geodesics (see (2.1)). Since the Legendrian geodesics
in $H_1$ are nothing but straight lines, a general complete
p-minimal surface is a complete ruled surface generated by
Legendrian rulings. We will discuss this and point out that a
known complete embedded non-planar p-minimal surface has no
singular points (characteristic points in the terminology of some
other authors) after $(4.10)$ in Section 4.

Since $(pMGE)$ is also a degenerate elliptic partial differential
equation, one can use non-degenerate elliptic equations to
approximate it. With this regularization method, Pauls ([Pau])
obtained a $W^{1,p}$ Dirichlet solution and showed that such
surfaces are the X-minimal surfaces in the sense of Garofalo and
Nhieu ([GN]). In general the solution to the Dirichlet problem may
not be unique. However, we can still establish a uniqueness
theorem by making use of a structural equality of "elliptic" type
(Lemma 5.1). More generally we have the following comparison
principle.

\bigskip
\noindent {\bf Theorem C.} {\it For a bounded domain $\Omega $ in
the $xy$-plane, let $u,v\in C^{2}(\Omega )\cap
C^{0}(\bar{\Omega})$ satisfy $div$$N(u)\geq$ $div$$N(v)$ in
$\Omega \backslash S$ and $u\leq v$ on $\partial \Omega $ where
$S=S(u)\cup S(v).$ Suppose $\mathcal{H}_{1}(\bar{S}),$ the
1-dimensional Hausdorff measure of $\bar S$, vanishes. Then $u\leq
v$ in $\Omega .$}
\bigskip

As an immediate consequence of Theorem C, we have the following
uniqueness result for the Dirichlet problem of $(pMGE)$ (see Section 5).

\bigskip
\noindent {\bf Corollary D.} {\it For a bounded domain $\Omega$ in
the $xy$-plane, let $u,v\in C^{2}(\Omega )\cap
C^{0}(\bar{\Omega})$ satisfy $div$$N(u)=div$$N(v)=0$ in $\Omega
\backslash S$ and $u=v$ on $\partial \Omega $ where $S=S(u)\cup
S(v).$ Suppose ${\mathcal H}_{1}(\bar{S}),$ the 1-dimensional
Hausdorff measure of $\bar S$, vanishes. Then $u=v$ in $\Omega .$}
\bigskip

We remark that the condition on ${\mathcal H}_{1}(\bar{S})$ in
Corollary D is necessary. A counterexample is given in [Pau] with
${\mathcal H}_{1}(\bar{S})\neq 0.$ We generalize Theorem C to
higher dimensions and for a class of general $N$ (see Section 5).
It is noticeable that we do not need the condition on the size control
of the singular set for the higher dimensional version of Theorem C
(see Theorem C').

We also study closed p-minimal surfaces in the standard
pseudohermitian 3-sphere. A characteristic curve of such a
p-minimal surface is part of a Legendrian great circle (see Lemma
7.1). Using this fact, we can describe the {\em extension}
theorems (Corollary 3.6, Theorem 3.10 or Lemma 7.3) in terms of
Legendrian great circles, and hence give a direct proof of the
nonexistence of hyperbolic p-minimal surfaces embedded in the
standard pseudohermitian 3-sphere (part of Corollary F). Then we
generalize to the situation that the ambient pseudohermitian
3-manifold is spherical and the immersed surface has bounded
p-mean curvature (see Section 7).

\bigskip
\noindent {\bf Theorem E.} {\it Let $M$ be a spherical
pseudohermitian 3-manifold. Let $\Sigma$ be a closed,
connected surface, $C^{2}$ smoothly immersed in $M$ with bounded p-mean
curvature. Then the genus of $\Sigma$ is less than or equal to 1.
In particular, there are no constant p-mean curvature or p-minimal
surfaces $\Sigma$ of genus greater than one in $M.$}
\bigskip

There are many examples of spherical CR manifolds, and there have
been many studies in this direction (e.g., [BS], [KT], [FG], [CT],
[S] ). We speculate that Theorem E might imply a topological
constraint on a spherical CR 3-manifold. The idea of the proof for
Theorem E goes as follows. A spherical pseudohermitian manifold is
locally the Heisenberg group with the contact form being a
multiple of the standard one. So locally near a singular point,
$\Sigma$ is a graph in $H_1$ having bounded p-mean curvature with
respect to the standard contact form. We can then apply Lemma 3.8
to conclude that the characteristic line field has index $+1$ at
every isolated singular point. Therefore the Euler characteristic
number (equals the index sum), hence the genus, of $\Sigma$ has
constraint in view of the Hopf index theorem for a line field.

\bigskip
\noindent {\bf Corollary F.} {\it There are no closed, connected,
$C^{2}$ smoothly immersed constant p-mean curvature or p-minimal
surfaces of genus $\geq 2$ in the standard pseudohermitian
3-sphere.}
\bigskip

Note that in the standard Euclidean 3-sphere, there exist many
closed $C^{\infty }$ smoothly embedded minimal surfaces of genus $\geq 2$ ([La]).

On a surface in a pseudohermitian 3-manifold, we define an
operator, called the tangential sublaplacian. The p-mean curvature
is related to this operator acting on coordinate functions (see
$(2.19a),(2.19b),(2.19c)$) for a graph in $H_{1}$. We therefore
obtain a "normal form" (see (2.20)) of $(pMGE)$. We also interpret
the notion of p-mean curvature in terms of calibration geometry.
{From} this we deduce the area-minimizing property for p-minimal
surfaces (see Proposition 6.2). Since the second variation formula
is important for later development, we derive it and discuss the
stability of a p-minimal surface in Sections 6 and 7.

We remark that, in the preprint ([GP1]), the authors claim
the vertical planes are the only complete p-minimal graphs having
no singular points (non-characteristic complete minimal graphs in
their terminology). This is faulty. For instance, $y=xz$ is a
complete (in fact, entire) p-minimal graph over the $xz$-plane and
has no singular points. In [CH], two of us classify all the entire
p-minimal graphs over any plane among other things. After our paper
was completed, we were informed that the above claim had been 
corrected and some similar results are also obtained in the new 
preprint ([GP2]).

\bigskip

\noindent {\bf Acknowledgments.} We would like to thank Ai-Nung
Wang for informative discussion of  Monge's ([Mo]) third order
equation for ruled surfaces (see Section 4). The first author
would also like to thank Yng-Ing Lee for showing him some basic
facts in calibration geometry. We began this research during the
first author's visit at the IAS, Princeton in the 2001-2002
academic year. He would therefore like to thank the faculty and
staff there as their hospitality has greatly facilitated this
collaboration.

\bigskip

\section{{\bf Surfaces in a 3-dimensional pseudohermitian manifold}}

Let $(M,J,\Theta )$ be a 3-dimensional oriented pseudohermitian manifold
with a $CR$ structure $J$ and a global contact form $\Theta $ (see the
Appendix). Let $\Sigma $ be a surface contained in $M.$ The singular set $%
S_{\Sigma }$ consists of those points where $\xi $ coincides with
the tangent bundle $T\Sigma $ of $\Sigma $. It is easy to see that
$S_{\Sigma }$ is a closed set. On the nonsingular (open) set
$\Sigma \backslash S_{\Sigma } $, we call the leaves of the
1-dimensional foliation $\xi \cap T\Sigma $ the {\em
characteristic curves}. These curves will play an important role
in our study. On $\xi $, we can associate a natural metric
$G=\frac{1}{2}d\Theta (\cdot ,J\cdot ),$ called the Levi metric.
For a vector $v\in \xi $, we define the length of $v$ by
$|v|_{G}=(G(v,v))^{1/2}.$ With respect to this
metric, we can take a unit vector field $e_{1}\in \xi \cap T\Sigma $ on $%
\Sigma \backslash S_{\Sigma }$, called the characteristic field. Also
associated to $(J,\Theta )$ is the so-called pseudohermitian connection,
denoted as $\nabla ^{p.h.}$ (see (A.2) in the Appendix). We can define a
notion of mean curvature for $\Sigma $ in this geometry as follows. Since $%
\nabla ^{p.h.}$ preserves the Levi metric $G$, $\nabla ^{p.h.}e_{1}$ is
perpendicular to $e_{1}$ with respect to $G.$ On the other hand, it is
obvious that $G(e_{1,}e_{2})=0$ where $e_{2}=Je_{1}.$ We call $e_{2}$ the
Legendrian normal or Gauss map. So we have

\begin{align}
    \nabla _{e_{1}}^{p.h.}e_{1}=He_{2}  \tag{2.1}
\end{align}

\noindent for some function $H$. We call $H$ the p(pseudohermitian)-mean curvature of $%
\Sigma $. Note that if we change the sign of $e_{1},$ then $e_{2}$ and $H$
change signs accordingly. If $H=0$, we call $\Sigma $ a p-minimal surface.
In this situation the characteristic curves are nothing but Legendrian
(i.e., tangent to $\xi $) geodesics with respect to the pseudohermitian
connection.

We are going to give a variational formulation for the p-mean
curvature $H$. First let us find a candidate area integral. Suppose $\Omega $
is a smooth domain in $M$ with boundary $\partial \Omega =\Sigma $.
Consider $V(\Omega )$, the volume of $\Omega ,$ given by

$$V(\Omega )=\frac{1}{2}\int_{\Omega }\Theta \wedge
d\Theta$$

\noindent ($\frac{1}{2}$ is a normalization constant. For $\Omega \subset H_{1}$, this
volume is just the usual Euclidean volume). Take Legendrian fields $%
e_{1},e_{2}=Je_{1}\in \xi$, orthonormal with respect to $G$, wherever
defined in a neighborhood of $\Sigma $ (note that we do not require $e_{1}$
to be characteristic, i.e. tangent along $\Sigma $ here)$.$ We consider a
variation of the surface $\Sigma $ in the direction $fe_{2}$ where $f$ is a
suitable function with compact support in $\Sigma \backslash S_{\Sigma }$.
The vector field $fe_{2}$ generates a flow $\varphi _{t}$ for $t$ close to $%
0 $. We compute

\begin{align}
    \delta _{fe_{2}}V(\Omega )&= \frac{d}{dt}|_{t=0}V(\varphi _{t}(\Omega ))
=\frac{1}{2}\frac{d}{dt}|_{t=0}\int_{\Omega
}\varphi _{t}^{\ast }(\Theta \wedge d\Theta )  \tag{2.2} \\
&= \frac{1}{%
2}\int_{\Omega }L_{fe_{2}}(\Theta \wedge d\Theta ). \nonumber
\end{align}

It follows from the formula $L_{X}=d\circ i_{X}+i_{X}\circ d$ ($%
i_{X} $ denotes the interior product in the direction $X$) that

\begin{align}
    L_{fe_{2}}(\Theta \wedge d\Theta
)=d\circ i_{fe_{2}}(\Theta \wedge d\Theta ). \tag{2.3}
\end{align}

\noindent Substituting $(2.3)$ in $(2.2)$ and making use of Stokes' theorem, we obtain

\begin{align}
    \delta _{fe_{2}}V(\Omega )=\frac{1}{2}%
\int_{\Sigma }i_{fe_{2}}(\Theta \wedge d\Theta )=\int_{\Sigma }f \Theta
\wedge e^{1}. \tag{2.4}
\end{align}

\noindent Here $e^{1}$ together with $e^{2},\Theta $ form a dual basis of $%
(e_{1},e_{2},T)$ where $T$ is the Reeb vector field (uniquely determined by $%
\Theta (T)=1$ and $i_{T}d\Theta =0).$ Note that $d\Theta =2e^{1}\wedge e^{2}$
(see $(A.1r)$). For $e_{1}$ being a characteristic field, we define the
p-area of a surface $\Sigma $ to be the surface integral of the 2-form $%
\Theta \wedge e^{1}:$

\begin{align}
    p-Area(\Sigma
)=\int_{\Sigma }\Theta \wedge e^{1}. \tag{2.5}
\end{align}

\noindent Note that $\Theta \wedge e^{1}$ continuously extends over the singular set $%
S_{\Sigma }$ and vanishes on $S_{\Sigma }.$ In fact, we can write $e^{1}$
with respect to a dual orthonormal basis $\{\hat{e}^{1},\hat{e}^{2}\}$ of $%
\xi ,$ which is smooth near a singular point, say $p_{0}$, as follows: $%
e^{1}=\cos \beta \hat{e}^{1}$ + sin$\beta \hat{e}^{2}.$ Here $\beta $ may
not be continuous at $p_{0}.$ Now $\Theta \wedge \hat{e}^{1}$ and $\Theta
\wedge \hat{e}^{2}$ tend to $0$ on $\Sigma $ as $p\in \Sigma $ tends to $%
p_{0}$ since $\Theta $ vanishes on $T_{p_{0}}\Sigma =\xi _{p_{0}}.$ It
follows that $\Theta \wedge e^{1}$ tends to $0$ on $\Sigma $ as $p\in \Sigma
$ tends to $p_{0}$ since $\cos \beta $ and $\sin \beta $ are bounded by $1.$

We can recover the p-mean curvature $H$ from the first variation
formula of the p-area functional $(2.5).$ We compute

\begin{align}
    \delta _{fe_{2}}\int_{\Sigma }\Theta
\wedge e^{1}=\int_{\Sigma }L_{fe_{2}}(\Theta \wedge e^{1})=\int_{\Sigma
}i_{fe_{2}}\circ d(\Theta \wedge e^{1}). \tag{2.6}
\end{align}

\noindent Here we have used the formula $L_{X}=d\circ i_{X}+i_{X}\circ d$ and the
condition that $f$ is a function with compact support away from the singular
set and the boundary of $\Sigma $. From the equations $d\Theta =2e^{1}\wedge
e^{2}$ and $de^{1}=-e^{2}\wedge \omega $ mod $\Theta $ (see (A.1r), (A.3r)),
we compute

\begin{align}
    d(\Theta \wedge e^{1})=d\Theta \wedge
e^{1}-\Theta \wedge de^{1}=\Theta \wedge e^{2}\wedge \omega . \tag{2.7}
\end{align}

\noindent Substituting $(2.7)$ into $(2.6)$, we obtain by the definition of the interior
product that

\begin{align}
    \delta _{fe_{2}}\int_{\Sigma }\Theta \wedge
e^{1}&=\int_{\Sigma }f(-\Theta \wedge \omega +\omega (e_{2})\Theta \wedge
e^{2}) \tag{2.8} \\
&=\int_{\Sigma }-f\omega
(e_{1})\Theta \wedge e^{1} \nonumber \\
(\Theta \wedge e^{2}&=0 \mbox{ on }
\Sigma  \mbox{ since } e_{1} \mbox{ is tangent along } \Sigma ) \nonumber \\
&=-\int_{\Sigma }fH\Theta \wedge e^{1}. \nonumber
\end{align}

\noindent In the last equality, we have used the fact that $H=\omega (e_{1})$ (obtained by
comparing (2.1) with (A.2r)). Similarly we can also compute the first
variation of (2.5) with respect to the field $gT$ where $g$ is a function
with compact support away from the singular set and the boundary of $\Sigma
. $ Together with $(2.8)$, the result reads

\begin{align}
\delta _{fe_{2}+gT}\int_{\Sigma }\Theta \wedge
e^{1}=-\int_{\Sigma }(f-\alpha g)H\Theta \wedge e^{1}. \tag{$2.8^{\prime}$}
\end{align}

\noindent Here we define the function $\alpha $ on $\Sigma \backslash S_{\Sigma }$ such that $%
\alpha e_{2}+T\in T\Sigma .$ We leave the deduction of $%
(2.8^{\prime})$ to the reader. Let $\psi $ be a defining function of $%
\Sigma $. It follows that the subgradient $\nabla _{b}\psi =(e_{1}\psi
)e_{1}+(e_{2}\psi )e_{2}=(e_{2}\psi )e_{2}$ since $e_{1}\in T\Sigma $, hence
$e_{1}\psi =0.$ So $\nabla _{b}\psi /|\nabla _{b}\psi |_{G}=e_{2}$ (change
the sign of $\psi $ if necessary)$.$ Next we compute the subdivergence of $%
e_{2}.$ Since $\nabla ^{p.h.}$ preserves the Levi metric $G,$ we can easily
obtain $G(\nabla ^{p.h.}e_{2},e_{2})=0$ and $G(\nabla
_{e_{1}}^{p.h.}e_{2},e_{1})=-G(e_{2},\nabla _{e_{1}}^{p.h.}e_{1})=-H$ by $%
(2.1).$ Therefore $div_{b}(e_{2})$ $\equiv $ $G(\nabla
_{e_{1}}^{p.h.}e_{2},e_{1})$ $+$ $G(\nabla _{e_{2}}^{p.h.}e_{2},e_{2})$ $=$ $%
-H.$ We have derived the p-mean curvature equation ($pMCE$) and the p-minimal
surface equation ($pMSE$):

\begin{align}
    -div_{b}(\nabla _{b}\psi /|\nabla _{b}\psi |_{G})=H \tag{pMCE}
\end{align}
and
\begin{align}
div_{b}(\nabla _{b}\psi /|\nabla _{b}\psi |_{G})=0, \tag{pMSE}
\end{align}

\noindent respectively.

For a graph $(x,y,u(x,y))$ in the 3-dimensional Heisenberg group $%
H_{1}$, we can take $\psi =z-u(x,y).$ Then (at a nonsingular point) up to
sign, $e_{1}$ is uniquely determined by the following equations:

\begin{align}
     &&\Theta _{0}\equiv dz+xdy-ydx=0, \tag{2.9a} \\
     &&d\psi =d(z-u(x,y))=0, \tag{2.9b} \\
     &&G(e_{1},e_{1})=\frac{1}{2}
d\Theta _{0}(e_{1},Je_{1})=1. \tag{2.9c}
\end{align}

\noindent Using $(2.9a)$, $(2.9c),$ we can write $e_{1}=f\hat{e}_{1}+g\hat{e}_{2}$
with $f^{2}+g^{2}=1,$ in which $\hat{e}_{1}=\frac{\partial }{\partial x}+y%
\frac{\partial }{\partial z},\hat{e}_{2}=\frac{\partial }{\partial y}-x\frac{%
\partial }{\partial z}$ are standard left-invariant Legendrian vector fields
in $H_{1}$ (see the Appendix)$.$ Applying $(2.9b)$ to this expression, we
obtain $(u_{x}-y)f+(u_{y}+x)g=0.$ So $(f,g)=\pm
(-(u_{y}+x),u_{x}-y)/[(u_{x}-y)^{2}+(u_{y}+x)^{2}]^{1/2}$ (positive sign so
that $\nabla _{b}\psi /|\nabla _{b}\psi |_{G}=e_{2}=-[(u_{x}-y)\hat{e}%
_{1}+(u_{y}+x)\hat{e}_{2}]/D$ where $D=[(u_{x}-y)^{2}+(u_{y}+x)^{2}]^{1/2}$%
). Now from $(pMCE)$ we obtain a formula for $H$ through a direct
computation:

\begin{align}
    H=D^{-3}%
\{(u_{y}+x)^{2}u_{xx}-2(u_{y}+x)(u_{x}-y)u_{xy}+(u_{x}-y)^{2}u_{yy}\}.\tag{2.10}
\end{align}

\noindent At a nonsingular point, the equation $(pMSE)$ reduces to the p-minimal graph
equation $(pMGE):$

\begin{align}
    \tag{pMGE}
(u_{y}+x)^{2}u_{xx}-2(u_{y}+x)(u_{x}-y)u_{xy}+(u_{x}-y)^{2}u_{yy}=0.
\end{align}

\noindent In fact, if $u$ is $C^{2}$ smooth, the p-mean curvature $H$ in (2.10)
vanishes on the nonsingular domain (where $D\neq 0$) if and only if $(pMGE)$
holds on the whole domain. We can also compute $e^{1}=D^{-1}%
\{-(u_{y}+x)dx+(u_{x}-y)dy\}$ and express the p-area 2-form as follows:

\begin{align}
    \tag{2.11} \Theta \wedge e^{1}=Ddx\wedge
dy=[(u_{x}-y)^{2}+(u_{y}+x)^{2}]^{1/2}dx\wedge dy.
\end{align}

\noindent At a singular point, the contact form $\Theta $ is proportional to $d\psi $ (see $%
(2.9a),(2.9b)$). Therefore $u_{x}-y=0,u_{y}+x=0$ describe the $xy-$plane
projection $S(u)$ of the singular set $S_{\Sigma }$:

\begin{align}
    \tag{S} S(u)=\{(x,y)\in
R^{2}:u_{x}-y=0,u_{y}+x=0\}.
\end{align}

\noindent From $(2.11)$ we see that the p-area form $\Theta \wedge e^{1}$ is
degenerate on $S(u)$ or $S_{\Sigma}$.

Let $e_{2}$ be the Legendrian normal of a family of deformed
surfaces foliating a neighborhood of $\Sigma .$ We define the tangential
subgradient $\nabla _{b}^{t}$ of a function $f$ defined near $\Sigma $ by
the formula: $\nabla _{b}^{t}f=\nabla _{b}f-G(\nabla _{b}f,e_{2})e_{2}$ (see
$(A.8)$ for the definition of $\nabla _{b}f$) and the tangential
pseudohermitian connection $\nabla ^{t.p.h.}$ of a Legendrian (i.e. in $\xi $%
) vector field $X$ by $\nabla ^{t.p.h.}X=\nabla ^{p.h.}X-G(\nabla
^{p.h.}X,e_{2})e_{2}.$ Then we define the tangential sublaplacian $\Delta
_{b}^{t}$ of $f$ by

\begin{align}
\tag{2.12} \Delta _{b}^{t}f=div_{b}^{t}(\nabla _{b}^{t}f)
\end{align}

\noindent where $div_{b}^{t}(X)$ is defined to be the trace of $\nabla ^{t.p.h.}X$
considered as an endomorphism of $\xi $ :$v\rightarrow \nabla
_{v}^{t.p.h.}X. $ Now for an orthonormal basis $e_{1},e_{2}$ of $\xi $ with
respect to $G,$ we have

\begin{align}
    \tag{2.13} div_{b}^{t}(X)=G(\nabla
_{e_{1}}^{t.p.h.}X,e_{1})+G(\nabla _{e_{2}}^{t.p.h.}X,e_{2})=G(\nabla
_{e_{1}}^{t.p.h.}X,e_{1})
\end{align}

\noindent since $\nabla _{e_{2}}^{t.p.h.}X$ is proportional to $e_{1}.$ On the other
hand, we write $\nabla _{b}f=(e_{1}f)e_{1}+(e_{2}f)e_{2}$ (cf. $(A.8r)$)$.$
It follows that $\nabla _{b}^{t}f=(e_{1}f)e_{1}.$ Setting $X=\nabla
_{b}^{t}f=(e_{1}f)e_{1}$ in $(2.13)$ and noting that $\nabla
_{e_{1}}^{t.p.h.}e_{1}=0$ by $(2.1)$ gives $div_{b}^{t}(\nabla
_{b}^{t}f)=(e_{1})^{2}f+(e_{1}f)G(\nabla
_{e_{1}}^{t.p.h.}e_{1},e_{1})=(e_{1})^{2}f$. Substituting this in $(2.12),$
we obtain

\begin{align}
    \tag{2.14} \Delta _{b}^{t}f=(e_{1})^{2}f.
\end{align}

\noindent Note that $(2.14)$ holds for a general surface $\Sigma $ contained in an
arbitrary pseudohermitian 3-manifold. We also note that
$\Delta _{b}^{t}+2\alpha e_{1}$ is
self-adjoint with respect to the p-area form $\Theta \wedge e^{1}$
as shown by the following integral formula:

\begin{align}
&\int_{\Sigma}[f\Delta _{b}^{t}g-g\Delta _{b}^{t}f]\Theta \wedge e^{1}
=2\int_{\Sigma}[g(e_{1}f)-f(e_{1}g)]\alpha\Theta \wedge e^{1}&\nonumber
\end{align}

\noindent for smooth functions $f,g$ with compact support away from the
singular set. The proof is left to the reader (Hint: observe that the adjoint
of $e_{1}$ is $-e_{1}-2\alpha$ by noting that $\Theta\wedge e^{2}=0$ and
$e^{1}\wedge e^{2}=\alpha e^{1}\wedge \Theta$ on the surface). When $\Sigma $
is a graph $%
(x,y,u(x,y)) $ in $H_{1},$ we can relate $(e_{1})^{2}f$ for $f=x,y,$ or $u$
to the p-mean curvature $H.$ Denote the projection of $-e_{2}$ ($-e_{1},$
respectively) onto the $xy$-plane by $N(u)$ or simply $N$ ($N^{\perp }(u)$
or simply $N^{\perp },$ respectively). Recall that $e_{1}=$ $[-(u_{y}+x)%
\hat{e}_{1}+(u_{x}-y)\hat{e}_{2}]/D$ where $%
D=[(u_{x}-y)^{2}+(u_{y}+x)^{2}]^{1/2}$ (see the paragraph after $(2.9)$)$.$
So $N^{\perp }=[(u_{y}+x)\partial _{x}-(u_{x}-y)\partial _{y}]/D.$ Write $%
(u_{y}+x)D^{-1}=\sin \theta ,$ $(u_{x}-y)D^{-1}=\cos \theta $ for some local
function $\theta .$ Then we can write

\begin{align}
    \tag{2.15a} &N=(\cos \theta )\partial _{x}+(\sin \theta
)\partial _{y}, \\
\tag{2.15b} &N^{\perp }=(\sin \theta )\partial _{x}-(\cos
\theta )\partial _{y}.
\end{align}

\noindent First note that from $(pMCE)$ we can express the p-mean curvature $H$ as
follows:

\begin{align}
    \tag{2.16} H=div N=(\cos \theta )_{x}+(\sin \theta
)_{y}=-(\sin \theta )\theta _{x}+(\cos \theta )\theta _{y}.
\end{align}

\noindent Now starting from $(2.15b)$ and using $(2.16),$ we can deduce

\begin{align}
    \tag{2.17} (N^{\perp })^{2}&=\sin ^{2}\theta \partial _{x}^{2}-2\sin
\theta \cos \theta \partial _{x}\partial _{y}+\cos ^{2}\theta \partial
_{y}^{2} \\
&\ \ -(\cos \theta )H\partial _{x}-(\sin \theta )H\partial _{y}. \nonumber
\end{align}

\noindent On the other hand, we can write $(2.10)$ in the following form:

\begin{align}
    \tag{2.18} H=D^{-1}(\sin ^{2}\theta \partial
_{x}^{2}u-2\sin \theta \cos \theta \partial _{x}\partial _{y}u+\cos
^{2}\theta \partial _{y}^{2}u).
\end{align}

\noindent Applying $(2.17)$ to $x,y,u(x,y),$ respectively and making use of $(2.18)$,
we obtain

\begin{align}
    \tag{2.19a} \Delta _{b}^{t}x &=(e_{1})^{2}x=(N^{\perp
})^{2}x=-(\cos \theta )H, \\
\tag{2.19b} \Delta _{b}^{t}y&=(e_{1})^{2}y=(N^{\perp
})^{2}y=-(\sin \theta )H, \\
\tag{2.19c} \Delta _{b}^{t}u&=(e_{1})^{2}u=(N^{\perp
})^{2}u=D^{-1}(xu_{y}-yu_{x}+x^{2}+y^{2})H \\
&=H(x\sin \theta
-y\cos \theta )=H\det \left[
\begin{array}{cc}
x & y \\
\cos \theta & \sin \theta%
\end{array}%
\right]  \nonumber
\end{align}

\noindent (here $"\det "$ means {\em determinant}). Formula
$(2.19c)$ gives the following {\em normal form} of $(pMGE):$

\bigskip

{\bf Lemma 2.1.} \textsl{For a }$C^{2}$ \textsl{smooth} \textsl{p-minimal graph }%
$u=u(x,y)$\textsl{\ in the 3-dimensional Heisenberg group }$H_{1},$\textsl{\
we have the equation}

\begin{align}
    \tag{2.20} \textsl{} \Delta
_{b}^{t}u=(e_{1})^{2}u=(N^{\perp })^{2}u=0
\end{align}

\noindent \textsl{on the nonsingular domain.}

\bigskip

Also from $(2.19a)$ and $(2.19b),$ we have $\ \Delta _{b}^{t}x=\Delta
_{b}^{t}y=0$ on a p-minimal graph $\Sigma $ $=$ $\{(x,y,u(x,y))\}.$ Together
with $(2.20)$, we have $\Delta _{b}^{t}(x,y,z)$ $\equiv $ $(\Delta
_{b}^{t}x,\Delta _{b}^{t}y,\Delta _{b}^{t}z)$ $=$ $(0,0,0)$ (i.e., $\Delta
_{b}^{t}$ annihilates the coordinate functions) on $\Sigma .$ This is a
property analogous to that for (Euclidean) minimal surfaces in $R^{3}$
([Os]). In general, we have 

$$\Delta _{b}^{t}(x,y,z)=He_{2}$$.

We will often call the $xy-$plane projection of characteristic curves
for a graph $(x,y,u(x,y))$ in $H_{1}$ still characteristic curves if no
confusion occurs. Note that the integral curves of $N^{\perp }$ are just the
$xy-$plane projection of integral curves of $e_{1}.$ So they are
characteristic curves. Along a characteristic curve $(x(s),y(s))$ where $s$
is a unit-speed parameter, we have the equations

\begin{align}
 \tag{2.21}
 \frac{dx}{ds}=\sin \theta ,\ \ \frac{dy}{ds}=-\cos \theta
\end{align}

\noindent by (2.15b). Noting that $u_{x}=(\cos \theta )D+y,$ $u_{y}=(\sin \theta )D-x,$
we compute

\begin{align}
 \tag{2.22}
 \frac{du}{ds}=u_{x}\frac{dx}{ds}+u_{y}\frac{dy}{ds}=[(\cos
 \theta )D+y]\sin \theta +[(\sin \theta )D-x](-\cos \theta )\\
=x\cos \theta +y\sin \theta ,\nonumber\\
\tag{2.23} \frac{d\theta }{ds}=\theta _{x}\frac{dx}{ds}+\theta _{y}%
\frac{dy}{ds}=\theta _{x}\sin \theta -\theta _{y}\cos \theta =-H
\end{align}

\noindent by (2.16). In general, we can consider $H$ as a function of $x,y,u,\theta $
in view of the O.D.E. system $(2.21),(2.22)$ and $(2.23).$ From $(2.21)$ and
$(2.23),$ we compute

\begin{align}
\tag{2.24} \frac{d^{2}x}{ds^{2}}=-H\cos \theta , \ \ \frac{d^{2}y}{%
ds^{2}}=-H\sin \theta .
\end{align}

\noindent Observe that $(\cos \theta ,\sin \theta )$ is the unit normal to the unit
tangent $(\frac{dx}{ds},\frac{dy}{ds})=(\sin \theta ,$ $-\cos \theta ).$ So $-H$
is just the curvature of a characteristic curve. In particular, when $H=0,$
characteristic curves are nothing but straight lines or line segments (see
Proposition 4.1 and Corollary 4.2).

\bigskip

\section{{\bf The singular set-proof of Theorem B}}

Let $\Omega $ be a domain (connected open subset) in the $xy$-plane,
and let $u\in C^{2}(\Omega )$. Let $\Sigma =\{(x,y,u(x,y)$ $|$ $(x,y)\in
\Omega \}$ $\subset $ $H_{1}$. In this section, we want to analyze $S(u)$
(still called the singular set), the $xy-$plane projection of the singular
set $S_{\Sigma }$, for the graph $\Sigma $.

First for $a,b\in R,$ $a^{2}+b^{2}=1,$ we define $F_{a,b}\equiv
a(u_{x}-y)+b(u_{y}+x)$ and $\Gamma _{a,b}\equiv \{(x,y)\in \Omega $ $|$ $%
F_{a,b}(x,y)=0\}.$ Then it is easy to see that $\Gamma _{a,b}=S(u)\cup
\{(x,y)\in \Omega $ $|$ $N(u)(x,y)$ $=$ $\pm (b,-a)\}.$ Let

\begin{align}
U = \left[
\begin{array}{cc}
u_{xx} & u_{xy}-1 \\
u_{xy}+1 & u_{yy}%
\end{array}%
\right] .\tag{U}
\end{align}

\bigskip

{\bf Lemma 3.1.} \textsl{Let }$u\in C^{2}(\Omega )$\textsl{\ and let }$p_{0}$%
\textsl{\ }$\in $\textsl{\ }$S(u)\subset \Omega .$\textsl{\ Then there
exists a small neighborhood of }$p_{0}$\textsl{, whose intersection with }$%
S(u)$\textsl{\ is contained in a }$C^{1}$\textsl{\ smooth curve.}

\bigskip

Proof. Compute the gradient of $F_{a,b}:$

\begin{align}
    \tag{3.1} \nabla
F_{a,b}&=(au_{xx}+b(u_{xy}+1),a(u_{xy}-1)+bu_{yy}) \\
&=(a,b)\left[
\begin{array}{cc}
u_{xx} & u_{xy}-1 \\
u_{xy}+1 & u_{yy}%
\end{array}%
\right] =(a,b)U.\nonumber
\end{align}

\noindent Note that $U$ is never a zero matrix since $u_{xy}+1$ and $u_{xy}-1$ can
never be zero simultaneously. So there exists at most one unit eigenvector $%
(a_{0},b_{0})$ of eigenvalue 0 up to a sign for $U$ at $p_{0}.$ For $%
(a,b)\neq \pm (a_{0},b_{0}),$ $\nabla F_{a,b}\neq 0$ at $p_{0}.$ Then by the
implicit function theorem, there exists a small neighborhood of $p_{0},$ in
which $\Gamma _{a,b}$ at least for $(a,b)\neq \pm (a_{0},b_{0})$ forms a $%
C^{1}$ smooth curve. On the other hand, it is obvious that $S(u)$\textsl{\ }%
is contained\textsl{\ }in $\Gamma _{a,b}$. We are done.
\begin{flushright}
Q.E.D.
\end{flushright}

\bigskip

{\bf Lemma 3.2.} \textsl{Let }$u\in C^{2}(\Omega )$ and suppose\textsl{\ }$p_{0}$%
\textsl{\ }$\in $\textsl{\ }$S(u)\subset \Omega $\textsl{\ is not isolated,
i.e., there exists a sequence of distinct points }$p_{j}\in $\textsl{\ }$%
S(u) $\textsl{\ approaching }$p_{0}$. Then\textsl{\ }$\det U(p_{0})$\textsl{%
\ is zero. }

\bigskip

Proof. Let $\Gamma _{a,b}$ be the $C^{1}$ smooth curve as in the proof of
Lemma 3.1. Since $\Gamma _{a,b}$ is $C^{1}$ smooth near $p_{0},$ we can take
a parameter $s$ of unit speed for $\Gamma _{a,b}$ near $p_{0}$, and find a
subsequence of $p_{j},$ still denoted as $p_{j}$, $%
p_{j}=(x,y)(s_{j}),p_{0}=(x,y)(s_{0})$ such that $s_{j}$ tends to $s_{0}$
monotonically. Since $(u_{x}-y)(p_{j})=(u_{x}-y)(x,y)(s_{j})=0$ for all $j,$
there exists $\tilde{s}_{j}$ between $s_{j}$ and $s_{j+1}$ such that $%
d(u_{x}-y)/ds=0$ at $\tilde{s}_{j}.$ So by the chain rule we obtain

\begin{align}
    \tag{3.2a} u_{xx}\frac{dx}{ds}+(u_{xy}-1)\frac{dy}{ds}=0
\end{align}

\noindent at $\tilde{s}_{j}$ and at $s_{0}$ by letting $\tilde{s}_{j}$ go to $s_{0}.$
Starting from $(u_{y}+x)(p_{j})=(u_{y}+x)(x,y)(s_{j})=0,$ we also obtain by
a similar argument that

\begin{align}
    \tag{3.2b} (u_{xy}+1)\frac{dx}{ds}+u_{yy}\frac{dy}{ds}=0
\end{align}

\noindent at $s_{0}.$ But $(\frac{dx}{ds},\frac{dy}{ds})\neq 0$ since $s$ is a
unit-speed parameter. It follows from $(3.2a)$, $(3.2b)$ that $\det
U(p_{0})=0.$
\begin{flushright}
Q.E.D.
\end{flushright}

Note that we do not assume any condition on $H$ in Lemma 3.1 and
Lemma 3.2. If we make a restriction on $H,$ we can obtain the converse of
Lemma 3.2.

\bigskip

{\bf Theorem 3.3.} \textsl{Let }$\Omega $\textsl{\ be a domain in
the }$xy-$\textsl{plane. Let }$u\in C^{2}(\Omega )$\textsl{\ be such that }$div$%
$N(u)=H$\textsl{\ in }$\Omega \backslash S(u)$\textsl{. Suppose }$%
|H|\leq C\frac{1}{r}$\textsl{\ near a singular point }$p_{0}\in S(u)$\textsl{%
\ where }$r(p)=|p-p_{0}|$\textsl{\ for }$p\in \Omega $\textsl{\ and }$C$%
\textsl{\ is a positive constant}$.$\textsl{\ Then the following are
equivalent:}

\textsl{\ \ \ \ \ \ (1) }$p_{0}$\textsl{\ is not isolated in }$S(u),$

\textsl{\ \ \ \ \ \ (2) }$\det U(p_{0})=0,$

\textsl{\ \ \ \ \ \ (3) }\textsl{there exists a small neighborhood of }$p_{0}$\textsl{%
\ which intersects with }$S(u)$\textsl{\ in exactly a }$C^{1}$\textsl{\
smooth curve past }$p_{0}.$

\bigskip

Proof.\ (1)$\Rightarrow $(2) by Lemma 3.2. (3)$\Rightarrow $(1) is obvious.
It suffices to show that (2)$\Rightarrow $(3). Suppose $\det U(p_{0})=0$%
. Note that $U(p_{0})\neq 0-$matrix since the off-diagonal terms
$u_{xy}+1$ and $u_{xy}-1$ in $(U)$ can never be zero
simultaneously. Therefore $rank(U(p_{0}))=1$, and there exists a
unique $(a_{0},b_{0}),$ up to sign, such that
$(a_{0},b_{0})U=0$ at $p_{0}.$ So for any $(a,b)\neq \pm (a_{0},b_{0}),$ $%
F_{a,b}\equiv a(u_{x}-y)+b(u_{y}+x)$ satisfies $\nabla F_{a,b}\neq 0$ at $%
p_{0}$ and for $(a_{1},b_{1})\neq \pm (a_{0},b_{0}),(a_{2},b_{2})\neq \pm (a_{0},b_{0})$ and
$(a_{1},b_{1})\neq \pm (a_{2},b_{2})$

\begin{align}
\tag{3.3}\nabla F_{a_{1},b_{1}}=c\nabla F_{a_{2},b_{2}}
\end{align}


\noindent at $p_{0}$ where $c$ is a nonzero constant. Therefore
$\Gamma _{a_{1},b_{1}}$ and $\Gamma _{a_{2},b_{2}}$ are $C^1$
smooth curves in a neighborhood of $p_{0}$ (recall that $\Gamma
_{a,b}$ is defined by $F_{a,b}=0$). Also $\Gamma _{a_{1},b_{1}}$
and $\Gamma _{a_{2},b_{2}}$ are tangent at $p_{0}$. Thus we can
take unit-speed arc-length parameters $s$, $t$ for $\Gamma
_{a_{1},b_{1}}$, $\Gamma _{a_{2},b_{2}}$ described by
$\gamma_{1}(s)$, $\gamma_{2}(t)$, respectively, so that
$\gamma_{1}(0)=\gamma_{2}(0)=p_{0}$ and
$\gamma_{1}^{\prime}(0)=\gamma_{2}^{\prime}(0)$. Observe that

\begin{align}
(r{\circ}\gamma_{1})^{\prime}(0+)=(r{\circ}\gamma_{2})^{\prime}(0+)=1 \nonumber\\
(r{\circ}\gamma_{1})^{\prime}(0-)=(r{\circ}\gamma_{2})^{\prime}(0-)=-1. \nonumber
\end{align}

\noindent Therefore we can find a small $\epsilon >0$ such that
$(r{\circ}\gamma_{1})^{\prime}(s),(r{\circ}\gamma_{2})^{\prime}(t)>0$
for $s,t>0$
($(r{\circ}\gamma_{1})^{\prime}(s),(r{\circ}\gamma_{2})^{\prime}(t)<0$
for $s,t<0$, respectively) whenever $\gamma_{1}(s),\gamma_{2}(t)$
$\in$ $B_{\epsilon }(p_{0})$, a ball of radius $\epsilon$ and
center $p_{0}$. Also note that $S(u)\cap {B_{\epsilon
}(p_{0})}\subset \Gamma _{a_{i},b_{i}},i=1,2$. We can write

\begin{align}
(\Gamma _{a_{1},b_{1}}\cap B_{\epsilon }(p_{0}))\backslash
S(u)=\cup _{j=1}^{\infty }\gamma_{1}(]s_{j}^{+},{\tilde{s}}_{j}^{+}[)\cup
\cup_{j=1}^{\infty }\gamma_{1}(]s_{j}^{-},{\tilde{s}}_{j}^{-}[)\nonumber\\
(\Gamma _{a_{2},b_{2}}\cap B_{\epsilon }(p_{0}))\backslash
S(u)=\cup _{j=1}^{\infty }\gamma_{2}(]t_{j}^{+},{\tilde{t}}_{j}^{+}[)\cup
\cup_{j=1}^{\infty }\gamma_{2}(]t_{j}^{-},{\tilde{t}}_{j}^{-}[).\nonumber
\end{align}

\noindent where
$]s_{j}^{+},{\tilde{s}}_{j}^{+}[,]s_{j}^{-},{\tilde{s}}_{j}^{-}[$,
etc. denote open intervals and we have arranged

\begin{align}
s_{1}^{-}<{\tilde{s}}_{1}^{-}{\leq}..s_{j}^{-}<{\tilde{s}}_{j}^{-}{\leq}
s_{j+1}^{-}<{\tilde{s}}_{j+1}^{-}{\leq}..{\leq}0{\leq}..{\leq}
s_{j+1}^{+}<{\tilde{s}}_{j+1}^{+}{\leq}s_{j}^{+}<{\tilde{s}}_{j}^{+}{\leq}..{\leq}
s_{1}^{+}<{\tilde{s}}_{1}^{+}\nonumber\\
t_{1}^{-}<{\tilde{t}}_{1}^{-}{\leq}..t_{j}^{-}<{\tilde{t}}_{j}^{-}{\leq}
t_{j+1}^{-}<{\tilde{t}}_{j+1}^{-}{\leq}..{\leq}0{\leq}..{\leq}
t_{j+1}^{+}<{\tilde{t}}_{j+1}^{+}{\leq}t_{j}^{+}<{\tilde{t}}_{j}^{+}
{\leq}..{\leq}t_{1}^{+}<{\tilde{t}}_{1}^{+}\nonumber
\end{align}

\noindent (if there are only finite number of open intervals for
negative (positive, respectively) $s$, then
${\tilde{s}}_{j}^{-}=s_{j}^{-}={\tilde{s}}_{m}^{-}$
($s_{j}^{+}={\tilde{s}}_{j}^{+}=s_{m}^{+}$, respectively) for
$j{\geq}m$, a certain integer, by convention. Note that
${\tilde{s}}_{j}^{-}$ or $s_{j}^{+}$ may equal $0$ for some finite
$j$. In this case, all the $s$ with subindex larger than $j$ equal
$0$. We apply the same convention to the parameter $t$). Since
$r{\circ}\gamma_{1}$ and $r{\circ}\gamma_{2}$ are monotonic
(increasing for positive parameters and decreasing for negative
parameters), we actually have
$\gamma_{1}(s_{j}^{\pm})=\gamma_{2}(t_{j}^{\pm})$,
$\gamma_{1}({\tilde{s}}_{j}^{\pm})=\gamma_{2}({\tilde{t}}_{j}^{\pm})$
except possibly $\gamma_{1}(s_{1}^{-})\neq\gamma_{2}(t_{1}^{-})$
or $\gamma_{1}({\tilde{s}}_{1}^{+})\neq
\gamma_{2}({\tilde{t}}_{1}^{+})$. Let
$e_{j}=\gamma_{1}(s_{j}^{+})$ and
$\tilde{e}_{j}=\gamma_{1}({\tilde{s}}_{j}^{+})$. Then
either $\Gamma _{a_{1},b_{1}}$ and $\Gamma _{a_{2},b_{2}}$ meet at
$\tilde{e}_{1}$ or $\tilde{e}_{2}$, or they do not meet in
$B_{\epsilon }(p_{0}){\backslash}\{ p_{0}\}$ for positive
parameters. In the former situation, we need to show that $e_{i}$,
$i\geq 1$, does not converge to $p_{0}$ as $i\rightarrow \infty $
(then there is a smaller ball $B_{\epsilon ^{\prime }}(p_{0})$
such that $S(u)\cap B_{\epsilon ^{\prime }}(p_{0})$ contains
$\gamma_{1}([0,{\bar s}^{+}[)$ for a small ${\bar s}^{+}>0$).
Suppose it is not so. Let $\Omega _{i}$ be the region surrounded
by $\Gamma _{a_{1},b_{1}}$ and $\Gamma _{a_{2},b_{2}}$ from
$e_{i}$ to $\tilde{e}_{i}$ for all $i\geq 1$ or $2$ (if
$\gamma_{1}({\tilde{s}}_{1}^{+})\neq
\gamma_{2}({\tilde{t}}_{1}^{+})$) (note that the curves $\Gamma
_{a_{1},b_{1}}$ and $\Gamma _{a_{2},b_{2}}$ only meet at singular
points since $(a_{1},b_{1})\neq \pm (a_{2},b_{2})$, i.e., points
of the arcs
$\gamma_{1}([\tilde{s}_{i+1}^{+},s_{i}^{+}])=\gamma_{2}([\tilde{t}_{i+1}^{+},t_{i}^{+}]).$)
Observe that $\Gamma _{a_{1},b_{1}}$ and $\Gamma _{a_{2},b_{2}}$
asymptotically approximate the same straight line by $(3.3).$ So
the distance function $r(p)\equiv |p-p_{0}|$ is one-to-one for
$p\in $ $\Gamma _{a_{1},b_{1}}(\Gamma _{a_{2},b_{2}},$
respectively) near $p_{0}$ with parameter $s>0$ ($t>0$,
respectively) since
$(r{\circ}\gamma_{1})^{\prime}(s),(r{\circ}\gamma_{2})^{\prime}(t)>0$
for $s,t>0$ as shown previously. Now we want to compare both sides
of

\begin{align}
    \tag{3.4} \oint_{\partial \Omega _{i}}N(u)\cdot \nu ds=\int
\int_{\Omega _{i}}divN(u)dxdy=\int \int_{\Omega _{i}}Hdxdy
\end{align}

\noindent where $\nu $ is the unit outward normal to $\Gamma _{a_{1},b_{1}}$ and $%
\Gamma _{a_{2},b_{2}}.$ On $\Gamma _{a_{1},b_{1}}(\Gamma _{a_{2},b_{2}},$
respectively), $N(u)\perp (a_{1},b_{1})$ ($(a_{2},b_{2}),$ respectively). So $%
N(u)$ is a constant unit vector field along $\Gamma _{a_{1},b_{1}}$ ($\Gamma
_{a_{2},b_{2}},$ respectively). On the other hand, $\nu $ approaches a fixed
unit vector $\nu _{0}=\nabla F_{a_{1},b_{1}}(p_{0})/|\nabla
F_{a_{1},b_{1}}(p_{0})|=\nabla F_{a_{2},b_{2}}(p_{0})/|\nabla
F_{a_{2},b_{2}}(p_{0})|$ as $e_{i},\tilde{e}_{i}$ tend to $p_{0}$ for
$i$ large. It follows that $N(u)\cdot \nu $ tends to a constant $c_{1}=\pm
(-b_{1},a_{1})\cdot \nu _{0}$ ($c_{2}=\pm (-b_{2},a_{2})\cdot \nu _{0}$,
respectively) along $\Gamma _{a_{1},b_{1}}(\Gamma _{a_{2},b_{2}},$
respectively) as $i$ goes to infinity. We choose (in advance) $(a_{1},b_{1})\neq \pm
(a_{2},b_{2})$ (also both $\neq \pm (a_{0},b_{0})$) such that $c_{1}\neq
c_{2}.$ Thus we can estimate

\begin{align}
    \tag{3.5} \mid \oint_{\partial \Omega _{i}}N(u)\cdot \nu ds\mid \geq
 (\mid c_{1}-c_{2}\mid -\delta _{i})\mid r(\tilde{e}_{i})-r(e_{i})\mid
\end{align}

\noindent for some small positive $\delta _{i}$ that goes to $0$ as $i\rightarrow
\infty .$ On the other hand, we can make $\Omega _{i}$ contained in a
fan-shaped region of angle $\theta _{i}$ with vertex $p_{0}$ so that $\theta
_{i}\rightarrow 0$ as $i\rightarrow \infty .$ We estimate

\begin{align}
    \tag{3.6} \mid \int \int_{\Omega _{i}}Hdxdy\mid \leq \mid
\int_{r(e_{i})}^{r(\tilde{e}_{i})}(\int_{0}^{\theta _{i}}\frac{C}{r}d\theta
)rdr\mid \leq C\theta _{i}\mid r(\tilde{e}_{i})-r(e_{i})\mid .
\end{align}

\noindent Substituting $(3.5)$ and $(3.6)$ into $(3.4),$ we obtain

$$(\mid c_{1}-c_{2}\mid -\delta _{i})\mid r(%
\tilde{e}_{i})-r(e_{i})\mid \leq C\theta _{i}\mid r(\tilde{e}%
_{i})-r(e_{i})\mid .$$

\noindent Hence $\mid c_{1}-c_{2}\mid -\delta _{i}\leq C\theta _{i}.$ But $\delta _{i}$
and $\theta _{i}$ tending to $0$ gives $c_{1}=c_{2},$ a contradiction.
Another situation is that $\Gamma _{a_{1},b_{1}}$ and $\Gamma _{a_{2},b_{2}}$
do not meet in a small neighborhood of $p_{0}$ except at $p_{0}$ for $s,t>0$. In this case, let $\Omega
_{i}$ be the region surrounded by $\Gamma _{a_{1},b_{1}},$ $\Gamma
_{a_{2},b_{2}}$, and $\partial B_{r_{i}}(p_{0})$ for large $i$ and contained
in a fan-shaped region of angle $\theta _{i}$ with vertex $p_{0}$ so that $%
r_{i}\rightarrow 0,$ $\theta _{i}\rightarrow 0$ as $i\rightarrow
\infty .$ Observing that the arc length of $\bar{\Omega}_{i}$
$\cap $ $\partial B_{r_{i}}(p_{0})$ is bounded by $\theta
_{i}r_{i},$ we can reach a contradiction by a similar argument as
above. For $s,t<0$, we apply a similar argument to conclude that
there is a small $\epsilon ^{\prime \prime }>0$ such that
$S(u)\cap B_{\epsilon ^{\prime \prime }}(p_{0})$ contains
$\gamma_{1}(]{\bar s}^{-},0])$ for a small ${\bar s}^{-}<0$. Now
an even smaller ball of radius $<min(\epsilon ^{\prime},\epsilon
^{\prime \prime })$ and center $p_{0}$ will serve our purpose.
\begin{flushright}
Q.E.D.
\end{flushright}

\noindent $\mathbf{Proof\ of\ Theorem\ B:}$

It is an immediate consequence of Theorem 3.3.
\begin{flushright}
Q.E.D.
\end{flushright}

We remark that Theorem B does not hold if we remove the condition on
$H$.

\medskip

{\bf Example 1.} Let $u=xg(y).$ Then $u_{x}=g(y),u_{y}=xg^{\prime }(y).$
It follows that the singular set $S(u)=\{g(y)=y$ and $g^{\prime
}(y)+1=0\}\cup \{g(y)=y$ and $x=0\}.$ Take $g(y)=\exp (-\frac{1}{y^{2}})\sin
(-\frac{1}{y})+y.$ Compute $g^{\prime }(y)$ $=$ $2\exp (-\frac{1}{y^{2}}%
)y^{-3}\sin (-\frac{1}{y})$ $+$ $\exp (-\frac{1}{y^{2}})y^{-2}\cos (-\frac{1%
}{y})$ $+$ $1.$ So for $y$ small, $g^{\prime }(y)+1=0$ has no solution.
Therefore $S(u)$ (when $y$ is small) $=$ $\{g(y)=y$, $x=0\}$ $=$ $\{\sin (-%
\frac{1}{y})=0$, $x=0\}$ has infinitely-many points near $(0,0).$ Note that $%
g$, hence $u$, is $C^{\infty }$ smooth. This example shows that even for $%
u\in C^{\infty }$, $S(u)$ may contain non-isolated points which do
not belong to curves in $S(u)$. On the other hand, the p-mean
curvature $H(u)$ has a blow-up rate $\exp (\frac{1}{r^{2}})$ for a
certain sequence of points $(x_{j},y_{j})$ satisfying $x_{j}$ $=$
$\exp (-\frac{1}{y_{j}^{2}})$ and converging to $(0,0).$

\medskip

{\bf Example 2.} Let $u=\pm \frac{1}{2}(x^{2}+y^{2}).$ Then $u_{x}=\pm x,$
$u_{y}=\pm y.$ So $S(u)=\{(0,0)\}.$ By $(2.10)$ we compute the p-mean
curvature $H=\pm 2^{-1/2}r^{-1}$ where $r=\sqrt{x^{2}+y^{2}}.$ This is the
case that the equality of the condition on $H$ holds.

\medskip

{\bf Example 3.} The following example shows that in Theorem 3.3, (2) does
not imply (1) if we remove the condition on $H.$ Let

$$u=\frac{1}{2}(x^{2}-y^{2})+\frac{1}{\beta }%
(sgn(x)|x|^{\beta }-sgn(y)|y|^{\beta })$$

\noindent for $\beta >2.$ Compute $u_{x}$ $=$ $x+|x|^{\beta -1},$ $u_{y}$ $=$ $%
-y-|y|^{\beta -1}.$ It follows that $(u_{x}-y)-(u_{y}+x)$ $=$ $|x|^{\beta
-1}+|y|^{\beta -1}.$ So $S(u)$ $\subset $ $\{|x|^{\beta -1}+|y|^{\beta
-1}=0\}$, and hence $S(u)$ $=$ $\{(0,0)\}.$ This means that $(0,0)$ is an
isolated singular point. Compute $u_{xx}$ $=$ $1+(\beta -1)sgn(x)|x|^{\beta
-2},$ $u_{xy}=0$ and $u_{yy}$ $=$ $-1-(\beta -1)sgn(y)|y|^{\beta -2}.$ It is
easy to see that $|u_{xx}|\leq \beta $ and $|u_{yy}|\leq \beta $ for $%
|x|\leq 1,$ $|y|\leq 1.$ So by $(2.10)$ we can estimate

$$H\leq \frac{|u_{xx}|\text{ }+\text{ }|u_{yy}|}{D}\leq \frac{\beta
\text{ }+\text{ }\beta }{(|x|^{\beta -1}+|y|^{\beta -1})/\sqrt{2}}\leq C%
\frac{1}{r^{\beta -1}}$$

\noindent near $(0,0)$ where $r=\sqrt{x^{2}+y^{2}}.$ In the second inequality above,
we have used the following estimate

\begin{align}
D^{2}\equiv (u_{x}-y)^{2}+(u_{y}+x)^{2}&=(x+|x|^{\beta
-1}-y)^{2}+(-y-|y|^{\beta -1}+x)^{2}\nonumber\\
&\geq \frac{1}{2}(|x|^{\beta -1}+|y|^{\beta
-1})^{2}\nonumber
\end{align}

\noindent (by $2(a^{2}+b^{2})$ $\geq $ $(a-b)^{2}$). On the other hand, we observe
that $U=\left(
\begin{array}{cc}
1 & -1 \\
1 & -1%
\end{array}%
\right) $ at $(0,0).$ It follows that $\det U$ $=$ $0$ at $(0,0).$

\medskip

According to Theorem B, $S(u)$ may contain some $C^{1}$ smooth
curves, called singular curves. We will study the behavior of
$N(u)$ near a point of a singular curve. First we show that for a
graph $t = u(x,y)$ the p-minimal graph equation $(pMGE)$ is
rotationally invariant.

\bigskip

{\bf Lemma 3.4.} \textsl{Suppose }$u\in C^{2}.$\textsl{\ Let }$x=a\tilde{x}-b\tilde{y},%
y=b\tilde{x}+a\tilde{y}$\textsl{\ where }$a^{2}+b^{2}=1$\textsl{, and let }$%
\tilde{u}(\tilde{x},\tilde{y})=u(x(\tilde{x},\tilde{y}),y(\tilde{x},\tilde{y})).$\textsl{\ Then }$\widetilde{div}N(%
\tilde{u})=divN(u)$\textsl{\ where }$\widetilde{div}$\textsl{\ denotes the
plane divergence with respect to }$(\tilde{x},\tilde{y}).$

\bigskip

Proof. First we observe that both $\nabla u$ (viewed as a row vector) and $%
(-y,x)$ satisfy the following transformation law (for a plane {\em
vector}):

$$\tilde{\nabla}\tilde{u}=(\nabla u)\left[
\begin{array}{cc}
a & -b \\
b & a%
\end{array}%
\right] , \ \ (-\tilde{y},\tilde{x})=(-y,x)\left[
\begin{array}{cc}
a & -b \\
b & a%
\end{array}%
\right] .$$

\noindent It follows that $N(u)=[\nabla u+(-y,x)]/D$ also obeys the same
transformation law by noting that $D$ is invariant. Then a direct
computation shows that $\widetilde{div}N(\tilde{u})=divN(u)$.
\begin{flushright}
Q.E.D.
\end{flushright}

\bigskip

Let $\Gamma _{s}$ be a singular curve contained in $S(u)$ for a $%
C^{2}$ smooth $u$ (defined on a certain domain). Let $p_{0}\in
\Gamma _{s}. $ Suppose there exists a ball $B_{\epsilon }(p_{0})$
such that $\Gamma _{s}\cap B_{\epsilon }(p_{0})$ divides
$B_{\epsilon }(p_{0})$ into two disjoint nonsingular parts (this
is true if $|H|\leq C\frac{1}{r}$ near $p_{0}$ in view of Theorem
B). That is to say, $B_{\epsilon }(p_{0})\backslash (\Gamma
_{s}\cap B_{\epsilon }(p_{0}))=B^{+}\cup B^{-}$ where $B^{+}$ and
$B^{-}$ are disjoint domains (proper open and connected).

\bigskip

{\bf Proposition 3.5.} \textsl{Suppose we have the situation as described above.
Then both }$N(u)(p_{0}^{+})$ $\equiv $ $\lim_{p\in B^{+}\rightarrow
p_{0}}N(u)(p)$\textsl{\ and }$N(u)(p_{0}^{-})$ $\equiv $ $\lim_{p\in
B^{-}\rightarrow p_{0}}N(u)(p)$\textsl{\ exist. Moreover, }$%
N(u)(p_{0}^{+})=-N(u)(p_{0}^{-}).$

\bigskip

Proof. By Lemma 3.4 we may assume the $x$-axis past $p_{0}=(x_{0},y_{0})$ is
transverse to $\Gamma _{s}$. Moreover, we may assume either $%
u_{xx}(x_{0},y_{0})\neq 0$ or $(u_{xy}+1)(x_{0},y_{0})\neq 0$ (Note that in $%
(3.2a)$ $\frac{dy}{ds}\neq 0$ at $p_{0}.$ So if $u_{xx}(x_{0},y_{0})=0,$
then $(u_{xy}-1)(x_{0},y_{0})=0.$ It follows that $(u_{xy}+1)(x_{0},y_{0})=2%
\neq 0$). Let $\Gamma _{s}\cap \{y=c\}=\{(x_{c},c)\}$ for $c$ close to $%
y_{0}.$ Since $(x_{c},c)\in \Gamma _{s}$, we have $(u_{y}+x)(x_{c},c)=0$ and
$(u_{x}-y)(x_{c},c)=0.$ So if $u_{xx}(x_{0},y_{0})\neq 0,$ we compute

\begin{align}
    \tag{3.7} \frac{(u_{y}+x)(x,c)}{(u_{x}-y)(x,c)}&=\frac{%
(u_{y}+x)(x,c)-(u_{y}+x)(x_{c},c)}{(u_{x}-y)(x,c)-(u_{x}-y)(x_{c},c)} \\
&=\frac{%
(x-x_{c})(u_{xy}+1)(x_{c}^{1},c)}{(x-x_{c})u_{xx}(x_{c}^{2},c)} \nonumber \\
&\mbox{(for }x_{c}^{1},x_{c}^{2} \mbox{ between } x_{c} \mbox{ and }
x \mbox{ by the mean value theorem)} \nonumber \\
&=\frac{
(u_{xy}+1)(x_{c}^{1},c)}{u_{xx}(x_{c}^{2},c)}. \nonumber
\end{align}

\noindent Letting $(x,c)$ go to $(x_{0},y_{0})$ in $(3.7),$ we obtain

\begin{align}
    \tag{3.8} \lim _{p\in B^{+}\rightarrow p_{0}}\frac{u_{y}+x}{%
u_{x}-y}=\frac{u_{xy}+1}{u_{xx}}(p_{0})=\lim_{p\in B^{-}\rightarrow p_{0}}%
\frac{u_{y}+x}{u_{x}-y}.
\end{align}

\noindent Therefore by $(3.8)$ two limits of the unit vector $%
N(u)=(u_{x}-y,u_{y}+x)D^{-1}$ from both sides exist, and their values can
only be different by a sign. That is to say,

\begin{align}
    \tag{3.9}  N(u)(p_{0}^{+})=\pm N(u)(p_{0}^{-}).
\end{align}

Now we observe that $%
(u_{x}-y)(x,y_{0})-(u_{x}-y)(x_{0},y_{0})=u_{xx}(\eta ,y_{0})(x-x_{0})$ for
some $\eta $ between $x_{0}$ and $x.$ Since $(u_{x}-y)(x_{0},y_{0})=0,$ $%
(u_{x}-y)(x,y_{0})$ and $u_{xx}(x_{0},y_{0})$ have the same (different,
respectively) sign for $x>x_{0}$ ($x<x_{0}$, respectively) and $x$ being
close enough to $x_{0}.$ Thus we should have the negative sign in $(3.9).$
We are done. In case $u_{xx}(x_{0},y_{0})=0,$ we have $%
(u_{xy}+1)(x_{0},y_{0})\neq 0.$ So we compute $\frac{(u_{x}-y)(x,c)}{%
(u_{y}+x)(x,c)}$ instead of $\frac{(u_{y}+x)(x,c)}{(u_{x}-y)(x,c)},$ and get
$\frac{u_{xx}}{u_{xy}+1}(p_{0})$ instead of $\frac{u_{xy}+1}{u_{xx}}(p_{0})$
in $(3.8).$ Then we still conclude $(3.9).$ Instead of $(u_{x}-y)(x,y_{0}),$
we consider $(u_{y}+x)(x,y_{0}).$ A similar argument as above shows that $%
(u_{y}+x)(x,y_{0})$ will have the same (different, respectively) sign as $%
(u_{xy}+1)(x_{0},y_{0})$ for $x>x_{0}$ ($x<x_{0}$, respectively) and $x$
being close enough to $x_{0}.$ So we still take the negative sign in $(3.9).$
\begin{flushright}
Q.E.D.
\end{flushright}

\bigskip

Note that we do not assume any condition on $H$ in Proposition 3.5.
We will study how two characteristic curves meet at a point of a singular
curve. We say a characteristic curve $\Gamma \subset B^{+}$ or $B^{-}$
touches $\Gamma _{s}$ at $p_{0}$ if $p_{0}\in \bar{\Gamma},$ the closure of $%
\Gamma $ in the $xy$-plane, and touches transversally if, furthermore, $p_{0}$
is the only intersection point of the tangent line of $\Gamma _{s}$ at $%
p_{0} $ and the tangent line of $\Gamma $ at $p_{0}$ (which makes
sense in view of Proposition 3.5). 

\bigskip

{\bf Corollary 3.6.} \textsl{Suppose }\textsl{we have the same
assumptions as in Proposition 3.5.} \textsl{Then there is a unique
characteristic curve }$\Gamma _{+}\subset B^{+}$\textsl{\ touching
}$\Gamma _{s}$\textsl{\ at }$p_{0}$\textsl{\ transversally (with
}$N^{\perp}(u)(p_{0}^{+})$\textsl{\ being the tangent vector of
}$\bar{\Gamma}_{+}$\textsl{\ at }$p_{0}$).\textsl{\ Similarly
there exists a unique characteristic curve }$\Gamma _{-}\subset
B^{-}$\textsl{\ touching }$\Gamma _{s}$\textsl{\ also at
}$p_{0}$\textsl{\
so that }$\Gamma _{-}\cup \{p_{0}\}\cup \Gamma _{+}$\textsl{\ forms a }$%
C^{1} $\textsl{\ smooth curve in }$B_{\epsilon }(p_{0}).$

\bigskip

First note that we can change the sign of $N^{\perp }(u)$ if
necessary to make a $C^{0}$ characteristic (i.e. tangent to
integral curves
of $N^{\perp }(u)$ where $N^{\perp }(u)$ is defined) vector field $\check{N}%
^{\perp }(u)$ in $B_{\epsilon }(p_{0})$ in view of Proposition
3.5. So $\bar{\Gamma}_{+}$ of any characteristic curve $\Gamma
_{+}\subset B^{+}$ touching $\Gamma _{s}$ at $p_{0}$ must have the
tangent vector $N^{\perp}(u)(p_{0}^{+})$ at $p_{0}$. To show the
uniqueness of $\Gamma _{+}$ and the transversality of
$N^{\perp}(u)(p_{0}^{+})$ to $\Gamma _{s}$ in Corollary 3.6, we
observe that

\bigskip

{\bf Lemma 3.7.} \textsl{Suppose }\textsl{we have the same assumptions as in Proposition
3.5.} \textsl{\ Then there hold}

\textsl{\ \ \ \ \ \ (1) }$N^{\perp}(u)(p_{0}^{+})U(p_{0})=0$,

\textsl{\ \ \ \ \ \ (2) }$(c,d)U^{T}(p_{0})=0$

\noindent \textsl{for a nonzero vector }$(c,d)$\textsl{ tangent to }$\Gamma_{s}$
\textsl{ at }$p_{0}$\textsl{, where we view }$N^{\perp}(u)(p_{0}^{+})$
\textsl{ as a row vector and }$U^{T}$\textsl{ denotes the transpose of }$U.$

\bigskip

Proof. Let $(x(s),y(s))$ describe $\Gamma_{s}$ so that
$p_{0}=(x(0),y(0))$ and $(x'(0),y'(0))=(c,d)$. Since $u_{x}-y=0$,
$u_{y}+x=0$ on $\Gamma_{s}$, we differentiate these two equations
to get (3.2a) and (3.2b) along $\Gamma_{s}$ by the chain rule. At
$s=0$, we obtain (2). From the proof of Proposition 3.5, we learn
that $N(u)(p_{0}^{+})$ is proportional to
$(u_{xx},u_{xy}+1)(p_{0})$ (if this is not a zero vector). A
similar argument (L'H$\hat{o}$pital's rule in the $y$-direction)
shows that $N(u)(p_{0}^{+})$ is also proportional to
$(u_{xy}-1,u_{yy})(p_{0})$ (if this is not a zero vector).
Therefore $N^{\perp}(u)(p_{0}^{+})$ is perpendicular to
$(u_{xx},u_{xy}+1)(p_{0})$ and $(u_{xy}-1,u_{yy})(p_{0})$, hence
(1) follows.
\begin{flushright}
Q.E.D.
\end{flushright}

We will give a proof of Corollary 3.6 after the proof of Theorem 3.10.



\medskip

Remark. If $u$ is not of class $C^{2},$ the extension theorem
(Corollary 3.6) may fail as the following example shows. Consider
the function

\[u(x,y)=\left\{%
\begin{array}{ll}
-xy, &y\geq 0 \\
-xy+y^{2}\cot \vartheta ,&y<0%
\end{array} \right.
\]

\noindent where $0<\vartheta <\frac{\pi }{2}.$ There holds

\[N^{\perp }(u)=\left \{%
\begin{array}{ll}
(0,1), & y>0 \\
(\cos \vartheta ,\sin \vartheta ),& y<0. %
\end{array} \right.%
\]

\noindent Note that the function $u$ is of class $C^{1,1}$ on
$R^{2}$ and satisfies $divN(u)=0$ in $R^{2}\backslash \{y=0\}$
where $u$ is of class $C^{2}.$

\medskip

Next we want to analyze the configuration of characteristic curves
near an isolated singularity. First observe that for a $C^{2}$ smooth $u$
defined on a domain $\Omega ,$ characteristic curves are also the integral
curves of the $C^{1}$ smooth vector field $N^{\perp }D$ $=$ $%
(u_{y}+x,-u_{x}+y)$ which vanishes at singular points. We think of $N^{\perp
}D$ as a mapping: $\Omega \subset R^{2}\rightarrow R^{2},$ so that the
differential $d(N^{\perp }D)_{p_{0}}$: $R^{2}\rightarrow R^{2}$ is defined
for $p_{0}\in \Omega .$

\bigskip

{\bf Lemma 3.8.} \textsl{Let }$u\in C^{2}(\Omega ).$\textsl{\ Suppose }$|H|=o(%
\frac{1}{r})$\textsl{\ (little ''}$o")$\textsl{\ near an isolated singular
point }$p_{0}\in \Omega $\textsl{\ where }$r(p)=|p-p_{0}|.$\textsl{\ Then }$%
d(N^{\perp }D)_{p_{0}}$\textsl{\ is the identity linear transformation and
the index of }$N^{\perp }D$\textsl{\ at the isolated zero }$p_{0}$\textsl{\
is }$+1.$\textsl{ Moreover,} $u_{xx}=u_{xy}=u_{yy}=0$\textsl{ at }$p_0$.

\bigskip

Proof. In view of (2.10), we write $H=D^{-3}P(u)$ where

\begin{align}
\tag{3.10}
P(u)\equiv
(u_{y}+x)^{2}u_{xx}-2(u_{y}+x)(u_{x}-y)u_{xy}+(u_{x}-y)^{2}u_{yy}.
\end{align}

\noindent Let $p-p_{0}=(\bigtriangleup x,\bigtriangleup y),$ $r=r(p)=|p-p_{0}|$ $=$ $%
[(\bigtriangleup x)^{2}+(\bigtriangleup y)^{2}]^{1/2}.$ Noting that $%
(u_{y}+x)(p_{0})$ $=$ $(u_{x}-y)(p_{0})$ $=$ $0$ by the definition of a
singular point$,$ we can express

\begin{align}
\tag{3.11} (u_{y}+x)(p)=(u_{yx}+1)(p_{0})\bigtriangleup
x+u_{yy}(p_{0})\bigtriangleup y+o(r),\\
\tag{3.12} (u_{x}-y)(p)=u_{xx}(p_{0})\bigtriangleup
x+(u_{xy}-1)(p_{0})\bigtriangleup y+o(r)
\end{align}

\noindent for $p$ near $p_{0}.$ Let $a=(u_{yx}+1)(p_{0}),$ $b=u_{yy}(p_{0}),$ $%
c=u_{xx}(p_{0}),$ and $d=(u_{xy}-1)(p_{0}).$ Substituting (3.11) and (3.12)
in (3.10), we compute the highest order term:

\begin{align}
\tag{3.13}
P(u)&=(a\bigtriangleup x+b\bigtriangleup
y)^{2}c\\
&\ \ -(a\bigtriangleup x+b\bigtriangleup y)(c\bigtriangleup
x+d\bigtriangleup y)(a+d)
+(c\bigtriangleup x+d\bigtriangleup y)^{2}b+o(r^{2})\nonumber\\
&=(bc-ad)\{\left(
\begin{array}{cc}
\bigtriangleup x & \bigtriangleup y%
\end{array}%
\right) \left[
\begin{array}{cc}
c & a \\
d & b%
\end{array}%
\right] \left(
\begin{array}{c}
\bigtriangleup x \\
\bigtriangleup y%
\end{array}%
\right) \}+o(r^{2})\nonumber
\end{align}

On the other hand, substituting (3.11), (3.12) into $D^{3}\equiv \lbrack
(u_{x}-y)^{2}+(u_{y}+x)^{2}]^{3/2},$ we obtain

\begin{align}
\tag{3.14} D^{3}=\mid \left[
\begin{array}{cc}
a & b \\
c & d%
\end{array}%
\right] \left(
\begin{array}{c}
\bigtriangleup x \\
\bigtriangleup y%
\end{array}%
\right) \mid ^{3}+o(r^{3}).
\end{align}

\noindent Note that $bc-ad=\det U(p_{0})\neq 0$ by Theorem 3.3. Letting $%
\bigtriangleup y=0$ and assuming $c\neq 0$, we estimate the highest order
term of $H=D^{-3}P(u):$

$$\frac{(bc-ad)c(\bigtriangleup x)^{2}}{[(a^{2}+c^{2})(%
\bigtriangleup x)^{2}]^{3/2}}=\frac{(bc-ad)c(\bigtriangleup x)^{2}}{%
[a^{2}+c^{2}]^{3/2}(\bigtriangleup x)^{3}}=\frac{(bc-ad)c}{[a^{2}+c^{2}]^{3/2}{\bigtriangleup x}}.$$

\noindent The assumption $|H|=o(\frac{1}{r})$ forces $c=0.$ On the other hand,
letting $\bigtriangleup x=0$ will force $b=0$ by a similar argument. Now we
can write

$$H=\frac{-ad(a+d)\bigtriangleup x\bigtriangleup y}{[a^{2}(%
\bigtriangleup x)^{2}+d^{2}(\bigtriangleup y)^{2}]^{3/2}}+o(\frac{1}{r}).$$

\noindent It follows from the assumption $|H|=o(\frac{1}{r})$ again that $a+d=0$
(note that $ad=ad-bc=-\det U(p_{0})\neq 0$). We have proved that $%
u_{xx}(p_{0})=c=0,$ $u_{yy}(p_{0})=b=0,$ $u_{xy}(p_{0})=\frac{a+d}{2}=0.$
Therefore in matrix form,

$$ d(N^{\perp }D)_{p_{0}}=\left[
\begin{array}{cc}
u_{yx}+1 & u_{yy} \\
-u_{xx} & -u_{xy}+1%
\end{array}%
\right] =\left[
\begin{array}{cc}
1 & 0 \\
0 & 1%
\end{array}%
\right] .$$

\noindent Note that $\det (d(N^{\perp }D)_{p_{0}})=1>0.$ So the index of $N^{\perp }D$
at $p_{0}$ is $+1$ (see, e.g., Lemma 5 in Section 6 in [Mil]).
\begin{flushright}
Q.E.D.
\end{flushright}

\bigskip

{\bf Lemma 3.9.} \textsl{Let }$u\in C^{2}(\Omega ).$\textsl{\ Suppose }$|H|=o(%
\frac{1}{r})$\textsl{\ (little ''}$o")$\textsl{\ near an isolated singular
point }$p_{0}\in \Omega $\textsl{\ where }$r(p)=|p-p_{0}|.$\textsl{\ Then
there exists a small neighborhood }$V\subset \Omega $\textsl{\ of }$p_{0}$%
\textsl{\ such that the characteristic curve past a point in }$V\backslash
\{p_{0}\}$\textsl{\ reaches }$p_{0}$\textsl{\ (towards the }$-N^{\perp }$%
\textsl{\ direction) in finite unit-speed parameter.}

\bigskip

Proof. Write $p-p_{0}=(\bigtriangleup x,\bigtriangleup y)=(r\cos \varphi
,r\sin \varphi )$ in polar coordinates. At $p$, we express

\begin{align}
\tag{3.15} N^{\perp }=\alpha (\cos \varphi ,\sin \varphi )+\beta (-\sin
\varphi ,\cos \varphi )
\end{align}

\noindent where

\begin{align}
\tag{3.16} \alpha = N^{\perp }\cdot (\cos \varphi ,\sin \varphi ), \ \
\beta =N^{\perp }\cdot (-\sin \varphi ,\cos \varphi ).
\end{align}

\noindent Noting that $u_{xx}=u_{xy}=u_{yy}=0$ at $p_{0}$ by Lemma 3.8$%
, $ we obtain that

\begin{align}
\tag{3.17} u_{y}+x=\bigtriangleup x+o(r), \ \ -u_{x}+y=\bigtriangleup
y+o(r)
\end{align}

\noindent near $p_{0}$ by (3.11) and (3.12). It follows that $D=[(\bigtriangleup
x)^{2}+(\bigtriangleup y)^{2}]^{1/2}+o(r)$ $=$ $r+o(r)$ near $p_{0}.$
{From} this and (3.17), we can estimate $N^{\perp }$ $=$ $(\frac{%
\bigtriangleup x}{r}+o(1),$ $\frac{\bigtriangleup y}{r}+o(1)).$ Substituting
this and $(\cos \varphi ,\sin \varphi )$ $=$ $(\frac{\bigtriangleup x}{r},%
\frac{\bigtriangleup y}{r})$ in (3.16) gives

\begin{align}
\tag{3.18}  \alpha = \frac{(\bigtriangleup x)^{2}+(\bigtriangleup
y)^{2}}{r^{2}}+o(1)=1+o(1).
\end{align}

Now let $(x(s),y(s))$ describe a characteristic curve in the $%
-N^{\perp }$ direction, i.e., $\frac{d(x(s),y(s))}{ds}=-N^{\perp }.$ We
compute

\begin{align}
 \frac{d(x(s),y(s))}{ds}&=\frac{d(r(s)\cos \varphi (s),r(s)\sin \varphi
(s))}{ds}\nonumber \\
&=\frac{dr}{ds}(\cos \varphi ,\sin
\varphi )+r\frac{d\varphi }{ds}(-\sin \varphi ,\cos \varphi ).\nonumber
\end{align}

\noindent Comparing with (3.15), we obtain

\begin{align}
\tag{3.19} \frac{dr}{ds}=-\alpha , \ \ r\frac{d\varphi
}{ds}=-\beta .
\end{align}

\noindent Observe that $\alpha $ tends to $1$ as $r$ goes to $0$
by (3.18) (hence ${\beta}=o(1)$ since $\alpha^{2}+\beta^{2}=1$).
So from (3.19) we can find a small neighborhood $V$ of $p_{0}$ so
that the distance
between $p_{0}$ and the characteristic curve $\Gamma $ past a point $%
p_{1}\in V\backslash \{p_{0}\}$ is decreasing towards the $-N^{\perp }$
direction. Let $s$ denote a unit-speed parameter of $\Gamma $ and $%
p_{1}=\Gamma (s_{1}).$ Then from the following formula

$$r(s)-r(s_{1})=\int_{s_{1}}^{s}\frac{dr}{ds}ds=\int_{s_{1}}^{s}(-%
\alpha )ds,$$

\noindent we learn that $r(s)$ reaches $0$ for a finite $s.$
\begin{flushright}
Q.E.D.
\end{flushright}

\bigskip

Let $B_{r}(p_{0})=\{p\in \Omega $ $|$ $|p-p_{0}|<r\}.$ Define $%
H_{M}(r)=\max_{p\in \partial B_{r}(p_{0})}|H(p)|.$

\bigskip

{\bf Theorem 3.10.} \textsl{Let }$u\in C^{2}(\Omega ).$\textsl{\ Suppose }$%
|H|=o(\frac{1}{r})$\textsl{\ (little ''}$o")$\textsl{\ near an isolated
singular point }$p_{0}\in \Omega $\textsl{\ where }$r(p)=|p-p_{0}|.$\textsl{%
\ Moreover, suppose there is }$r_{0}>0$\textsl{\ such that}

\begin{align}
\tag{3.20} \int_{0}^{r_{0}}H_{M}(r)dr<\infty .
\end{align}

\noindent \textsl{Then for any unit tangent vector }$v$\textsl{\ at }$p_{0},$\textsl{\
there exists a unique characteristic curve }$\Gamma $\textsl{\ touching }$%
p_{0}$\textsl{\ (i.e. }$p_{0}\in \bar{\Gamma},$\textsl{\ the closure of }$%
\Gamma $\textsl{) such that }$N^{\perp }(u)(p_{0}^{+}),$\textsl{\ the limit
of }$N^{\perp }(u)$\textsl{\ at }$p_{0}$\textsl{\ along }$\Gamma ,$\textsl{\
equals }$v.$\textsl{\ Moreover, there exists a neighborhood }$V$\textsl{\ of
}$p_{0}$\textsl{\ such that }$V\backslash \{p_{0}\}$\textsl{\ is contained
in the union of all such }$\Gamma $\textsl{'}$s.$

\bigskip

Proof. Take $\delta >0$ small enough so that all characteristic curves
past points on $\partial B_{\delta }(p_{0})$ reach $p_{0}$ in finite
unit-speed parameter in view of Lemma 3.9. Let $\Gamma $ denote the
characteristic curve past a point $p_{1}\in \partial B_{\delta }(p_{0}).$
Let $(s_{0},s_{1}]$ be the interval of unit-speed paramater describing
points of $\Gamma $ between $p_{0}$ and $p_{1}.$ Take a sequence of points $%
p_{j}\in \Gamma \rightarrow $ $p_{0}$ with parameter $s_{j}\rightarrow s_{0}$.
We compute

$$\theta (p_{j})-\theta (p_{k})=\int_{s_{k}}^{s_{j}}\frac{d\theta }{ds%
}ds=-\int_{s_{k}}^{s_{j}}Hds$$

\noindent by (2.23) (recall that $\theta $ is defined by $N=N(u)=(\cos \theta ,\sin
\theta )$). It follows from $\frac{dr}{ds}\rightarrow -1$ as $s\rightarrow
s_{0}$ or $r\rightarrow 0$ and (3.20) that

$$|\theta (p_{j})-\theta (p_{k})|\leq \int_{r_{k}}^{r_{j}}H_{M}(r)|%
\frac{ds}{dr}|dr\leq 2\int_{r_{k}}^{r_{j}}H_{M}(r)dr\rightarrow 0$$

\noindent as $p_{j},p_{k}\rightarrow p_{0}.$ This means that \{$\theta (p_{j})\}$ is a
Cauchy sequence. Therefore it converges to some number, denoted as $\theta
(p_{0};p_{1}).$ Define a map $\Psi $ $:$ $\partial B_{\delta }(p_{0})$ $%
\rightarrow $ $S^{1}$ by $\Psi (q)=\theta
(p_{0};q).$ We claim that $\Psi $ is a homeomorphism. Take a sequence of
points $q_{j}\in \partial B_{\delta }(p_{0})$ converging to $\hat{q}.$ We
want to show that $\theta (p_{0};q_{j})$ converges to $\theta (p_{0};\hat{q}%
).$ Without loss of generality, we may assume all $q_{j}^{\prime }s$ are
sitting on one side of $\hat{q}$ so that $\varphi (q_{1})>\varphi
(q_{2})>...>\varphi (\hat{q})$ where $\varphi $ is the angle in polar
coordinates centered at $p_{0}$ ranging in $[0,2\pi ).$ Observe that

\begin{align}
\tag{3.21} \theta (p_{0};q_{j})\geq \theta (p_{0};q_{j+1})\geq
...\geq \theta (p_{0};\hat{q})\\
\mbox{(letting }\theta \mbox{ take values in }
[0,2\pi ))\nonumber
\end{align}

\noindent for $j$ large since two distinct characteristic curves can not intersect in $%
B_{\delta }(p_{0})\backslash \{p_{0}\}$. Let $\hat{\theta}$ be the limit of $%
\theta (p_{0};q_{j})$ as $j\rightarrow \infty .$ Now suppose $\hat{\theta}%
\neq $ $\theta (p_{0};\hat{q})$ (hence $\hat{\theta}>\theta (p_{0};\hat{q})).
$ Let $\Gamma _{j}$ ($\hat{\Gamma}$, resp.) denote the characteristic curve
connecting $q_{j}$ ($\hat{q},$ resp.) and $p_{0}.$ Then we can find two rays
emitting from $p_{0}$ with angle smaller than $\hat{\theta}-\theta (p_{0};\hat{q})$
and a small positive $\hat{\delta}<\delta $ so that $\Gamma _{j}$ and $\hat{%
\Gamma}$ do not meet a fan-shaped region $\hat{\Omega}$ surrounded by these two
rays and $\partial B_{\hat{\delta}}(p_{0})$ for $j$ large. Take a point $%
\check{p}$ $\in $ $\hat{\Omega}.$ Consider the characteristic curve $\check{%
\Gamma}$ past $\check{p}.$ Then $\check{\Gamma}$ must intersect $\partial
B_{\delta }(p_{0})$ at a point $\check{q}$ while reaching $p_{0}$ with $%
\theta =\theta (p_{0};\check{q}).$ Since $\check{\Gamma}$ does not intersect
with any $\Gamma _{j},$ we have $\theta (p_{0};q_{j})>\theta (p_{0};\check{q}%
)$ for large $j.$ On the other hand, $\check{q}$ must coincide with $\hat{q}$
for the same reason. So $\check{\Gamma}=\hat{\Gamma},$ an obvious
contradiction. Thus $\hat{\theta}$ $=$ $\theta (p_{0};\hat{q}).$ We have
proved the continuity of $\Psi .$

Next we claim that $\Psi$ is surjective. If not, $S^{1}\backslash
\Psi (\partial B_{\delta }(p_{0}))$ is a nonempty open set. Then by a
similar fan-shaped region argument as shown above, we can reach a
contradiction. Let $\Gamma _{1},$ $\Gamma _{2}$ be two characteristic curves
past $q_{1},$ $q_{2}\in \partial B_{\delta }(p_{0})$ touching $p_{0}$ with $%
\theta (p_{0};q_{1})$ $=$ $\theta (p_{0};q_{2}).$ We want to show
that $q_{1}$ $=$ $q_{2}.$ Suppose $q_{1}$ $\neq $ $q_{2}.$ So
$\Gamma _{1}$ and $\Gamma _{2}$ are distinct (with empty
intersection in $B_{\delta }(p_{0})\backslash \{p_{0}\}$) and
tangent at $p_{0}.$ Let $\Omega _{r}$ denote the smaller domain,
surrounded by $\Gamma _{1},$ $\Gamma _{2}$, and $\partial
B_{r}(p_{0})$ for small $r>0.$ Then $\Omega _{r}$ is contained in
a fan-shaped region with vertex $p_{0}$ and angle $\theta _{r}$
such that $\theta _{r}\rightarrow 0$ as $r\rightarrow 0.$ Let
$\Gamma _{r}$ $\equiv $ $\partial \Omega _{r}\cap
\partial B_{r}(p_{0}).$ It follows that

\begin{align}
\tag{3.22} |\Gamma _{r}|\leq r\theta _{r}
\end{align}

\noindent where $|\Gamma _{r}|$ denotes the length of the arc $\Gamma _{r}.$ Since $%
N^{\perp }(u)$ is perpendicular to the unit outward normal $\nu $
($=\pm N(u)$) on $\Gamma _{1}$ and $\Gamma _{2},$ we obtain

\begin{align}
\tag{3.23} g(r)\equiv \oint_{\partial \Omega _{r}}(u_{y}+x,-u_{x}+y)\cdot
\nu ds=\int_{\Gamma _{r}}(u_{y}+x,-u_{x}+y)\cdot \nu ds.
\end{align}

\noindent Observing that $(u_{y}+x,-u_{x}+y)$ $=$ $p-p_{0}+o(r)$ by (3.17) and $%
\nu $ $=$ $\frac{p-p_{0}}{r}$ on $\Gamma _{r}$, we deduce from $(3.23)$ that

\bigskip

\begin{align}
\tag{3.24} g(r)=[r+o(r)]|\Gamma _{r}|
\end{align}

\noindent On the other hand, the divergence theorem tells us that

\begin{align}
\tag{3.25} g(r)&=\int \int_{\Omega
_{r}}[(u_{y}+x)_{x}+(-u_{x}+y)_{y}]dx\wedge dy\\
&=2\int \int_{\Omega _{r}}dx\wedge
dy=2\int_{0}^{r}|\Gamma _{\tau }|d\tau .\nonumber
\end{align}

\noindent It follows from (3.25) that $g^{\prime }(r)=2|\Gamma
_{r}|.$ Comparing this with (3.24), we obtain $\frac{g^{\prime
}(r)}{g(r)}=\frac{2}{r}+o(\frac{1}{r}).$ Therefore $g(r)$ $=$
$cr^{2}+o(r^{2})$ for some constant $c>0.$ However, inserting
(3.22) into (3.24) shows that $g(r)$ $=$ $o(r^{2}),$ i.e.,
$\frac{g(r)}{r^{2}}\rightarrow 0$ as $r\rightarrow 0.$ We have
reached a contradiction. So $q_{1}$ $=$ $q_{2}$ and hence $\Psi $
is injective. Next we will show that $\Psi ^{-1}$ is continuous.
Suppose this is not true. Then we can find a sequence of $q_{j}\in
\partial B_{\delta }(p_{0})$ converging to $\check{q}\neq \hat{q}$
while $\theta (p_{0};q_{j})$ converges to $\theta (p_{0};\hat{q})$
(may assume monotonicity (3.21) or reverse order for large $j$).
Take a point $\bar{q}$ $\in $ $\partial B_{\delta }(p_{0})$,
$\bar{q}\neq \check{q}$ and $\bar{q}\neq \hat{q},$ such that
$\theta (p_{0};q_{j})$ $\geq $ $\theta
(p_{0};\bar{q})$ $\geq $ $\theta (p_{0};\hat{q})$ for all large $j.$ Since $%
\lim \theta (p_{0};q_{j})$ $=$ $\theta (p_{0};\hat{q}),$ we must have $%
\theta (p_{0};\bar{q})$ $=$ $\theta (p_{0};\hat{q})$ contradicting the
injectivity of $\Psi .$ Altogether we have shown that $\Psi $ is a
homeomorphism. The theorem follows from this fact.

\begin{flushright}
Q.E.D.
\end{flushright}

\bigskip

\noindent {\bf Proof of Corollary 3.6:}

Let $N^{\perp}(p_{0}^{+})$ denote $N^{\perp}(u)(p_{0}^{+})$ for
simplicity. First we claim that $N^{\perp}(p_{0}^{+})$ can not be
tangent to $\Gamma_{s}$ at $p_{0}$. If yes, we have
$N^{\perp}(p_{0}^{+})\ (U(p_{0})-U^{T}(p_{0}))=0$ by Lemma 3.7
(1), (2). However, $U-U^{T}=\left(
\begin{array}{cc}
0 & -2 \\
2 & 0%
\end{array}%
\right) $. It follows that $N^{\perp}(p_{0}^{+})=0$, a
contradiction. So $N^{\perp}(p_{0}^{+})$ is transversal to
$\Gamma_{s}$, and $N^{\perp}(p_{0}^{+})U^{T}(p_{0})\neq 0$. In
fact, we have

\begin{align}
|N^{\perp}(p_{0}^{+})U^{T}(p_{0})|=2\tag{3.26}
\end{align}

\noindent by noting that the unit-length vector
$N^{\perp}(p_{0}^{+})$ is proportional to $(u_{xy}+1,-u_{xx})$ (if
$\neq 0$) and $(u_{yy},-u_{xy}+1)$ (if $\neq 0$) from the proof of
Lemma 3.7 or Proposition 3.5. Now suppose $\Gamma _{+}$ and
$\Gamma _{+}^{\prime}$ are two distinct (never intersect in
$B_{\epsilon}(p_{0})\backslash \{p_{0}\}$ for some small
$\epsilon$) characteristic curves contained in $B^{+}$ touching
$p_{0}$ (hence with the same tangent vector $N^{\perp}(p_{0}^{+})$
at $p_{0}$). Let $\Omega _{r}$ denote the domain, surrounded by
$\Gamma _{+},$ $\Gamma _{+}^{\prime}$, and $\partial B_{r}(p_{0})$
for $0<r<{\epsilon}.$ Then $\Omega _{r}$ is contained in a
fan-shaped region with vertex $p_{0}$ and angle $\theta _{r}$ such
that $\theta _{r}\rightarrow 0$ as $r\rightarrow 0.$ Let $\Gamma
_{r}$ $\equiv $ $\partial \Omega _{r}\cap
\partial B_{r}(p_{0}).$ Then (3.22) holds.

On the other hand,
$(u_x-y,u_y+x)=({\bigtriangleup}x,{\bigtriangleup}y)U^{T}(p_{0})+o(r^{2})$
by (3.11) and (3.12) while
$({\bigtriangleup}x,{\bigtriangleup}y)=rN^{\perp}(p_{0}^{+})+o(r)$
and ${\nu}=({\bigtriangleup}x,{\bigtriangleup}y)r^{-1}$ on $\Gamma
_{r}$ tends to $N^{\perp}(p_{0}^{+})$ as $r\rightarrow 0$. So from
(3.23) we compute

\begin{align}
g(r)&=\int_{\Gamma _{r}}(u_{y}+x,-u_{x}+y)\cdot \nu ds\tag{3.27}\\
&=\int_{\Gamma _{r}}DN^{\perp}\cdot \nu ds\nonumber\\
&=(2r+o(r))|\Gamma_{r}|. \nonumber
\end{align}

\noindent Here we have used (3.26) in the last equality. Now a
similar argument as in the proof of Theorem 3.10 by comparing
(3.25) with (3.27) gives $g(r)=cr+o(r)$ for a positive constant
$c$. However, substituting (3.22) into (3.27) shows that
$g(r)=o(r^{2})$. We have reached a contradiction. Therefore
$\Gamma _{+}$ must coincide with $\Gamma _{+}^{\prime}$. Similarly
we have a unique characteristic curve $\Gamma _{-}\subset B^{-}$
touching $\Gamma _{s}$ also at $p_{0}$ so that $\Gamma _{-}\cup
\{p_{0}\}\cup \Gamma _{+}$ forms a $C^{1}$ smooth curve in
$B_{\epsilon }(p_{0})$.

\begin{flushright}
Q.E.D.
\end{flushright}

The line integral in (3.23) has a geometric interpretation. Recall
that the standard contact form in the Heisenberg group $H_1$ is
$\Theta _{0}=dz+xdy-ydx.$ Let $\tilde{u}$ denote the
map: $(x,y)$ $\rightarrow $ $(x,y,u(x,y)).$ It is easy to see that $\tilde{u}%
^{\ast }\Theta _{0}$ $=$ $(u_{y}+x)dy$ $-$ $(y-u_{x})dx.$ Now it is clear
that the line integral in (3.23) is exactly the line integral of $\tilde{u}%
^{\ast }\Theta _{0}.$ Note that $\tilde{u}^{\ast }\Theta _{0}$ vanishes
along any characteristic curve. If we remove the condition (3.20) in
Theorem 3.10, $\theta (p_{0};p_{1})$ will not exist as shown in the
following example.

\medskip

{\bf Example.} Let $p_{0}=(0,0).$ Let $u=-\frac{r^{2}}{\log r^{2}}$ ($=0$ at
$p_{0})$ where $r^{2}=x^{2}+y^{2}.$ Write $u=f(r^{2}).$ A direct computation
shows that

\begin{align}
\tag{3.28} u_{x}=2xf^{\prime }(r^{2}),\ \ u_{y}=2yf^{\prime }(r^{2}),\\
\tag{3.29} D=r\sqrt{1+4(f^{\prime })^{2}},\\
\tag{3.30} u_{xx}=2f^{\prime }+4x^{2}f^{\prime \prime },
u_{xy}=4xyf^{\prime \prime },u_{yy}=2f^{\prime }+4y^{2}f^{\prime \prime }.
\end{align}

\noindent It is easy to see that $p_{0}$ is an isolated singularity (for a general $f).
$ Also $u$ is $C^{2}$ at $p_{0}$ and $u_{xx}$ $=$ $u_{xy}$ $=$ $u_{yy}$ $=$ $%
0$ at $p_{0}$ (for $f(t)=-\frac{t}{\log t},$ $t=r^{2}).$ Therefore $%
d(N^{\perp }D)_{p_{0}}$ is the identity transformation and the index of $%
N^{\perp }D$ is $+1.$ Noting that $(-\sin \varphi ,\cos \varphi
)=(-y,x)r^{-1},$ we compute $\beta =D^{-1}r^{-1}(u_{y},-u_{x})\cdot
(-y,x)=-2f^{\prime }/\sqrt{1+4(f^{\prime })^{2}}$ by (3.16), (3.28), and
(3.29). Therefore along a characteristic curve reaching $p_{0}$ in the $%
-N^{\perp }$ direction, we can estimate

\begin{align}
\tag{3.31} \frac{d\varphi }{ds}&=-\frac{1}{r}\beta \mbox{ (by (3.19))}\\
&=\frac{2f^{\prime }}{r\sqrt{1+4(f^{\prime })^{2}}}\approx \frac{2}{-r\log r^{2}}\nonumber
\end{align}

\noindent as $r\rightarrow 0$ for $f(t)=-\frac{t}{\log t},$ $t=r^{2}.$ Since $\frac{dr%
}{ds}\rightarrow -1$ by (3.18) and

$$ \int_{0}^{r_{1}}\frac{2}{-r\log r^{2}}dr=\infty ,$$

\noindent we conclude from (3.31) that $\varphi \rightarrow \infty ,$ hence $\theta
\rightarrow \infty $ as the point on the characteristic curve approaches $%
p_{0}$ (Observing that $\alpha \rightarrow 1$ and $\beta \rightarrow 0$ in
(3.15) as $r\rightarrow 0,$ we have the limit of $\theta $ equal to the
limit of $\varphi $ plus $\pi /2$ if one of the limits exists, hence another
limit exists too).

Next substituting (3.28) and (3.30) into (3.10) gives

\begin{align}
\tag{3.32} P(u)=2r^{2}[f^{\prime }+4(f^{\prime })^{3}+2r^{2}f^{\prime
\prime }].
\end{align}

\noindent By (2.10), (3.29), and (3.32), we obtain the p-mean curvature

\begin{align}
\tag{3.33} H=\frac{2(f^{\prime }+4(f^{\prime })^{3}+2r^{2}f^{\prime
\prime })}{r(1+4(f^{\prime })^{2})^{3/2}}.
\end{align}

\noindent For $f(t)=-\frac{t}{\log t},$ $t=r^{2},$ $f^{\prime }\approx -\frac{1}{\log t%
},$ $2r^{2}f^{\prime \prime }\approx \frac{2}{(\log t)^{2}}$ near $r=0.$
Inserting these estimates into (3.33) gives $H\approx -\frac{1}{r\log r}.$
It is now a straightforward computation to verify that $H=o(\frac{1}{r})$
and the integral in (3.20) for such an $H$ diverges.

\bigskip

\section{{\bf A Bernstein-type theorem and properly embedded p-minimal surfaces}}

Recall that the characteristic curves for a p-minimal surface $%
\Sigma $ in a pseudohermitian 3-manifold $M$ are Legendrian geodesics in $M$
by $(2.1).$ For $M=H_{1},$ we have

\bigskip

{\bf Proposition 4.1.} \textsl{The Legendrian geodesics in }$H_{1},$ \textsl{%
identified with }$R^{3},$\textsl{\ with respect to }$\nabla ^{p.h.}$\textsl{%
\ are straight lines}.

\bigskip

Proof. Write a unit Legendrian vector field $e_{1}=f\hat{e}_{1}+g\hat{e}_{2}$
with $f^{2}+g^{2}=1.$ Since $\nabla ^{p.h.}\hat{e}_{1}=\nabla ^{p.h.}\hat{e}%
_{2}=0,$ the geodesic equation $\nabla _{e_{1}}^{p.h.}e_{1}=0$ implies that $%
e_{1}f=e_{1}g=0.$ This means that $f=c_{1},g=c_{2}$ for some constants $%
c_{1},c_{2}$ along a geodesic $\Gamma $ (integral curve) of $e_{1}.$ We
compute

\begin{align}
\tag{4.1} e_{1}&=c_{1}\hat{e}_{1}+c_{2}\hat{e}%
_{2}=c_{1}(\partial _{x}+y\partial _{z})+c_{2}(\partial _{y}-x\partial _{z})\\
&=c_{1}\partial _{x}+c_{2}\partial
_{y}+(c_{1}y-c_{2}x)\partial _{z}.\nonumber
\end{align}

\noindent So $\Gamma $ is described by the following system of ordinary differential
equations:

\begin{align}
\tag{4.2a} \frac{dx}{ds}&=c_{1},\\
\tag{4.2b} \frac{dy}{ds}&=c_{2},\\
\tag{4.2c} \frac{dz}{ds}&=c_{1}y-c_{2}x.
\end{align}

\noindent By $(4.2a),(4.2b)$ we get $x=c_{1}s+d_{1},y=c_{2}s+d_{2}$ for some constants
$d_{1},d_{2}.$ Substituting into $(4.2c)$ gives $%
z=(c_{1}d_{2}-c_{2}d_{1})s+d_{3}$ for some constant $d_{3}.$ So $\Gamma $ is
a straight line in $R^{3}$.

\begin{flushright}
Q.E.D.
\end{flushright}

\bigskip

{\bf Corollary 4.2.} \textsl{The characteristic curves of a p-minimal surface in }$%
H_{1}$\textsl{\ are straight lines or line segments.}
\textsl{\ In particular, a characteristic curve (line) of a p-minimal surface in }
$H_{1}$\textsl{\ past a point }$q$\textsl{\ is contained in the contact plane past }$q$.

\bigskip

Recall that $F_{a,b}\equiv a(u_{x}-y)+b(u_{y}+x)$ for real constants $%
a,b$ with $a^{2}+b^{2}=1.$

\bigskip

{\bf Lemma 4.3.} \textsl{Suppose }$u\in C^{2}$\textsl{\ defines a p-minimal graph
near }$p\in S(u)$\textsl{, an isolated singular point. Then }$F_{a,b}=0,$%
\textsl{\ for} $a,b\in R,$\textsl{\ }$a^{2}+b^{2}=1$,\textsl{\ define all
straight line segments passing through }$p$\textsl{\ which are all
characteristic curves in a neighborhood of }$p$\textsl{\ with }$p$\textsl{\
deleted.}

\bigskip

Proof. First we claim $\nabla F_{a,b}(p)\neq 0$ for all $(a,b)$ with $%
a^{2}+b^{2}=1.$ If not, there exists $(a_{0},b_{0})$ such that $\nabla
F_{a_{0},b_{0}}(p)=0.$ So $detU(p)=0$ (see the paragraph after $(3.1)$)$.$
It follows from the proof of Theorem B that there is a small neighborhood of
$p$ which intersects with $S(u)$ in exactly a $C^{1}$ smooth curve past $p.$
This contradicts $p$ being an isolated singular point. We have shown that $%
\nabla F_{a,b}(p)\neq 0$ for all $(a,b)$ with $a^{2}+b^{2}=1.$ Therefore $%
F_{a,b}=0$ defines a $C^{1}$ smooth curve past $p$ for all $(a,b).$

In a neighborhood of $p$ with $p$ deleted, we observe that $%
F_{a,b}D^{-1}\equiv \sin \theta _{0}\cos \theta -\cos \theta _{0}\sin \theta
.$ Here we write $a=\sin \theta _{0},b=-\cos \theta _{0}.$ Recall that $%
(u_{x}-y)D^{-1}=\cos \theta ,$ $(u_{y}+x)D^{-1}=\sin \theta $ (see Section
2). So $\theta =\theta _{0}$ on $\{F_{a,b}=0\}$, and hence by $(2.18b)$ $%
N^{\perp }=(\sin \theta ,-\cos \theta )=(\sin \theta _{0},-\cos \theta _{0})$
is a constant unit vector field along $\{F_{a,b}=0\}.$ On the other hand, $%
\nabla (F_{a,b}D^{-1})=(-a\sin \theta +b\cos \theta )\nabla \theta $ is
parallel to $N$ $=(\cos \theta ,\sin \theta )$ since $N^{\perp }\cdot \nabla
\theta =0$ is our equation. It follows that $\{F_{a,b}=0\}$ is a straight
line segment and an integral curve of $N^{\perp }$ in a $p$-deleted
neighborhood.

\begin{flushright}
Q.E.D.
\end{flushright}

We remark that Lemma 4.3 provides a more precise description of
Theorem 3.10 in the case of $H=0$.

Since the characteristic curves are straight lines, we will often
call them characteristic lines (line segments). From Corollary
3.6, we know that a characteristic line keeps straight after it
{\em goes through} a singular curve. Note that two characteristic
line segments $\Gamma _{1},\Gamma _{2}$ can not touch a singular
curve at the same point $p_{0}$ unless they lie on a straight line
by Proposition 3.5 (the limits of $N(u)$ at $p_{0}$ along $\Gamma
_{1},\Gamma _{2}$ must be either the same or different by a sign).
We say a graph is entire if it is defined on the whole $xy$-plane.

\bigskip

{\bf Lemma 4.4.} \textsl{Suppose }$u\in C^{2}$\textsl{\ defines an entire
p-minimal graph. Then }$S(u)$\textsl{\ contains no more than one isolated
singular point.}

\bigskip

Proof. Suppose we have two such points $p_{1},p_{2}\in S(u).$ Then there
exist two distinct straight lines passing through $p_{1},p_{2}$,
respectively and intersecting at a third point $q$ such that $q\notin $ $%
S(u) $ in view of Theorem B and Corollary 3.6. From the proof
of Lemma 4.3 and Corollary 3.6, these two straight lines are
characteristic curves in the complement of $S(u)$, namely integral curves of
$N^{\perp }(u).$ But then at $q,$ $N^{\perp }(u)$ has two values, a
contradiction.

\begin{flushright}
Q.E.D.
\end{flushright}

On the other hand, remember that we can change the sign of $N^{\perp }(u)$ if
necessary to make a $C^{0}$ characteristic (i.e. tangent to integral curves
of $N^{\perp }(u)$ where $N^{\perp }(u)$ is defined) vector field $\check{N}%
^{\perp }(u)$ on the whole $xy$-plane except possibly one isolated
singular point in view of Proposition 3.5. Moreover, we have a
unique characteristic curve "going through" a point of a singular
curve by Corollary 3.6. So we can conclude that the following
result holds.

\bigskip

{\bf Lemma 4.5.} \textsl{Suppose }$u\in C^{2}$\textsl{\ defines an entire
p-minimal graph and $S(u)$ contains no isolated singular point. Then all
integral curves (restrict to characteristic lines of }$N^{\perp }(u)$\textsl{%
) of }$\check{N}^{\perp }(u)$ \textsl{are parallel.}

\bigskip

\noindent $\mathbf{Proof\ of\ Theorem\ A:}$

According to Lemma 4.4, we have the following two cases.

Case 1. $S(u)$ contains one isolated singular point.

In this case, we claim the solution $u$ is nothing but $(1.1)$. Let $%
p_{0}$ be the singular point. Let $r,\vartheta $ denote the polar
coordinates with center $p_{0}.$ We can write $\pm \check{N}^{\perp }(u)=%
\frac{\partial }{\partial r}$ in view of Lemma 4.3. By $(2.20)$ we have the
equation

\begin{align}
\tag{4.3} u_{rr}&=\frac{\partial ^{2}u}{\partial r^{2}}=0
\end{align}

\noindent defined on the whole $xy$-plane except $p_{0}$. It follows from $(4.3)$ that
$u=rf(\vartheta )+g(\vartheta )$ for some $C^{2}$ functions $f,g$ in $%
\vartheta .$ Since $u$ is continuous at $p_{0}=(x_{0},y_{0})$ (where $r=0$)$%
, $ $u(x_{0},y_{0})=g(\vartheta )$ for all $\vartheta .$ So $g$ is a
constant function, say $g=c.$ Also $f(\vartheta )=f(\vartheta +2\pi )$
implies that we can write $f(\vartheta )=\tilde{f}(\cos \vartheta ,\sin
\vartheta )$ where $\tilde{f}$ is $C^{2}$ in $\alpha =\cos \vartheta $ and $%
\beta =\sin \vartheta $. Compute $u_{x}=u_{r}r_{x}+u_{\vartheta }\vartheta
_{x}=\alpha \tilde{f}+\beta ^{2}\tilde{f}_{\alpha }-\alpha \beta \tilde{f}%
_{\beta }$ in which $\tilde{f}_{\alpha }=\partial \tilde{f}/\partial \alpha ,%
\tilde{f}_{\beta }=\partial \tilde{f}/\partial \beta $ , etc. and we have
used $\vartheta _{x}=-(\sin \vartheta )/r.$ Similarly we obtain $u_{y}=\beta
\tilde{f}+\alpha ^{2}\tilde{f}_{\beta }-\alpha \beta \tilde{f}_{\alpha }.$
Since $u_{x}$ and $u_{y}$ are continuous at $(x_{0},y_{0}),$ we immediately
have the following identities:

\begin{align}
\tag{4.4} \beta ^{2}\tilde{f}_{\alpha }-\alpha \beta \tilde{f}%
_{\beta }+\alpha \tilde{f}&=a \\
\tag{4.5} -\alpha \beta \tilde{f}_{\alpha }+\alpha ^{2}\tilde{%
f}_{\beta }+\beta \tilde{f}&=b
\end{align}

\noindent for all $\alpha ,\beta .$ Here $a=u_{x}(x_{0},y_{0}),b=u_{y}(x_{0},y_{0}).$
Multiplying $(4.4),(4.5)$ by $\alpha$, $\beta$, respectively and adding
the resulting identities, we obtain $(\alpha ^{2}+\beta ^{2})\tilde{f}%
=a\alpha +b\beta .$ It follows that $\tilde{f}=a\alpha +b\beta $ since $%
\alpha ^{2}+\beta ^{2}=1.$ We have shown that $u(x,y)$ $=$ $r(a\cos
\vartheta +b\sin \vartheta )+c$ $=$ $a(x-x_{0})+b(y-y_{0})+c_{0}$ $=$ $%
ax+by+(c-ax_{0}-by_{0})$ $=$ $ax+by+c.$ (In fact $(x_{0},y_{0})$
$=$ $(-b,a)$ from the definition of a singular point and the plane
$\{(x,y,u(x,y)\}$ is just the contact plane passing through
$(x_{0},y_{0})$). We can also give a geometric proof for Case 1 as
follows. Let $\xi_{0}$ denote the standard contact bundle over
$H_1$ (see the Appendix). Let $\Sigma$ denote the p-minimal
surface defined by $u$. Observe that the union of all
characteristic lines "going through" $p_0$, the isolated singular
point, (together with $p_0$) constitutes the contact plane
$\xi_{0}(p_{0})$ in view of Corollary 4.2 and Lemma 4.3. It
follows that $\xi_{0}(p_{0})\subset \Sigma$. So
$\Sigma=\xi_{0}(p_{0})$, an entire plane, since $\Sigma$ is also
an entire graph. We are done.

Case 2. $S(u)$ contains no isolated singular point.

In this case we claim $u$ is nothing but $(1.2).$ By Lemma 4.5 and
Lemma 3.4 we can find a rotation $\tilde{x}$ $=$ $ax+by,$ $\tilde{y}$ $=$ $%
-bx+ay$\textsl{\ }with\textsl{\ }$a^{2}+b^{2}$ $=$ $1$ such that

\begin{align}
    \tag{4.6} \check{N}^{\perp }(u)=\pm \frac{\partial }{\partial
\tilde{x}}.
\end{align}

\noindent By $(2.20)$ our equation reads $\tilde{u}_{\tilde{x}\tilde{x}}=0$ where $%
\tilde{u}(\tilde{x},\tilde{y})$ $=$ $u(x,y).$ It follows that

\begin{align}
    \tag{4.7} \tilde{u}=
\tilde{x}\tilde{y}+g(\tilde{y}),
\end{align}

\noindent for some $C^{2}$ smooth functions $f,g.$ From $(4.6)$ we know $N(u)=(0,\pm 1).$ By
the definition of $N(u)$ we obtain $\tilde{u}_{\tilde{x}}-$ $\tilde{y}=0.$
So $f(\tilde{y})=\tilde{y}.$ Substituting this into $(4.7)$ gives $\tilde{u}=%
\tilde{x}\tilde{y}+g(\tilde{y}),$ and hence $%
u=-abx^{2}+(a^{2}-b^{2})xy+aby^{2}+g(-bx+ay).$

\begin{flushright}
Q.E.D.
\end{flushright}

\medskip

We remark that the singular curve in Case 2 is defined by $\tilde{x}%
=-g^{\prime }(\tilde{y})/2,$ and this curve has only one connected component.

Next we will describe a general properly embedded p-minimal surface
in $H_{1},$ which may not be a graph. According to Proposition 4.1, such a
surface must be a properly embedded ruled surface with Legendrian (tangent
to contact planes) rulings when we view $H_{1}$ as $R^{3}$. We call a ruled
surface with Legendrian rulings a Legendrian ruled surface. Conversely, we
claim that a properly embedded Legendrian ruled surface is a properly
embedded p-minimal surface. First observe that a straight line $L$ past $%
p_{0}$ $=$ $(x_{0},y_{0},z_{0})$ pointing in a contact direction $c_{1}%
\hat{e}_{1}(p_{0})$ $+$ $c_{2}\hat{e}_{2}(p_{0}),$ $c_{1}^{2}+c_{2}^{2}=1,$
is tangent to the contact plane everywhere. Here $\hat{e}_{1}(p_{0})$ $=$ $%
\partial _{x}+y_{0}\partial _{z}$ or $(1,0,y_{0})$ and $\hat{e}_{2}(p_{0})$ $%
=$ $\partial _{y}-x_{0}\partial _{z}$ or $(0,1,-x_{0}).$ In fact we can
parametrize any point $p$ $=$ $(x,y,z)$ $\in L$ as follows:

\begin{align}
    \tag{4.8} (x,y,z)=(x_{0},y_{0},z_{0})+s[c_{1}\hat{e}_{1}(p_{0})+c_{2}%
\hat{e}_{2}(p_{0})]
\end{align}

\noindent for some $s\in R.$ The tangent vector at $p$ is just $c_{1}\hat{e}%
_{1}(p_{0}) $ $+$ $c_{2}\hat{e}_{2}(p_{0})$ which exactly equals $c_{1}%
\hat{e}_{1}(p)$ $+ $ $c_{2}\hat{e}_{2}(p)$ by a simple computation. So it is
a vector in the contact plane at $p.$ And $L$ is a Legendrian line. A
Legendrian ruled surface is generated by such Legendrian lines with its
characteristic field $e_{1}(p)$ $=$ $c_{1}\hat{e}_{1}(p_{0})$ $+$ $c_{2}%
\hat{e}_{2}(p_{0})$ $=$ $c_{1}\hat{e}_{1}(p)$ $+$ $c_{2}\hat{e}_{2}(p)$ with
$c_{1},c_{2}$ being constant along the characteristic line (or line segment)
past a nonsingular point $p.$ It follows that $\nabla _{e_{1}}^{p.h.}e_{1}$ $%
=$ $c_{1}\nabla _{e_{1}}^{p.h.}\hat{e}_{1}$ $+$ $c_{2}\nabla _{e_{1}}^{p.h.}%
\hat{e}_{2}$ $=$ $0$ since $\nabla _{{}}^{p.h.}\hat{e}_{j}=0,$ $j=1,2.$ By
(2.1) the p-mean curvature $H$ vanishes. So we have shown that a Legendrian
ruled surface is a p-minimal surface. Also an immersed Legendrian ruled
surface is the union of a family of curves of the form $(4.8),$ and has the following
expression:

\begin{align}
    \tag{4.9} (x_{0}(\tau ),y_{0}(\tau ),z_{0}(\tau ))+s[\sin \theta (\tau
)(1,0,y_{0}(\tau ))-\cos \theta (\tau )(0,1,-x_{0}(\tau ))].
\end{align}

\noindent Here $(x_{0}(\tau ),y_{0}(\tau ),z_{0}(\tau ))$ is a curve transverse to
rulings, and we have written $c_{1}(\tau )$ $=$ $\sin \theta (\tau )$ and $%
c_{2}(\tau )$ $=$ $-\cos \theta (\tau ).$

\medskip

{\bf Example.} In $(4.9)$ we take $\gamma (\tau )$ $\equiv $ $(x_{0}(\tau
),y_{0}(\tau ),z_{0}(\tau ))$ $=$ $(\cos \tau ,\sin \tau ,0)$ and $\theta
(\tau )$ $=$ $\tau ,$ $0\leq \tau <2\pi .$ It is easy to see that $%
e_{1}(\tau )$ $=$ $(\sin \tau ,-\cos \tau ,1)$ (note that $e_{1}$ is
independent of $s$). Compute $e_{1}(\tau _{1})\times e_{1}(\tau _{2})\cdot
(\gamma (\tau _{2})-\gamma (\tau _{1}))$ $=$ $(\sin \tau _{2}-\sin \tau
_{1})^{2}+(\cos \tau _{2}-\cos \tau _{1})^{2}.$ So $e_{1}(\tau _{1})\times
e_{1}(\tau _{2})\cdot (\gamma (\tau _{2})-\gamma (\tau _{1}))$ $=$ $0$ if
and only if $\tau _{1}=\tau _{2}.$ Now it is easy to see that this
Legendrian ruled surface is embedded. Let us write down the $x,y,z$
components as follows:

\begin{align}
    \tag{4.10} x(\tau ,s)=\cos \tau +(\sin \tau )s, y(\tau ,s)=\sin \tau
-(\cos \tau )s, z(\tau ,s)=s.
\end{align}

\noindent So $\partial _{\tau }(x,y,z)$ $=$ $(\frac{\partial x}{\partial \tau },\frac{%
\partial y}{\partial \tau },\frac{\partial z}{\partial \tau })$ $=$ $(-\sin
\tau +(\cos \tau )s,\cos \tau +(\sin \tau )s,0)$ and $\Theta _{0}\ (\partial
_{\tau }(x,y,z))$ $=$ $1+s^{2}\neq 0.$ This means that the tangent vector $%
\partial _{\tau }(x,y,z))$ is not annihilated by the contact form $\Theta
_{0}.$ Therefore $(4.10)$ defines a properly embedded p-minimal surface in $%
H_{1}$ with no singular points, which is not a vertical plane (i.e. having
the equation $ax+by$ $=$ $c$). In fact, eliminating the parameters $\tau $
and $s$ in $(4.10)$ gives the equation $z^{2}$ $=$ $x^{2}+y^{2}-1.$

For a Legendrian ruled surface of graph type, we can have an
alternative approach to show that it is p-minimal. Let $(x,y,u(x,y))$
describe such a Legendrian ruled surface. Suppose we can take $x$ as the
parameter of the rulings (straight lines) for simplicity. Then $%
d^{2}/dx^{2}\{u(x,y(x))\}$ $=$ $0$ along a ruling$.$ By the chain rule we have

\begin{align}
    \tag{4.11} r+2sa+ta^{2}=0
\end{align}

\noindent where $a$ $=$ $\frac{dy}{dx},$ $r=u_{xx},$ $s=u_{xy},$ and $t=u_{yy}.$ On
the other hand, along a Legendrian line, we have the contact form $%
dz+xdy-ydx $ $=$ $0.$ It follows that $\frac{dz}{dx}+xa-y$ $=$ $p+qa+xa-y$ $%
= $ $0$ where $p=u_{x},$ $q=u_{y}.$ So $a$ $=$ $-\frac{p-y}{q+x}$ (if $%
q+x=0, $ then $p-y=0).$ Substituting this into $(4.11)$ gives

$$(q+x)^{2}r-2(q+x)(p-y)s+(p-y)^{2}t=0$$

\noindent which is exactly $(pMGE).$ We remark that a general ruled surface satisfies
a third order partial differential equation (see page 225 in [Mo]. Solving $%
(4.11)$ for ''$a$'' in terms of $r,s,t,$ and substituting the result into $%
d^{3}/dx^{3}\{u(x,y(x))\}$ $=$ $0$ give such an equation).

\bigskip

\section{{\bf Comparison principle and uniqueness for the Dirichlet problem}}

Let $\Omega $ be a domain (connected and proper open subset) in the $xy$%
-plane. Let $u,v$ : $\Omega \rightarrow R$ be two $C^{1}$ functions. Recall
the singular set $S(u)$ $=$ $\{(x,y)\in \Omega $ $\mid $ $u_{x}-y=0,$ $%
u_{y}+x=0\}$ and $N(u)=[\nabla u+(-y,x)]D_{u}^{-1}$ where $D_{u}$ $=$ $\sqrt{%
(u_{x}-y)^{2}+(u_{y}+x)^{2}}$ (e.g. see $(2.15a)$)$.$

\bigskip

{\bf Lemma 5.1.} \textsl{Suppose we have the situation described above. Then the
equality }

\begin{align}
    \tag{5.1}(\nabla u-\nabla v)\cdot (N(u)-N(v))
 = \frac{D_{u}+D_{v}}{2}\mid N(u)-N(v)\mid ^{2}
\end{align}

\noindent \textsl{holds on }$\Omega \backslash (S(u)\cup S(v))$\textsl{. In
particular, }$(\nabla u-\nabla v)\cdot (N(u)-N(v))$\textsl{\ }$=$\textsl{\ }$%
0$\textsl{\ if and only if }$N(u)=N(v).$

\bigskip

Proof. Let $\vec{\alpha}=\nabla u+(-y,x),$ $\vec{\beta}=\nabla v+(-y,x).$
Noting that $N(u)=\frac{\vec{\alpha}}{|\vec{\alpha}|},$ $N(v)=\frac{\vec{%
\beta}}{|\vec{\beta}|}$ ($D_{u}=|\vec{\alpha}|,$ $D_{v}=|\vec{\beta}|$), we
compute

\begin{align}
    \tag{5.2} (\nabla u-\nabla v) \cdot (N(u)-N(v))=(\vec{\alpha}-\vec{%
\beta})\cdot (\frac{\vec{\alpha}}{|\vec{\alpha}|}-\frac{\vec{\beta}}{|\vec{%
\beta}|}) \\
=|\vec{\alpha}|+|%
\vec{\beta}|-\frac{\vec{\alpha}\cdot \vec{\beta}}{|\vec{\beta}|}-\frac{\vec{%
\alpha}\cdot \vec{\beta}}{|\vec{\alpha}|}=(|\vec{\alpha}|+|\vec{\beta}%
|)(1-\cos \vartheta ) \nonumber
\end{align}

\noindent in which $\cos \vartheta =\frac{\vec{\alpha}\cdot \vec{\beta}}{|\vec{\alpha}%
||\vec{\beta}|}.$ On the other hand, we compute $\mid N(u)-N(v)\mid ^{2}$ $=$
$\mid \frac{\vec{\alpha}}{|\vec{\alpha}|}-\frac{\vec{\beta}}{|\vec{\beta}|}%
\mid ^{2}$ $=$ $2-2\frac{\vec{\alpha}\cdot \vec{\beta}}{|\vec{\alpha}||\vec{%
\beta}|}$ $=$ $2(1-\cos \vartheta ).$ Substituting this into the right-hand
side of $(5.2)$ gives $(5.1).$

\begin{flushright}
Q.E.D.
\end{flushright}

Remark. For the prescribed mean curvature equation $divTu=H$ in $R^n$ where
$Tu=\frac{\nabla u}{\sqrt{1+|{\nabla u}|^{2}}}$, we have the following structural
inequality:

\begin{align}
(\nabla u-\nabla v)\cdot (Tu-Tv)&\geq
\frac{\sqrt{1+|{\nabla u}|^{2}}+\sqrt{1+|{\nabla v}|^{2}}}{2}|Tu-Tv|^{2}\nonumber\\
&\geq |Tu-Tv|^{2}.\nonumber
\end{align}

\noindent The above inequality was discovered by Miklyukov [Mik],
Hwang [Hw1], and Collin-Krust [CK] independently. Here we have
adopted Hwang's method to prove Lemma 5.1.

Next let $u\in C^{0}(\bar{\Omega}\backslash S_{1}),$ $v\in C^{0}(%
\bar{\Omega}\backslash S_{2}),$ i.e., $u,v$ are not defined (may blow up) on
sets $S_{1},$ $S_{2}\subset \Omega ,$ respectively. Let $S$ $\equiv $ $%
S_{1}\cup S_{2}\cup S(u)\cup S(v)$ where $S(u)\subset \Omega \backslash
S_{1},$ $S(v)\subset \Omega \backslash S_{2}.$

\bigskip

{\bf Theorem 5.2.} \textsl{Suppose }$\Omega $\textsl{\ is a bounded domain in the }%
$xy$-\textsl{plane and }${\mathcal H}_{1}(\bar{S})$\textsl{, the 1-dimensional
Hausdorff measure of }$\bar{S}$\textsl{, vanishes.}\textsl{\ Let }$u\in C^{0}(\bar{\Omega}%
\backslash S_{1})\cap C^{2}(\Omega \backslash S),$\textsl{\ }$v\in C^{0}(%
\bar{\Omega}\backslash S_{2})\cap C^{2}(\Omega \backslash S)$\textsl{\ such
that}

\begin{align}
\tag{5.3} divN(u)&\geq divN(v)\mbox{ in }
\Omega \backslash S,\\
\tag{5.4} u&\leq v \; \; \; \; \;
\mbox{    on }\partial \Omega \backslash S.
\end{align}

\noindent \textsl{Then }$N(u)=N(v)$\textsl{\ in }$\Omega ^{+}\backslash S$ \textsl{%
where }$\Omega ^{+}$\textsl{\ }$\equiv $\textsl{\ }$\{p\in \Omega $\textsl{\
}$\mid $ $u(p)-v(p)>0\}.$

\bigskip

Proof. First ${\mathcal H}_{1}(\bar{S})=0$ means that given any $\epsilon >0,$ we can find
countably many balls $B_{j,\epsilon }$, $j=1,2,...$ such that ${\bar S}\subset \cup
_{j=1}^{\infty }B_{j,\epsilon }$ and $\Sigma _{j=1}^{\infty }length(\partial
B_{j,\epsilon })$ $<$ $\epsilon $ and we can arrange $\cup _{j=1}^{\infty
}B_{j,\epsilon _{1}}\subset \cup _{j=1}^{\infty }B_{j,\epsilon _{2}}$ for $%
\epsilon _{1}<\epsilon _{2}.$ Since $\bar S$ is compact, we can
find finitely many $B_{j,\epsilon }$'s, say
$j=1,2,...,n({\epsilon})$, still covering $\bar S$. Suppose
$\Omega ^{+}\neq \emptyset .$ Then by Sard's theorem there exists
a sequence of positive number $\delta _{i}$
converging to $0$ as $i$ goes to infinity, such that $\Omega _{i}^{+}$ $%
\equiv $ $\{p\in \Omega $ $\mid $ $u(p)-v(p)>\delta _{i}\}$ $\neq $ $%
\emptyset $ and $\partial \Omega _{i}^{+}\backslash S$ is $C^{2}$ smooth.
Note that $\partial \Omega _{i}^{+}\cap \partial \Omega \subset S$ by $%
(5.4). $ Now we consider

$$I_{\epsilon }^{i}=\oint_{\partial (\Omega _{i}^{+}\backslash
\cup _{j=1}^{n({\epsilon})}B_{j,\epsilon })}\tan ^{-1}(u-v)(N(u)-N(v))\cdot \nu
ds$$

\noindent where $\nu$, $s$ denote the outward unit normal vector and the arc length
parameter, respectively. By the divergence theorem we have

\begin{align}
    \tag{5.5} I_{\epsilon }^{i}= & \int \int_{\Omega_{i}^{+}\backslash \cup
_{j=1}^{n({\epsilon})}B_{j,\epsilon }}\{\frac{1}{1+(u-v)^{2}}(\nabla u-\nabla
v)\cdot (N(u)-N(v))+ \\
& \tan ^{-1}(u-v) div (N(u)-N(v))\}dxdy. \nonumber
\end{align}

Observe that the second term in the right hand side of $(5.5)$ is
nonnegative by $(5.3).$ It follows from $(5.1)$\ and $(5.5)$ that

\begin{align}
\tag{5.6} I_{\epsilon }^{i}\geq \int \int_{\Omega _{i}^{+}\backslash
\cup _{j=1}^{n({\epsilon})}B_{j,\epsilon }}\{\frac{1}{1+(u-v)^{2}}(\frac{%
D_{u}+D_{v}}{2})|N(u)-N(v)|^{2}\}dxdy.
\end{align}

\noindent On the other hand, we can estimate

\begin{align}
\tag{5.7} I_{\epsilon }^{i}&\leq (\tan ^{-1}\delta _{i})\int_{\partial
\Omega _{i}^{+}\backslash (\cup _{j=1}^{n({\epsilon})}B_{j,\epsilon
})}(N(u)-N(v))\cdot \nu ds \\
&\ \ +\frac{\pi }{2}\cdot 2\cdot \Sigma _{j=1}^{n({\epsilon})}length(\partial B_{j,\epsilon })
\leq \pi \cdot \Sigma _{j=1}^{\infty
}length(\partial B_{j,\epsilon })<\pi \epsilon \nonumber
\end{align}

\noindent by noting that $\nu =-\frac{\nabla (u-v)}{|\nabla (u-v)|}$ and hence $%
(N(u)-N(v))\cdot \nu \leq 0$ by (5.1). If $N(u)\neq N(v)$ at some point $p$
in $\Omega ^{+}\backslash {\bar S}$, then $N(u)\neq N(v)$ in an open neighborhood $%
V $ of $p,$ contained in $\Omega _{i}^{+}$ for all large $i.$ Observe that
the measure of $V\backslash \cup _{j=1}^{\infty }B_{j,\epsilon }$ is bounded
from below by a positive constant independent of small enough $\epsilon $
and $i.$ Thus from $(5.6)$ $I_{\epsilon }^{i}\geq c,$ a positive constant
independent of small enough $\epsilon $ and large enough $i.$ Letting $%
\epsilon $ go to $0$ in $(5.7)$ will give us a contradiction. So $N(u)=N(v)$
in $\Omega ^{+}\backslash {\bar S}$ and hence in $\Omega ^{+}\backslash S$ by continuity.
\begin{flushright}
Q.E.D.
\end{flushright}

Remark. Theorem 5.2 is an analogue of Concus and Finn's general
comparison principles for the prescribed mean curvature equation
(cf. Theorem 6 in [CF]). In [Hw2] Hwang invoked the "$tan^{-1}$"
technique to simplify the proof of [CF]. Here we have followed the
idea of Hwang in [Hw2] to prove Theorem 5.2.

\bigskip

{\bf Lemma 5.3.} \textsl{Let }$u,v\in C^{2}(\Omega )\cap C^{0}(\bar{\Omega})$%
\textsl{\ where }$\Omega $\textsl{\ is a bounded domain in the }$xy$\textsl{%
-plane. Suppose }$N(u)=N(v)$\textsl{\ in }$\Omega \backslash (S(u)\cup
S(v)), $\textsl{\ }$u=v$\textsl{\ on }$\partial \Omega .$ \textsl{Then }$u=v$%
\textsl{\ in }$\Omega .$

\bigskip

Proof. Suppose $u\neq v$ in $\Omega .$ We may assume the set $\{ p
\in \Omega \vert u(p) > v(p) \}$ $\neq $ $\emptyset $ (otherwise,
interchange $u$ and $v).$ By Sard's theorem (e.g.,[St], noting
that $C^2$ is essential), there exists $\epsilon >0$
such that $\Omega _{\epsilon }$ $\equiv $ $\{p\in \Omega $ $|$ $%
u(p)-v(p)-\epsilon >0\}$ $\neq $ $\emptyset $ and $\Gamma _{\epsilon }$ $%
\equiv $ $\{p\in \Omega $ $|$ $u(p)-v(p)=\epsilon \}$ is $C^{2}$ smooth.
Note that $\bar{\Gamma}_{\epsilon }\cap \partial \Omega $ $=$ $%
\emptyset $ since $u=v$ on $\partial \Omega $ by assumption. Also $%
\Gamma _{\epsilon }$ is closed and bounded, hence compact. Therefore $\Gamma
_{\epsilon }$ is the union of (finitely-many) $C^{2}$ smooth loops. Choose one of them, and
denote it as $\Gamma _{\epsilon}^{\prime}.$ We claim

\begin{align}
    \tag{5.8} \frac{du}{ds}+x\frac{dy}{ds}-y\frac{dx}{ds}=0 {}
\end{align}

\noindent on $\Gamma _{\epsilon}^{\prime}$ where $s$ is a unit-speed parameter of $\Gamma
_{\epsilon}^{\prime}.$ For $p\in \Gamma _{\epsilon}^{\prime}\cap S(u),$ $(5.8)$ holds by the
definition of a singular point. For $p\in \Gamma _{\epsilon}^{\prime}\backslash
(S(u)\cup S(v)),$ we compute $N^{\perp }(u)u$ ($N^{\perp }(u)$ as an
operator acting on $u$) as follows:

\begin{align}
    \tag{5.9} N^{\perp }(u)u&=N^{\perp }(u)\cdot \nabla u=D^{-1}\{(\nabla
u)^{\perp }+(x,y)\}\cdot \nabla u \\
&=D^{-1}(x,y)\cdot \nabla u=D^{-1}(x,y)\cdot
\{\nabla u+(-y,x)\} \nonumber \\
&=(x,y)\cdot N(u). \nonumber
\end{align}

\noindent Similarly we can show

\begin{align}
    \tag{5.10} N^{\perp }(v)v=(x,y)\cdot N(v).
\end{align}

Since $N(u)=N(v),$ hence $N^{\perp }(u)=N^{\perp }(v)$ at $p,$ we conclude
that $N^{\perp }(u)(u-v)=0$ at $p$ by $(5.9)$ and $(5.10).$ This means that $%
N^{\perp }(u)$ is tangent to $\Gamma _{\epsilon}^{\prime}$ at $p.$ So $(5.8)$ holds
at $p$. For $p\in (\Gamma _{\epsilon}^{\prime}\backslash S(u))\cap S(v),$ $\vec{\beta%
}$ ($\equiv \nabla v+(-y,x)$) $=0.$ We observe that $\vec{\alpha}$ ($\equiv
\nabla u+(-y,x)$) $=\vec{\alpha}-\vec{\beta}$ $=$ $\nabla (u-v).$ This means
that $N(u)$ ($\equiv \frac{\vec{\alpha}}{|\vec{\alpha}|}$) is normal to $%
\Gamma _{\epsilon}^{\prime}.$ So again $N^{\perp }(u)$ is tangent to $\Gamma
_{\epsilon}^{\prime}$ at $p.$ Thus $(5.8)$ holds at $p.$ We have shown that $(5.8)$
holds for all $p\in \Gamma _{\epsilon}^{\prime}$. $\Gamma _{\epsilon}^{\prime}$ bounds a
domain, denoted as $\Omega _{\epsilon}^{\prime}.$ Now integrating $(5.8)$ over
$\Gamma _{\epsilon}^{\prime}$, we obtain that the area of $\Omega _{\epsilon}^{\prime}$
vanishes by the divergence theorem, an obvious contradiction.
\begin{flushright}
Q.E.D.
\end{flushright}

\bigskip

\noindent $\mathbf{Proof\ of\ Theorem\ C:}$

It follows from Theorem 5.2 and Lemma 5.3.
\begin{flushright}
Q.E.D.
\end{flushright}

\medskip

We can generalize Lemma 5.1 in the following form. Let $\Omega $ be
a domain in $R^{n}$. Let $u,v$ : $\Omega \rightarrow R$ be two $C^{1}$
functions. Let $\vec{F}$ be a $C^{0}$ vector field in $R^{n}.$ Define $S(u,%
\vec{F})$ $=$ $\{p\in \Omega $ $|$ $\nabla u+\vec{F}=0$ at $p\}$ and $S(v,%
\vec{F})$ similarly.

\bigskip

{\bf Lemma 5.1'.} \textsl{On }$\Omega \backslash \lbrack S(u,\vec{F})\cup S(v,\vec{%
F})],$\textsl{\ we have the following identity:}

$$(\nabla u-\nabla v)\cdot (\frac{\vec{\alpha}}{|\vec{%
\alpha}|}-\frac{\vec{\beta}}{|\vec{\beta}|})=(\frac{|\vec{\alpha}|+|\vec{%
\beta}|}{2})\mid \frac{\vec{\alpha}}{|\vec{\alpha}|}-\frac{\vec{\beta}}{|%
\vec{\beta}|}\mid ^{2}$$

\noindent \textsl{where }$\vec{\alpha}=\nabla u+\vec{F},$\textsl{\ }$\vec{\beta}%
=\nabla v+\vec{F}.$\textsl{\ }

\bigskip

In general a contact form $dz+\Sigma _{j=1}^{j=n}f_{j}dx^{j}$ in $R^{n+1}$
gives rise to an $\vec{F}$ $=$ $(f_{1},f_{2},...,f_{n})$ such that $\nabla
u+\vec{F}$ is the $R^{n}$-projection of the Legendrian normal to the graph $%
z=u(x^{1},x^{2},...,x^{n}).$ To generalize Theorem 5.2 to a domain $\Omega $
in $R^{n}$ and replace $N(u),$ $N(v)$ by $\frac{\vec{\alpha}}{|\vec{\alpha}|}%
,$ $\frac{\vec{\beta}}{|\vec{\beta}|}$, we will use $S_{\vec{F}}=S_{1}\cup
S_{2}\cup S(u,\vec{F})\cup S(v,\vec{F})$ instead of $S.$

\bigskip

{\bf Theorem 5.2'.} \textsl{Suppose }$\Omega $\textsl{\ is a bounded domain in }$%
R^{n}$\textsl{\ and }${\mathcal H}_{n-1}({\bar S}_{\vec{F}})=0.$\textsl{\ Let }$u\in C^{0}(%
\bar{\Omega}\backslash S_{1})\cap C^{2}(\Omega \backslash S_{\vec{F}}),$%
\textsl{\ }$v\in C^{0}(\bar{\Omega}\backslash S_{2})\cap C^{2}(\Omega
\backslash S_{\vec{F}}),$\textsl{\ and }$\vec{F}\in C^{1}(\Omega )\cap C^{0}(%
\bar{\Omega})$\textsl{\ such that}

$$div(\frac{\nabla u+\vec{F}}{|\nabla u+\vec{F}|}%
)\geq div(\frac{\nabla v+\vec{F}}{|\nabla v+\vec{F}|})\textsl{\ \ \ in }\Omega \backslash S_{\vec{F}},$$

$$ u\leq v\textsl{\ \ \ \ \ \
\ \ \ \ \ \ \ \ \ \ \ \ \ \ \ on }\partial \Omega \backslash S_{\vec{F}}.$$

\noindent \textsl{Then }$\frac{\nabla u+\vec{F}}{|\nabla u+\vec{F}|}=\frac{\nabla v+%
\vec{F}}{|\nabla v+\vec{F}|}$ \textsl{in }$\Omega ^{+}\backslash S_{\vec{F}%
} $ \textsl{where }$\Omega ^{+}\equiv \{p\in \Omega \ \mid u(p)-v(p)>0\}.$

\bigskip

The proof of Lemma 5.1' (Theorem 5.2', respectively) is similar to
that of Lemma 5.1 (Theorem 5.2, respectively). We can also generalize Lemma
5.3. Let $\Omega $ be a bounded domain in $R^{2m},$ $m\geq 1.$ Take two real
functions $u,v$ $\in C^{2}(\Omega )\cap C^{0}(\bar{\Omega})$. Let $\vec{%
\alpha}\equiv \nabla u$ $+$ $\vec{F}$ \ where $\vec{F}$ $=$ $%
(f_{1},f_{2},...,f_{2m})$\ is a $C^{1}$ smooth vector field on $\Omega .$
Define $\vec{F}$ $^{\ast }$ $\equiv $ $(f_{2},-f_{1},f_{4},-f_{3},$ $%
...,f_{2m},-f_{2m-1}).$ Denote $\frac{\vec{\alpha}}{|\vec{\alpha}|}$ by $N_{%
\vec{F}}(u)$ and the set $\{p\in \Omega $ $|$ $\vec{\alpha}=0\}$ by $S_{\vec{%
F}}(u).$

\bigskip

{\bf Lemma 5.3'.} \textsl{Suppose we have the situation as described above.
Suppose }$N_{\vec{F}}(u)=N_{\vec{F}}(v)$\textsl{\ in }$\Omega \backslash (S_{%
\vec{F}}(u)\cup S_{\vec{F}}(v)),$\textsl{\ }$u=v$\textsl{\ on }$\partial
\Omega ,$ \textsl{and }$ div \vec{F}^{\ast }>0$\textsl{\ a.e. in }
$\Omega .$ \textsl{Then }$u=v$\textsl{\ on } $\bar{\Omega}.$

\bigskip

Proof. Suppose the conclusion is not true. We may assume the set $\{p\in
\Omega $ $|$ $u(p)>v(p)\}\neq \emptyset .$ By Sard's theorem we can find a small
$\epsilon >0$ such that $\Omega _{\epsilon }$ $\equiv $ $\{p\in \Omega $ $|$ $%
u(p)-v(p)-\epsilon >0\}$ $\neq $ $\emptyset $ and $\Gamma _{\epsilon }$ $%
\equiv $ $\{p\in \Omega $ $|$ $u(p)-v(p)=\epsilon \} =\partial \Omega _{\epsilon}$
is $C^{2}$ smooth.
Let $\vec{\alpha}^{\ast }$ $\equiv $
$(u_{y_{1}},-u_{x_{1}},u_{y_{2}},-u_{x_{2}},...,u_{y_{m}},-u_{x_{m}})$ $+$ $%
\vec{F}$ $^{\ast }$(so that $\vec{\alpha}^{\ast }\cdot \vec{\alpha}=0$)$.$
Let $\nu =-\frac{\nabla (u-v)}{|\nabla (u-v)|}$ denote the outward unit
normal to $\Gamma _{\epsilon }$. We claim

\begin{align}
\tag{5.11} \vec{\alpha}^{\ast }\cdot \nu =0
\end{align}

\noindent on $\Gamma _{\epsilon }.$ Note that $\vec{\alpha}=0$ if and only if $\vec{%
\alpha}^{\ast }=0.$ So it is obvious that $(5.11)$ holds for $p\in S_{\vec{F}%
}(u)$. Let $N_{\vec{F}}^{\ast }(u)\equiv \frac{\vec{\alpha}^{\ast }}{|\vec{%
\alpha}^{\ast }|}$ for $p\in \Omega \backslash S_{\vec{F}}(u).$ In case $%
p\in \Gamma _{\epsilon }\backslash (S_{\vec{F}}(u)\cup S_{\vec{F}}(v)),$ we
can generalize $(5.9)$, $(5.10)$ as follows:

\begin{align}
\tag{$5.9^{\prime}$} N_{\vec{F}}^{\ast }(u)u&=\vec{F}^{\ast }\cdot N_{%
\vec{F}}(u),\\
\tag{$5.10^{\prime}$} N_{\vec{F}}^{\ast }(v)v&=\vec{F}^{\ast }\cdot N_{%
\vec{F}}(v).
\end{align}

\noindent Since $N_{\vec{F}}(u)=N_{\vec{F}}(v)$ by assumption, and hence $N_{\vec{F}%
}^{\ast }(u)=N_{\vec{F}}^{\ast }(v),$ we deduce from $(5.9^{\prime })$ and $%
(5.10^{\prime })$ that $N_{\vec{F}}^{\ast }(u)(u-v)=0.$ So $N_{\vec{F}%
}^{\ast }(u)$ is tangent to $\Gamma _{\epsilon }$ (at $p).$ This implies $%
(5.11).$ For $p\in (\Gamma _{\epsilon }\backslash S_{\vec{F}}(u))\cap S_{%
\vec{F}}(v),$ we still have $(5.11)$ by a similar argument as in the proof
of Lemma 5.3. We have proved $(5.11)$ for all $p$ $\in $ $\Gamma _{\epsilon
}.$ Let $dA$ denote the volume element of $\Gamma _{\epsilon },$ induced
from $R^{2m}.$ Now we compute

\begin{align}
\tag{5.12} 0&=\int_{\Gamma _{\epsilon }}\vec{\alpha}^{\ast }\cdot \nu dA
\mbox{ (by (5.11))}\\
&=\int_{\Omega _{\epsilon }}div(\vec{\alpha}%
^{\ast })\ d(volume)\mbox{ (by the divergence theorem)}\nonumber \\
&=\int_{\Omega _{\epsilon }}div\vec{F}^{\ast }d(volume)>0\nonumber
\end{align}

\noindent by assumption. We have reached a contradiction.
\begin{flushright}
Q.E.D.
\end{flushright}

\bigskip

We remark that the condition $div \vec{F}^{\ast }>0$ is essential in
Lemma 5.3'. Consider the case $\vec{F}=0$. Let $\Omega =B_{2}-{\bar B}_{1}$
where $B_{r}$ denotes the open ball of radius $r$. Let $u=f(r),v=g(r)$, and
$f\neq g$ with
the properties that $f(1)=g(1)$, $f(2)=g(2)$, and $f^{\prime}>0,g^{\prime}>0$
for $1\leq r \leq 2$. It follows that $S_{\vec{F}}(u)=\{ \nabla u =0\}
=\phi$, $S_{\vec{F}}(v)=\{ \nabla u =0\}=\phi$, and $\nabla u=f^{\prime}(r)
\nabla r$, $\nabla v=g^{\prime}(r)\nabla r$. Therefore we have

$$\frac{\nabla u}{|\nabla u|}=\nabla r=\frac{\nabla v}{|\nabla v|}$$

\noindent by noting that $|\nabla r|=1$. We have constructed a counterexample
for the statement of Lemma 5.3' if $div \vec{F}^{\ast }>0$ is not satisfied.

\medskip

For $\vec{F}$ $=$ $(-y_{1},x_{1},-y_{2},x_{2},...,-y_{m},x_{m}),$ we
have $\vec{F}^{\ast }$ $=$ $(x_{1},y_{1},x_{2},y_{2},...,x_{m},y_{m}).$ In
this case, we can view the integrand in $(5.12)$ geometrically:
($\widehat{dx_{j}\wedge dy_{j}}$ means deleting $dx_{j}\wedge dy_{j}$)

\begin{align}
\vec{\alpha}^{\ast }\cdot \nu dA&=\Sigma
_{j=1}^{j=m}[(u_{y_{j}}+x_{j})dy_{j}+(u_{x_{j}}-y_{j})dx_{j}]\wedge
\nonumber \\
&\ \ dx_{1}\wedge dy_{1}\wedge ...\wedge \widehat{dx_{j}\wedge dy_{j}}%
\wedge ...\wedge dx_{m}\wedge dy_{m}\nonumber \\
&=[du+\Sigma
_{j=1}^{j=m}(x_{j}dy_{j}-y_{j}dx_{j})]\wedge \nonumber \\
&\ \ (\Sigma
_{j=1}^{j=m}dx_{1}\wedge dy_{1}\wedge ...\wedge \widehat{dx_{j}\wedge dy_{j}}%
\wedge ...\wedge dx_{m}\wedge dy_{m})\nonumber \\
&=\frac{1}{2^{m-1}(m-1)!}\Theta
_{(m)}\wedge (d\Theta _{(m)})^{m-1}.\nonumber
\end{align}

\noindent Here $\Theta _{(m)}\equiv du+\Sigma _{j=1}^{j=m}(x_{j}dy_{j}-y_{j}dx_{j})$
is the standard contact form of the $2m+1$-dimensional Heisenberg group,
restricted to the hypersurface $\{(x_{1},y_{1},$ $x_{2},y_{2},$ $...,$ $%
x_{m},y_{m},$ $u(x_{1},y_{1},x_{2},y_{2},...,x_{m},y_{m})\}$. Integrating
the above form gives

\begin{align}
\int_{\Gamma _{\epsilon }}\vec{\alpha}^{\ast }\cdot \nu dA&=\frac{1}{%
2^{m-1}(m-1)!}\int_{\partial \Omega _{\epsilon }}\Theta _{(m)}\wedge
(d\Theta _{(m)})^{m-1}\nonumber \\
&=\frac{1}{2^{m-1}(m-1)!}%
\int_{\Omega _{\epsilon }}(d\Theta _{(m)})^{m} \mbox{ (by Stokes' Theorem)}
\nonumber \\
&=\frac{1}{2^{m-1}(m-1)!}%
2^{m}m!\ Volume(\Omega _{\epsilon })=2m\ Volume(\Omega _{\epsilon }).
\nonumber
\end{align}

\noindent In the last equality, we have used $d\Theta _{(m)}=2\Sigma
_{j=1}^{j=m}dx_{j}\wedge dy_{j}$ and hence $(d\Theta
_{(m)})^{m}=2^{m}m!dx_{1}\wedge dy_{1}\wedge ...\wedge dx_{m}\wedge dy_{m}.$
Note that $div\vec{F}^{\ast }$ $=$ $2m$ in this case. In general, let $%
\Theta _{\vec{F}}\equiv dz+\Sigma _{j=1}^{j=2m}f_{j}dx_{j}$ for $\vec{F}$ $=$
$(f_{1},f_{2},...,f_{2m}).$ We can easily compute

\begin{align}
&{d\Theta _{\vec{F}}\wedge \Sigma _{j=1}^{j=m}(dx_{1}\wedge
dx_{2}\wedge ...\wedge \widehat{dx_{2j-1}\wedge dx_{2j}}\wedge ...\wedge
dx_{2m-1}\wedge dx_{2m})}&\nonumber \\
&=(div\vec{F}^{\ast
})dx_{1}\wedge dx_{2}\wedge ...\wedge dx_{2m-1}\wedge dx_{2m}.&
\nonumber
\end{align}

\noindent Note that in case $\Theta _{\vec{F}}=\Theta _{(m)}$ (with $u,x_{j},y_{j}$
replaced by $z,x_{2j-1},x_{2j},$ respectively), we have

\begin{align}
 &(d\Theta _{\vec{F}})^{m-1}=(d\Theta _{(m)})^{m-1}& \nonumber \\
&=2^{m-1}(m-1)!\Sigma _{j=1}^{j=m}(dx_{1}\wedge dx_{2}\wedge
...\wedge \widehat{dx_{2j-1}\wedge dx_{2j}}\wedge ...\wedge dx_{2m-1}\wedge
dx_{2m}).&\nonumber
\end{align}

We can generalize Theorem C to higher dimensions without the condition on the
singular set.
Let $N(u)=$ $N_{\vec{F}}(u)$ and $S(u)=S_{\vec{F}}(u)$ for $\vec{F}$ $=$ $%
(-y_{1},x_{1},-y_{2},x_{2},...,-y_{m},x_{m})$.

\bigskip

{\bf Theorem C'.} \textsl{For a bounded domain }$\Omega $ \textsl{in }$%
R^{2m},\ m\geq 2,$\textsl{\ let }$u,v\in C^{2}(\Omega )\cap C^{0}(\bar{\Omega})$%
\textsl{\ satisfy }$divN(u)$\textsl{\ }$\geq $\textsl{\ }$divN(v)$\textsl{\
in }$\Omega \backslash S$\textsl{\ and }$u\leq v$\textsl{\ on }$\partial
\Omega $\textsl{\ where }$S=S(u)\cup S(v).$.\textsl{\ Then }$u\leq v$\textsl{\ in }$\Omega .$

\bigskip

First we observe the following size control of the singular set in general dimensions.

\bigskip

{\bf Lemma 5.4.} {\sl Suppose $u\in C^{2}(\Omega )$ where $\Omega\subset{R^{2m}}$. Then for any singular point $p\in S(u)$, there exists an open neighborhood $V\subset\Omega$
such that the $m$-dimensional Hausdorff measure of $S(u)\cap V$ is finite, and hence ${\mathcal H}_{2m-1}(S(u))=0$
for $m\geq 2$.}

\bigskip

Proof. Consider the map $G:p\in\Omega \rightarrow (\nabla u+\vec{F})(p)$.
Computing the differential $dG$ of $G$ at a singular point $p$
(where $G(p)=0$), we obtain

\begin{align}
\left( u_{ij}\right) +\left( \begin{array}{cccccc}
0 & -1 & 0 & 0 & . & .\\
1 & 0  & 0 & 0 & . & .\\
0 & 0 & 0 & -1 & . & .\\
0 & 0 & 1 & 0  & . & .\\
.&.&.&.&.&.\\
.&.&.&.&.&. \end{array}\right) \nonumber
\end{align}

\noindent in matrix form, where $(u_{ij})$ is the Hessian. Let $(dG)^{T}$
denote the transpose of $dG$. It is easy to see
that $2m=rank(dG-(dG)^{T})$ since $u_{ij}=u_{ji}$. On the other hand, $
rank(dG-(dG)^{T})\leq rank(dG)+rank(-dG)^{T}=2rank(dG)$. Therefore
$rank(dG)\geq m$. Hence the kernel of $dG$ has dimension $\leq m$. It follows
by the implicit function theorem that there exists an open neighborhood $V$
of $p$ such that $G^{-1}(0)\cap V=S(u)\cap V$ is a submanifold of $V$, having
dimension $\leq m$.

\begin{flushright}
Q.E.D.
\end{flushright}

\bigskip

\noindent $\mathbf{Proof\ of\ Theorem\ C'}$:

Observe that the condition ${\mathcal H}_{2m-1}(\bar{S})=0$ (the dimension
$n=2m$) in the proof of Theorem 5.2 (and Theorem 5.2') can be replaced by
the following condition: for any subdomain $O\subset\subset\Omega$, i.e.,
$\bar{O}\subset\Omega$, ${\mathcal H}_{2m-1}(\bar{O}\cap\bar{S})=0$. Since
$S=S(u)\cup S(v)$ is closed in the compact set $\bar{O}$,
$\bar{O}\cap\bar{S}=\bar{O}\cap S$. Now Theorem C' follows from Theorem 5.2'
(with the size control condition on $\bar{S}$ replaced by the above-mentioned
one), Lemma 5.3', and Lemma 5.4.

\begin{flushright}
Q.E.D.
\end{flushright}

\bigskip

\section{{\bf Second variation formula and area-minimizing property}}

In this section we will derive the second variation formula for the p-area
functional $(2.5)$ and examine the p-mean curvature $H$ from the viewpoint
of calibration geometry ([HL]). As a result we can prove the area-minimizing
property for a p-minimal graph in $H_{1}.$

We follow the notation in Section 2. We assume the surface $\Sigma $
is p-minimal. Let $f,g$ be functions with compact support away from the
singular set and the boundary of $\Sigma .$ Recall $T$ denotes the Reeb
vector field of $\Theta $ (see Section 2 or the Appendix). We compute the
second variation of $(2.5)$ in the direction $V=fe_{2}+gT:$

\begin{align}
    \tag{6.1} \delta _{V}^{2}\int_{\Sigma }\Theta \wedge
e^{1}=\int_{\Sigma }L_{V}^{2}(\Theta \wedge e^{1})
=\int_{\Sigma }i_{V}\circ d\{i_{V}\circ d(\Theta \wedge e^{1})\}.
\end{align}

\noindent Here we have used Stokes' theorem and the formula $L_{V}=i_{V}\circ d+d\circ
i_{V}$ and $d^{2}=0.$ By $(2.7)$ and $H$ $=$ $\omega (e_{1})$, we get

\begin{align}
    \tag{6.2} d(\Theta \wedge e^{1})=-H\Theta \wedge e^{1}\wedge
e^{2}.
\end{align}

We recall (see Section 2) to define a function $\alpha $ on $\Sigma \backslash S_{\Sigma }$ such that $%
\alpha e_{2}+T\in T\Sigma .$ Observe that $\{\alpha e_{2}+T,e_{1}\}$ is a
basis of $T(\Sigma \backslash S_{\Sigma }).$ So on $\Sigma \backslash
S_{\Sigma }$ we have

\begin{align}
    \tag{6.3} e^{2}\wedge e^{1}=\alpha \Theta \wedge e^{1}.
\end{align}

\noindent From $(6.2)$ it is easy to see that $i_{V}\circ d(\Theta \wedge
e^{1})=gHe^{2}\wedge e^{1}-fH\Theta \wedge e^{1}.$ Then taking $i_{V}\circ d$
of this expression and making use of (A.1r), (A.3r), $(6.3)$ and $H=0$ on $%
\Sigma $, we obtain

\begin{align}
    \tag{6.4} i_{V}\circ d\{i_{V}\circ d(\Theta \wedge
e^{1})\}&=(g\alpha -f)(gT+fe_{2})(H)\Theta \wedge e^{1}\\
&=-(g\alpha -f)^{2}e_{2}(H)\Theta \wedge e^{1}\nonumber
\end{align}

\noindent on $\Sigma .$ For the last equality we have used $T(H)=-\alpha e_{2}(H)$
since $\alpha e_{2}+T\in T\Sigma $ and $H=0$ on $\Sigma .$ Expanding the
left-hand side of (A.5r) gives

\begin{align}
    \tag{6.5} e_{2}(H)=2W+e_{1}(\omega (e_{2}))+2\omega (T)+(\omega
(e_{2}))^{2}.
\end{align}

\noindent Here we have used (A.6r) and $\omega (e_{1})=H=0$ on $\Sigma .$ The surfaces
$\varphi _{t}(\Sigma \backslash S_{\Sigma })$ are the level sets of a
defining function $\rho .$ Here $\dot{\varphi}_{t}$ $=$ $fe_{2}+gT.$ It
follows that $(fe_{2}+gT)(\rho )$ $=$ $1.$ On the other hand, $(\alpha
e_{2}+T)(\rho )$ $=$ $0$ from the definition of $\alpha .$ So $T(\rho )$ $=$
$-\alpha e_{2}(\rho )$ and $e_{2}(\rho )$ $=$ $(f-\alpha g)^{-1}$ (where $%
f-\alpha g$ $\neq $ $0$)$.$ Applying (A.6r) and (A.7r) to $\rho ,$ and using
the above formulas, we obtain

\begin{align}
    \tag{6.6a} \omega (e_{2})&=h^{-1}e_{1}(h)+2\alpha ,\\
\tag{6.6b} \omega (T)&=e_{1}(\alpha )-\alpha h^{-1}e_{1}(h)-%
\text{Im}A_{11}
\end{align}

\noindent where $h=f-\alpha g.$ Now substituting $(6.6a)$, $(6.6b)$ into $(6.5)$, we get

\begin{align}
    \tag{6.7} e_{2}(H)&=2W-2\text{Im}A_{11}+4e_{1}(\alpha )+4\alpha
^{2}\\
&\ \ +h^{-1}e_{1}^{2}(h)+2\alpha h^{-1}e_{1}(h).\nonumber
\end{align}

Observing that $e_{1}(e_{1}(h^{2}))\Theta \wedge e^{1}$ $=$ $\Theta
\wedge d(e_{1}(h^{2}))$ $=$ $-d(e_{1}(h^{2})\Theta )+2e_{1}(h^{2})\alpha
e^{1}\wedge \Theta $ on $\Sigma $ by (A.1r) and (6.3), we integrate $\frac{1%
}{2}e_{1}(e_{1}(h^{2}))$ $=$ $(e_{1}(h))^{2}+he_{1}^{2}(h)$ to obtain

$$-\int_{\Sigma }he_{1}^{2}(h)\Theta \wedge e^{1}=\int_{\Sigma
}[(e_{1}(h))^{2}+2\alpha he_{1}(h)]\Theta \wedge e^{1}.$$

\noindent Substituting $(6.7)$ into $(6.4)$ and $(6.4)$ into $(6.1)$ and using the above
formula$,$ we finally reach the following second variation formula.

\bigskip

{\bf Proposition 6.1.} \textsl{Suppose the surface }$\Sigma $\textsl{\ is
p-minimal as defined in Section 2. Let }$f,g$\textsl{\ be functions with
compact support away from the singular set and the boundary of }$\Sigma .$%
\textsl{\ Then}

\begin{align}
    \tag{6.8} &\delta _{fe_{2}+gT}^{2}\int_{\Sigma }\Theta
\wedge e^{1}& \\
&=\int_{\Sigma }\{(e_{1}(f-\alpha g))^{2}+(f-\alpha
g)^{2}[-2W+2\text{Im}A_{11}-4e_{1}(\alpha )-4\alpha ^{2}]\}\Theta \wedge
e^{1}.&\nonumber
\end{align}

Note that the Webster-Tanaka curvature $W$ and the torsion $A_{11}$
are geometric quantities of the ambient pseudohermitian 3-manifold $M.$ When
the torsion $A_{11}$ vanishes and $W$ is positive, we can easily discuss the
stability of a p-minimal surface (see Example 2 in Section 7). If both $W$
and $A_{11}$ vanish, e.g. in the case of $M$ $=$ $H_{1}$ (see the Appendix)$%
,$ we can compute $\alpha $, $e_{1}(\alpha )$ for a graph $z=u(x,y)$ as
follows. First note that the defining function $\rho =(z-u(x,y))D^{-1}$
satisfies the condition $e_{2}(\rho )=1$ (recall $e_{2}=-[(u_{x}-y)\hat{e}%
_{1}+(u_{y}+x)\hat{e}_{2}]D^{-1}$ in Section 2). So $\alpha =-T(\rho
)=-\partial \rho /\partial z=-D^{-1}$ and a direct computation shows (recall
$e_{1}=[-(u_{y}+x)\hat{e}_{1}+(u_{x}-y)\hat{e}_{2}]D^{-1}$ in Section 2) that

\begin{align}
    \tag{6.9} -4e_{1}(\alpha )-4\alpha
^{2}&=4\{(u_{x}-y)(u_{y}+x)(u_{xx}-u_{yy})+\\
&\ \ [(u_{y}+x)^{2}-(u_{x}-y)^{2}]u_{xy}\}D^{-4}.\nonumber
\end{align}

\noindent For example if $u=xy,$ the right-hand side of $(6.9)$ equals $1/x^{2}.$ So
away from the singular set, the second variation of the p-area is
nonnegative according to $(6.8)$ (it is easy to see that the second
variation in the $e_{1}$ direction always vanishes). Note that $\{x=0\}$ is
the singular set in this example.\ From the p-area minimizing property shown
below (Proposition 6.2), we know the second variation of the p-area for any
p-minimal graph over the $xy$-plane with no singular points must be
nonnegative.

If we consider only the variation in the $T$ direction, i.e., $f=0,$
we should combine the term in (6.9) with terms involving $e_{1}(\alpha )$ in
the expansion of $(e_{1}(-\alpha g))^{2}$ (to get a better expression of
(6.8))$.$ For instance, take a graph $(x,y,u(x,y))\in H_{1}$ over a domain $%
\Omega $ in the $xy$-plane$.$ We denote the energy functional for the p-area
by

$$E(u)=\int_{\Omega }Ddx\wedge dy$$

\noindent in view of (2.5) and (2.11). A direct computation shows that $\frac{d^{2}}{%
d\varepsilon ^{2}}\mid _{\varepsilon =0}E(u+\varepsilon
v)=\int_{\Omega}D^{-1}(e_{1}(v))^{2}dx\wedge dy$ for a variation
$v=v(x,y).$ On the other hand, this should be obtained from (6.8)
by letting $f=0$ and $g=v$ (with compact support away from the
singular set and $\partial \Omega ).$ It turns out to be
equivalent to verifying the following integral formula (note that
$\alpha $ $=$ $-D^{-1},$ $\Theta \wedge e^{1}$ $=$ $Ddx\wedge dy)$

$$\int_{\Omega
}\{e_{1}(v^{2})D^{-1}e_{1}(D^{-1})+v^{2}[(e_{1}(D^{-1}))^{2}+4D^{-2}e_{1}(D^{-1})-4D^{-4}]\}Ddx\wedge dy=0.$$

\noindent We leave this verification to the reader (Hint: we need an integration by
parts formula- $\int e_{1}(\varphi )\psi \Theta \wedge e^{1}=-\int [\varphi
e_{1}(\psi )+2\varphi \psi \alpha ]\Theta \wedge e^{1}$. Express $D^{2}$ $=$
$(e_{1}(\sigma ))^{2}+(e_{2}(\sigma ))^{2}$ where $\sigma $ $=$ $z-u(x,y).$
The following formulas: $e_{1}^{2}(\sigma )$ $=$ $e_{1}(\sigma )$ $=$ $0,$ $%
e_{2}(\sigma )$ $=$ $D,$ $[e_{1},e_{2}]$ $=$ $-2\partial _{z}-e_{2}(\theta
)e_{2},$ $e_{1}(D^{2})$ $=$ $-4D$ $-$ $2e_{2}(\theta )D^{2},$ and $%
e_{1}^{2}(D^{2})$ $=$ $2(e_{1}(D))^{2}+4e_{2}(\theta )D+4(e_{2}(\theta
)D)^{2}$ are useful).

For later use, we deduce a different expression for $\Xi $ $%
\equiv $ $[-2W+2\text{Im}A_{11}-4e_{1}(\alpha )-4\alpha ^{2}]\Theta \wedge
e^{1}$ in $(6.8).$ Noting that $e^{2}\wedge \Theta $ $=$ $0$ on $\Sigma ,$
we can easily get

\begin{align}
    \tag{6.10} d(\alpha \Theta )=[-e_{1}(\alpha )-2\alpha ^{2}]\Theta
\wedge e^{1}
\end{align}

\noindent by (A.1r) and (6.3). From (A.3r), (6.3), and (6.6a), we can relate $\text{Im}%
A_{11}$ to $\omega (T)$ as follows:

\begin{align}
    \tag{6.11} (\text{Im}A_{11})\Theta \wedge e^{1}=-de^{2}-[\omega
(T)+\alpha h^{-1}e_{1}(h)+2\alpha ^{2}]\Theta \wedge e^{1}.
\end{align}

\noindent In view of $(6.10)$ and $(6.11),$ we can express $\Xi $ in the following form:

\begin{align}
    \tag{6.12} \Xi =-2[W+\omega (T)+\alpha h^{-1}e_{1}(h)]\Theta \wedge
e^{1}+d(4\alpha \Theta -2e^{2}).
\end{align}

In Euclidean 3-geometry, we take the interior product of the
volume form with a vector field normal to a family of surfaces as
a calibrating form ([HL]). This 2-form restricts to the area form
on surfaces, and its exterior differentiation equals the mean
curvature times the volume form along a surface. We have analogous
results. Suppose $M$ is foliated by a family of surfaces $\Sigma
_{t}$, $-\varepsilon <t<\varepsilon .$ Let $e_{1}$ be a vector
field which is characteristic along any surface $\Sigma _{t}$. We
are assuming $\Sigma _{t}$ 's have no singular points. Let
$e_{2}=Je_{1}$ denote the Legendrian normal along each $\Sigma
_{t}$. Then the 2-form $\Phi =\frac{1}{2}i_{e_{2}}(\Theta \wedge
d\Theta )$ satisfies the following properties. First, a direct
computation shows that $\Phi =\Theta \wedge e^{1}$, our area
2-form from formula (A.1r). Secondly, $d\Phi =-H\Theta \wedge
e^{1}\wedge e^{2}$ by $(6.2)$. So $\{\Sigma _{t}\}$ are p-minimal
surfaces if and only if $d\Phi =0$. Now suppose this is the case
and $\Sigma ^{\prime }$ is a deformed surface with no singular
points near a p-minimal surface $\Sigma =\Sigma _{0}$ having the
same boundary. Also
suppose the Poincar\'{e} lemma holds. That is to say, there is a 1-form $%
\Psi $ such that $\Phi =d\Psi $. Then by Stokes' theorem, we have

\begin{align}
    \tag{6.13} p-Area(\Sigma )=\int_{\Sigma }\Phi
=\int_{\partial \Sigma }\Psi =\int_{\partial \Sigma ^{\prime }}\Psi
=\int_{\Sigma ^{\prime }}\Phi .
\end{align}

\noindent For $\Sigma ^{\prime }$, we have corresponding $e_{1}^{\prime
},e_{2}^{\prime },e^{1\prime },e^{2\prime }$. There is a function $\alpha
^{\prime }$ such that $T+\alpha ^{\prime }e_{2}^{\prime }$ is tangent to $%
\Sigma ^{\prime }$. Applying $\Phi =\Theta \wedge e^{1}$ to the basis $%
(T+\alpha ^{\prime }e_{2}^{\prime },e_{1}^{\prime })$ of $T\Sigma ^{\prime }$%
, we obtain $e^{1}(e_{1}^{\prime })$. It follows that $\Phi
=e^{1}(e_{1}^{\prime })\Theta \wedge e^{1\prime }$ when restricted to $%
\Sigma ^{\prime }$. So we have

\begin{align}
    \tag{6.14} \int_{\Sigma ^{\prime }}\Phi
&=\int_{\Sigma ^{\prime }}e^{1}(e_{1}^{\prime })\Theta \wedge e^{1\prime }\\
&\leq \int_{\Sigma
^{\prime }}\Theta \wedge e^{1\prime }=p-Area(\Sigma ^{\prime })
\mbox{ (since }e^{1}(e_{1}^{\prime })\leq 1).\nonumber
\end{align}

\noindent From $(6.13)$ and $(6.14),$ we have shown that

\begin{align}
    \tag{6.15} p-Area(\Sigma )\leq p-Area(\Sigma
^{\prime }).
\end{align}

Let us summarize the above arguments in the following

\bigskip

{\bf Proposition 6.2.} \textsl{Suppose we can foliate an open neighborhood of a
p-minimal surface }$\Sigma $\textsl{\ by a family of p-minimal surfaces with
no singular points, and in this neighborhood the Poincar\'{e} lemma holds
(i.e., any closed 2-form is exact). Then }$\Sigma $\textsl{\ has the
p-area-minimizing property. That is to say, if }$\Sigma ^{\prime }$\textsl{\
is a deformed surface with no singular points near }$\Sigma $\textsl{\
having the same boundary, then }$(6.15)$\textsl{\ holds.}

\bigskip

We remark that a p-minimal surface in $H_{1}$ with no singular
points, which is a graph over the $xy$-plane, satisfies the assumption in
Proposition 6.2. Note that a translation of such a p-minimal graph in the $%
z- $axis is still p-minimal (quantitatively $u+c$ is again a solution if $%
u=u(x,y)$ is a solution to $(pMGE))$. Also a vertical (i.e. perpendicular to
the $xy$-plane) plane in $H_{1}$ satisfies the assumption in Proposition
6.2. Note that a vertical plane is a p-minimal surface with no singular
points, and a family of parallel such surfaces surely foliates an open
neighborhood of a given one.

\bigskip

\section{{\bf Closed p-minimal surfaces in the standard $S^3$ and proof of Theorem E}}

First let us describe the standard pseudohermitian 3-sphere $(S^{3},\hat{J},%
\hat{\Theta})$ (see the Appendix for the definition of basic notions). The
unit 3-sphere $S^{3}$ in $C^{2}$ inherits a standard contact structure $\xi $
$=$ $TS^{3}\cap J_{C^{2}}(TS^{3})$ where $J_{C^{2}}$ denotes the almost
complex structure of $C^{2}.$ The standard $CR$ structure $\hat{J}$
compatible with $\xi $ is nothing but the restriction of $J_{C^{2}}$ on $\xi
.$ Let $r=|\zeta ^{1}|^{2}+|\zeta ^{2}|^{2}-1$ where $(\zeta ^{1},\zeta
^{2})\in C^{2}.$ The contact form $\hat{\Theta}$ $\equiv $ $-i\partial r$ $=$
$-i(\bar{\zeta}^{1}d\zeta ^{1}+\bar{\zeta}^{2}d\zeta ^{2})$ restricted to $%
S^{3}$ $\equiv $ $\{r=0\}$ gives rise to the Reeb vector field $\hat{T}$ $=$
$i\zeta ^{1}\partial _{\zeta ^{1}}+i\zeta ^{2}\partial _{\zeta ^{2}}-i\bar{%
\zeta}^{1}\partial _{\bar{\zeta}^{1}}-i\bar{\zeta}^{2}\partial _{\bar{\zeta}%
^{2}}.$ Take the complex vector field $\hat{Z}_{1}$ $=$ $\bar{\zeta}%
^{2}\partial _{\zeta ^{1}}-\bar{\zeta}^{1}\partial _{\zeta ^{2}}$ and the
complex 1-form $\hat{\theta}^{1}$ $=$ $\zeta ^{2}d\zeta ^{1}-\zeta
^{1}d\zeta ^{2}$ such that $\{\hat{\Theta},\hat{\theta}^{1},\hat{\theta}^{%
\bar{1}}\}$ is dual to $\{\hat{T},\hat{Z}_{1},\hat{Z}_{\bar{1}}\}$ and $d%
\hat{\Theta}=i\hat{\theta}^{1}\wedge \hat{\theta}^{\bar{1}}.$ It follows
that ${{\hat{\omega}}_{1}}^{1}$ $=$ $-2i\hat{\Theta},$ ${{\hat{A}}^{1}}_{\bar{1}%
} $ $=$ $0$ in the corresponding $(A.3)$ and $(A.4)$ in the Appendix$.$ Also
in the corresponding $(A.5),$ $\hat{W}$ $=$ $2.$ Write $\hat{Z}_{1}$ $=$ $%
\frac{1}{2}(\hat{e}_{1}-i\hat{e}_{2})$ for real vector fields $\hat{e}_{1},$
$\hat{e}_{2}.$ Let $\hat{\nabla}^{p.h.}$ denote the pseudohermitian
connection of $(\hat{J},\hat{\Theta}).$ From $\hat{\nabla}^{p.h.}\hat{Z}_{1}$
$=$ ${\hat{\omega}_{1}}^{1}$ $\otimes $ $\hat{Z}_{1}$ (see $(A.2)$), we have

\bigskip
\begin{align}
    \tag{7.1} \hat{\nabla}^{p.h.}\hat{e}_{1}=-2\hat{\Theta}\otimes \hat{e}%
_{2},\hat{\nabla}^{p.h.}\hat{e}_{2}=2\hat{\Theta}\otimes \hat{e}_{1}.
\end{align}

Recall that a Legendrian geodesic (with respect to $\hat{\nabla}%
^{p.h.}$) is a Legendrian curve $\gamma $ such that $\hat{\nabla}_{\dot{%
\gamma}}^{p.h.}\dot{\gamma}$ $=$ $0.$ Here $\dot{\gamma}=\frac{d\gamma }{ds}$
is the unit tangent vector with respect to the Levi metric and $s$ is a
parameter of unit speed. A Legendrian great circle of $(S^{3},\hat{J},\hat{%
\Theta})$ is a great circle in the usual sense, whose tangents belong to the
kernel of $\hat{\Theta}.$

\bigskip

{\bf Lemma 7.1.} \textsl{In }$(S^{3},\hat{J},\hat{\Theta})$\textsl{\ a Legendrian
geodesic is a part of a Legendrian great circle, and vice versa.}

\bigskip

Proof. Suppose $\gamma $ is a Legendrian geodesic\textsl{. }Write $\dot{%
\gamma}$ $=$ $a(s)\hat{e}_{1}+b(s)\hat{e}_{2}$ (note that $\hat{e}_{1}$ and $%
\hat{e}_{2}=\hat{J}\hat{e}_{1}$ belong to, and form a basis of, the kernel
of $\hat{\Theta}$)$.$ Compute $0$ $=$ $\hat{\nabla}_{\dot{\gamma}}^{p.h.}%
\dot{\gamma}$ $=$ $\dot{a}\hat{e}_{1}+\dot{b}\hat{e}_{2}$ since $\hat{\nabla}%
_{\dot{\gamma}}^{p.h.}\hat{e}_{1}$ $=$ $\hat{\nabla}_{\dot{\gamma}}^{p.h.}%
\hat{e}_{2}$ $=$ $0$ by $(7.1)$ and $\hat{\Theta}(\dot{\gamma})$ $=$ $0.$ So
$a$ ($b,$ respectively) is a constant $c_{1}$ ($c_{2}$, respectively) along $%
\gamma .$ Note that $c_{1}^{2}+c_{2}^{2}$ $=$ $1$ since $a^{2}+b^{2}$ $=$ $1$
by the unity of $\dot{\gamma}.$ Now write $\zeta ^{1}$ $=$ $x^{1}+iy^{1},$ $%
\zeta ^{2}$ $=$ $x^{2}+iy^{2}.$ From the definition we can express

\begin{align}
    \tag{7.2} \hat{e}_{1}&=x^{2}\partial _{x^{1}}-y^{2}\partial
_{y^{1}}-x^{1}\partial _{x^{2}}+y^{1}\partial _{y^{2}},\\
\hat{e}_{2}&=y^{2}\partial
_{x^{1}}+x^{2}\partial _{y^{1}}-y^{1}\partial _{x^{2}}-x^{1}\partial
_{y^{2}}.\nonumber
\end{align}

Writing $\gamma (s)=(x^{1}(s),y^{1}(s),x^{2}(s),y^{2}(s)),$ we can express
the equation $\dot{\gamma}$ $=$ $c_{1}\hat{e}_{1}+c_{2}\hat{e}_{2}$ by $%
(7.2) $ as

\begin{align}
    \tag{7.3}\dot{x}^{1}&=c_{1}x^{2}+c_{2}y^{2}, \dot{y}%
^{1}=c_{2}x^{2}-c_{1}y^{2},\\
\dot{x}^{2}&=-c_{1}x^{1}-c_{2}y^{1}, \dot{y}%
^{2}=-c_{2}x^{1}+c_{1}y^{1}.\nonumber
\end{align}

\noindent It is easy to see from $(7.3)$ that $\ddot{x}^{1}=-x^{1},$ $\ddot{y}%
^{1}=-y^{1},$ $\ddot{x}^{2}=-x^{2},$ $\ddot{y}^{2}=-y^{2}.$ Therefore

\begin{align}
    \tag{7.4} (x^{1}(s),y^{1}(s),x^{2}(s),y^{2}(s))=\cos (s) (\alpha
_{1},\alpha _{2},\alpha _{3},\alpha _{4})+\sin (s)(\beta _{1},\beta
_{2},\beta _{3},\beta _{4}).
\end{align}

\noindent Here the constant vector $(\beta _{1},\beta _{2},\beta _{3},\beta _{4})$ is
determined by the constant vector $(\alpha _{1},\alpha _{2},\alpha
_{3},\alpha _{4})$ as follows:

\begin{align}
    \tag{7.5} (\beta _{1},\beta _{2})=(\alpha _{3},\alpha _{4})\left(
\begin{array}{cc}
c_{1} & c_{2} \\
c_{2} & -c_{1}%
\end{array}%
\right) , (\beta _{3},\beta _{4})=(\alpha _{1},\alpha _{2})\left(
\begin{array}{cc}
-c_{1} & -c_{2} \\
-c_{2} & c_{1}%
\end{array}%
\right) .
\end{align}

\noindent Using $c_{1}^{2}+c_{2}^{2}$ $=$ $1,$ we have $\beta _{1}^{2}+\beta _{2}^{2}$
$=$ $\alpha _{3}^{2}+\alpha _{4}^{2},$ $\beta _{3}^{2}+\beta _{4}^{2}$ $=$ $%
\alpha _{1}^{2}+\alpha _{2}^{2}$ by $(7.5).$ Denote $(\alpha _{1},\alpha
_{2},\alpha _{3},\alpha _{4})$ by $\vec{\alpha}$ and $(\beta _{1},\beta
_{2},\beta _{3},\beta _{4})$ by $\vec{\beta}.$ We can write $(7.4)$ as

\begin{align}
    \tag{$7.4^{\prime}$} \gamma (s) = \cos (s)\vec{\alpha}+\sin (s)\vec{%
\beta}.
\end{align}

A direct computation using $(7.5)$ and $c_{1}^{2}+c_{2}^{2}$ $=$ $1$ shows
that $\vec{\alpha}$ is perpendicular to $\vec{\beta}$ and $|\vec{\alpha}|=|%
\vec{\beta}|=1$ in $R^{4}$ since $\gamma (s)\in S^{3}.$ It is now clear from
$(7.4^{\prime})$ that the Legendrian geodesic $\gamma (s),$ $0\leq s\leq
2\pi ,$ is a great circle. Moreover, $\vec{\beta}$ sits in the contact plane
at the point $\vec{\alpha}.$ Write $\hat{\Theta}$ $=$ $%
x^{1}dy^{1}-y^{1}dx^{1}$ $+$ $x^{2}dy^{2}-y^{2}dx^{2}$. Define $(e,f)^{\perp
}$ $=$ $(-f,e)$ for a plane vector $(e,f).$ Write $\vec{\alpha}$ $=$ $(\vec{%
\alpha}_{1},\vec{\alpha}_{2})$ and $\vec{\beta}$ $=$ $(\vec{\beta}_{1},\vec{%
\beta}_{2})$ for plane vectors $\vec{\alpha}_{1},\vec{\alpha}_{2},\vec{\beta}%
_{1},\vec{\beta}_{2}.$ Then $\vec{\beta}$ $\in $ $\ker \hat{\Theta}$ at the
point $\vec{\alpha}$ if and only if

\begin{align}
\tag{7.6} (\vec{\alpha}_{1}^{\perp },\vec{\alpha}_{2}^{\perp })\cdot (%
\vec{\beta}_{1},\vec{\beta}_{2})=\vec{\alpha}_{1}^{\perp }\cdot \vec{\beta}%
_{1}+\vec{\alpha}_{2}^{\perp }\cdot \vec{\beta}_{2}=0
\end{align}

\noindent (in which $"\cdot "$ denotes the inner product). Now it is easy to see that $%
(7.5)$ implies $(7.6)$ for $\vec{\alpha}_{1}$ $=$ $(\alpha _{1},\alpha
_{2}), $ $\vec{\alpha}_{2}$ $=$ $(\alpha _{3},\alpha _{4}),$ $\vec{\beta}%
_{1} $ $=$ $(\beta _{1},\beta _{2}),$ $\vec{\beta}_{2}$ $=$ $(\beta
_{3},\beta _{4}).$ Conversely, given an arbitrary point $\vec{\alpha}\in
S^{3}$ and a unit tangent vector $\vec{\beta}$ in the contact plane at $\vec{%
\alpha},$ we claim that the great circle $\gamma (s)$ defined by $%
(7.4')$ is Legendrian and a Legendrian geodesic. From $\dot{\gamma}%
(s)$ $=$ $-\sin (s)\vec{\alpha}+\cos (s)\vec{\beta},$ we compute

\bigskip

\ \ $\{\cos (s)(\vec{\alpha}_{1}^{\perp },\vec{\alpha}_{2}^{\perp })+\sin
(s)(\vec{\beta}_{1}^{\perp },\vec{\beta}_{2}^{\perp })\}\cdot \{-\sin (s)(%
\vec{\alpha}_{1},\vec{\alpha}_{2})+\cos (s)(\vec{\beta}_{1},\vec{\beta}%
_{2})\}$

\ \ \ \ \ \ \ $=-\sin ^{2}(s)(\vec{\beta}_{1}^{\perp }\cdot \vec{\alpha}_{1}+%
\vec{\beta}_{2}^{\perp }\cdot \vec{\alpha}_{2})+\cos ^{2}(s)(\vec{\alpha}%
_{1}^{\perp }\cdot \vec{\beta}_{1}+\vec{\alpha}_{2}^{\perp }\cdot \vec{\beta}%
_{2})$

\ \ \ \ \ \ \ $=\vec{\alpha}_{1}^{\perp }\cdot \vec{\beta}_{1}+\vec{\alpha}%
_{2}^{\perp }\cdot \vec{\beta}_{2}$ \ \ \ \ (since $\vec{\eta}\cdot \vec{%
\zeta}=\vec{\eta}^{\perp }\cdot \vec{\zeta}^{\perp }$ and $(\vec{\eta}%
^{\perp })^{\perp }=-\vec{\eta})$

\ \ \ \ \ \ \ $=0$ \ \ \ \ (by $(7.6)).$

\bigskip

\noindent So $\gamma (s)$ is Legendrian by $(7.6)$ again. From $(7.2)$ we can express $%
\hat{e}_{1},\hat{e}_{2}$ at $\gamma (s)$ as follows:

\bigskip

\ \ \ \ $\hat{e}_{1}(\gamma (s))=\cos (s)(\alpha _{3},-\alpha _{4},-\alpha
_{1},\alpha _{2})+\sin (s)(\beta _{3},-\beta _{4},-\beta _{1},\beta _{2}),$

\ \ \ \ $\hat{e}_{2}(\gamma (s))=\cos (s)(\alpha _{4},\alpha _{3},-\alpha
_{2},-\alpha _{1})+\sin (s)(\beta _{4},\beta _{3},-\beta _{2},-\beta _{1}).$

\bigskip

Recall that we write $\vec{\alpha}=(\alpha _{1},\alpha _{2},\alpha
_{3},\alpha _{4}),$ $\vec{\beta}=(\beta _{1},\beta _{2},\beta _{3},\beta
_{4}).$ Equating $(\dot{\gamma}(s)$ $=)$ $-\sin (s)\vec{\alpha}+\cos (s)\vec{%
\beta}$ $=$ $c_{1}\hat{e}_{1}(\gamma (s))$ $+$ $c_{2}\hat{e}_{2}(\gamma (s))$
with $c_{1}^{2}+c_{2}^{2}$ $=$ $1$ gives the equations $(7.5)$ (requiring $%
c_{1}^{2}+c_{2}^{2}$ $=$ $1$ gets rid of other equivalent equations).
Solving $(7.5)$ for $c_{1},c_{2},$ we obtain $c_{1}$ $=$ $(\beta _{1}\alpha
_{3}-$ $\beta _{2}\alpha _{4})(\alpha _{3}^{2}+\alpha _{4}^{2})^{-1},$ $%
c_{2} $ $=$ $(\beta _{1}\alpha _{4}+\beta _{2}\alpha _{3})(\alpha
_{3}^{2}+\alpha _{4}^{2})^{-1}$ if $\alpha _{3}^{2}+\alpha _{4}^{2}$ $\neq $
$0$ or $c_{1}$ $=$ $(-\alpha _{1}\beta _{3}+\alpha _{2}\beta _{4})(\alpha
_{1}^{2}+\alpha _{2}^{2})^{-1},$ $c_{2}$ $=$ $(-\alpha _{1}\beta _{4}-\alpha
_{2}\beta _{3})(\alpha _{1}^{2}+\alpha _{2}^{2})^{-1}$ if $\alpha
_{1}^{2}+\alpha _{2}^{2}$ $\neq $ $0$. Note that two expressions for $c_{1}$
($c_{2},$ respectively) are equal where $\alpha _{1}^{2}+\alpha _{2}^{2}$ $%
\neq $ $0$ and $\alpha _{3}^{2}+\alpha _{4}^{2}$ $\neq $ $0$ by the
condition $(7.6)$ and $\vec{\alpha}\cdot \vec{\beta}=0.$ Now $\hat{\nabla}_{%
\dot{\gamma}}^{p.h.}\dot{\gamma}$ $=$ $c_{1}\hat{\nabla}_{\dot{\gamma}%
}^{p.h.}\hat{e}_{1}(\gamma (s))$ $+$ $c_{2}\hat{\nabla}_{\dot{\gamma}}^{p.h.}%
\hat{e}_{2}(\gamma (s))$ $=$ $0$ by $(7.1).$ We have proved that $\gamma (s)$
is a Legendrian geodesic.
\begin{flushright}
Q.E.D.
\end{flushright}

\bigskip

Recall (see Section 2) that the p-mean curvature $H$ of a surface $%
\Sigma $ in a pseudohermitian 3-manifold $(M,J,\Theta )$ depends on $%
(J,\Theta ).$ Let $\tilde{H}$ denote the p-mean curvature associated to
another contact form $\tilde{\Theta}$ $=$ $\lambda ^{2}\Theta $, $\lambda >0$
with $J$ fixed. Let $e_{2}$ denote the Legendrian normal to $\Sigma .$ Then
we have the following transformation law.

\bigskip

{\bf Lemma 7.2.} \textsl{Suppose }$\tilde{\Theta}$\textsl{\ }$=$\textsl{\ }$%
\lambda ^{2}\Theta $\textsl{, }$\lambda >0.$\textsl{\ Then }$\tilde{H}%
=\lambda ^{-2}(\lambda H-3e_{2}(\lambda )).$

\bigskip

Proof. Let $\tilde{e}_{1}$ denote the characteristic field with respect to $%
\tilde{\Theta}.$ Then it follows from the definition that $\tilde{e}_{1}$ $=$
$\lambda ^{-1}e_{1}.$ Applying $(5.7)$ in [Lee] to $e_{1}$ (in our case, $%
n=1,$ $Z_{1}$ $=\frac{1}{2}(e_{1}-ie_{2})$), we obtain

\begin{align}
    \tag{7.7} \lambda {\tilde{\omega}_{1}}^{1}(\tilde{e}_{1})={%
\omega _{1}}^{1}(e_{1})-3i\lambda ^{-1}e_{2}(\lambda ).
\end{align}

\noindent Note that $H$ $=$ $\omega (e_{1})$ $=$ $-i{\omega _{1}}^{1}(e_{1})$ (see the
remark after (2.8)). Rewriting (7.7) in terms of $H$ and $\tilde{H}$ gives
what we want.
\begin{flushright}
Q.E.D.
\end{flushright}

\bigskip

We define the Cayley transform $F$: $S^{3}\backslash \{(0,-1)\}$ $%
\rightarrow $ $H_{1}$ by

$$x=\text{Re}(\frac{\zeta ^{1}}{1+\zeta ^{2}}),y=\text{Im}(\frac{%
\zeta ^{1}}{1+\zeta ^{2}}),z=\frac{1}{2}\text{Re}[i(\frac{1-\zeta ^{2}}{%
1+\zeta ^{2}})]$$

\noindent where $(\zeta ^{1},\zeta ^{2})\in S^{3}\subset C^{2}$ satisfies $|\zeta
^{1}|^{2}+|\zeta ^{2}|^{2}=1$ (see, e.g., [JL])$.$ A direct computation
shows that

\begin{align}
    \tag{7.8} \hat{\Theta}=F^{\ast }(\lambda ^{2}\Theta _{0})
\end{align}

\noindent where $\lambda ^{2}$ $=$ $4[4z^{2}+(x^{2}+y^{2}+1)^{2}]^{-1}$ (recall that $%
\Theta _{0}$ $=$ $dz+xdy-ydx$ is the standard contact form for $H_{1}$)$.$

\bigskip

{\bf Lemma 7.3.} \textsl{Let }$\Sigma $\textsl{\ be a }$C^{2}$\textsl{\ smoothly
embedded p-minimal surface in }$(S^{3},\hat{J},\hat{\Theta}).$\textsl{\
Suppose }$p\in \Sigma $\textsl{\ is an isolated singular point. Then there
exists a neighborhood }$V$\textsl{\ of }$p$\textsl{\ in }$\Sigma $\textsl{\
such that }$V$\textsl{\ is contained in the union of all Legendrian great
circles past }$p.$

\bigskip

Proof. Without loss of generality, we may assume $(0,-1)$ $\notin $ $\Sigma
$. Consider $\Sigma _{0}$ $=$ $F(\Sigma )$ $\subset $ $H_{1}.$ Let $H_{0}$
denote the p-mean curvature of $\Sigma _{0}$ in $H_{1}.$ By $(7.8)$ and
Lemma 7.2, we obtain $H_{0}$ $=$ $3\lambda ^{-1}e_{2}(\lambda )$. Here $%
e_{2} $ is the Legendrian normal of $\Sigma _{0}.$ Recall (Section 2) that $%
e_{2}$ $=$ $-(\cos \theta )\hat{e}_{1}$ $-$ $(\sin \theta )\hat{e}_{2}.$ So $%
e_{2}(\lambda )$ $=$ $-(\cos \theta )\hat{e}_{1}(\lambda )$ $-$ $(\sin
\theta )\hat{e}_{2}(\lambda )$ is bounded near the isolated singular point $%
F(p)$ since $\hat{e}_{1}(\lambda )$ and $\hat{e}_{2}(\lambda )$ are global $%
C^{\infty }$ smooth functions. Thus $H_{0}$ is bounded near $F(p)$ (Note
that near a singular point, $\Sigma _{0}$ is a graph over the $xy$-plane).
Observing that characteristic curves and the singular set are preserved
under the contact diffeomorphism $F,$ we conclude the proof of the lemma by
Theorem 3.10 and Lemma 7.1 .
\begin{flushright}
Q.E.D.
\end{flushright}

\noindent $\mathbf{A\ direct\ proof\ of\ Corollary\ F}$ (for the nonexistence
of hyperbolic p-minimal surfaces):

Let $\Sigma $\ be a closed, connected, $C^{2}$ smoothly embedded
p-minimal surface in $(S^{3},\hat{J},\hat{\Theta}).$ Without loss of
generality, we may assume $(0,-1)$ $\notin $ $\Sigma .$ Consider $\Sigma
_{0} $ $=$ $F(\Sigma )$ $\subset $ $H_{1}.$ As argued in the proof of Lemma
7.3, $H_{0}$, the p-mean curvature of $\Sigma _{0},$ is bounded$.$ According
to Theorem B any singular point of $\Sigma _{0}$ is either isolated or
contained in a $C^{1}$ smooth singular curve with no other singular points
near it (note that near a singular point, $\Sigma _{0}$ is a graph over the $%
xy$-plane). So back to $\Sigma $ through $F^{-1},$ we can only
have isolated singular points or closed singular curves on $\Sigma
.$ Similarly via the Cayley transform we can have an extension
theorem analogous to Corollary 3.6 with the characteristic curve
replaced by a
characteristic (Legendrian) great circle (arc) in the case of $(S^{3},\hat{J}%
,\hat{\Theta}).$

Now suppose $\Sigma $ has an isolated singular point $p.$ By Lemma
7.3 there exists a neighborhood $V$ of $p$ in $\Sigma $ such that $V$\textsl{%
\ }is contained in the union of all Legendrian great circles past
$p.$ But this union simply forms a p-minimal 2-sphere. This
2-sphere must be the whole $\Sigma $ since $\Sigma $ is connected.
Next suppose $\Sigma $ does not have any isolated singular point.
Then in view of the extension theorem, the space of leaves
(Legendrian great circles) of the characteristic foliation
(including touching points on singular curves) forms a closed,
connected, 1-dimensional manifold. So it must be homeomorphic to
$S^{1}.$ In this case, $\Sigma $ is topologically a torus.
\begin{flushright}
Q.E.D.
\end{flushright}

{\bf Example 1.} Every coordinate sphere (defined by $x^{1},y^{1},x^{2}$
or $y^{2}=0$) is a closed, connected, embedded p-minimal surface of genus 0
in $(S^{3},\hat{J},\hat{\Theta}).$ For instance, we can write $\{y^{2}=0\}$
as the union of Legendrian great circles: $\gamma _{t}(s)$ $=$ $\cos
(s)(0,0,1,0)$ $+$ $\sin (s)(\cos (t),\sin (t),0,0)$ parametrized by $t$ (it
is a simple exercise to verify $(7.6)$).

\medskip

{\bf Example 2.} Write $\zeta ^{1}=\rho _{1}e^{i\varphi _{1}},$ $\zeta
^{2}=\rho _{2}e^{i\varphi _{2}}$ in polar coordinates with $\rho
_{1}^{2}+\rho _{2}^{2}=1$ on $S^{3}.$ We consider the surface $\Sigma _{c}$
defined by $\rho _{1}=c,$ a constant between $0$ and $1$. Note that $\hat{%
\Theta}$ $=$ $\rho _{1}^{2}d\varphi _{1}+\rho _{2}^{2}d\varphi _{2}$, $%
\hat{J}\partial _{\varphi _{j}}$ $=$ $-\rho _{j}\partial _{\rho _{j}}$ and $%
\hat{J}\partial _{\rho _{j}}$ $=$ $\rho _{j}^{-1}\partial _{\varphi _{j}},$ $%
j=1,2$ in polar coordinates. Here we have used $\partial _{\varphi _{j}}$ $=$
$i\zeta ^{j}\partial _{\zeta ^{j}}-i\bar{\zeta}^{j}\partial _{\bar{\zeta}%
^{j}} $, $\partial _{\rho _{j}}$ $=$ $\rho _{j}^{-1}(\zeta ^{j}\partial
_{\zeta ^{j}}+\bar{\zeta}^{j}\partial _{\bar{\zeta}^{j}})$ for $j=1,2$ (no
summation convention here). Next we compute
the Reeb vector field $T,$ the characteristic field $e_{1},$ and the
Legendrian normal $e_{2}$ as follows:

$$T=\partial _{\varphi _{1}}+\partial _{\varphi _{2}},e_{1}=\frac{\rho
_{2}}{\rho _{1}}\partial _{\varphi _{1}}-\frac{\rho _{1}}{\rho _{2}}\partial
_{\varphi _{2}},e_{2}=\hat{J}e_{1}=-\rho _{2}\partial _{\rho _{1}}+\rho
_{1}\partial _{\rho _{2}}.$$

\noindent We then have $e^{1}$ $=$ $\rho _{1}\rho _{2}(d\varphi _{1}-d\varphi _{2})$
and $e^{2}$ $=$ $-\rho _{2}d\rho _{1}+\rho _{1}d\rho _{2}.$ So we can
compute the p-area 2-form $\hat{\Theta}\wedge e^{1}$ $=$ $-\rho _{1}\rho
_{2}d\varphi _{1}\wedge d\varphi _{2},$ the volume form $\hat{\Theta}\wedge
e^{1}\wedge e^{2}$ $=$ $\rho _{1}d\rho _{1}\wedge d\varphi _{1}\wedge
d\varphi _{2}$ and $d(\hat{\Theta}\wedge e^{1})$ $=$ $-(\rho _{2}^{2}-\rho
_{1}^{2})\rho _{2}^{-1}d\rho _{1}\wedge d\varphi _{1}\wedge d\varphi _{2}$
(noting that $\rho _{1}d\rho _{1}+\rho _{2}d\rho _{2}$ $=$ $0).$ By (6.2) we
obtain the p-mean curvature $H$ $=$ $(\rho _{2}^{2}-\rho _{1}^{2})(\rho
_{1}\rho _{2})^{-1}$ $=$ $\rho _{2}\rho _{1}^{-1}-\rho _{1}\rho _{2}^{-1}$
for $\Sigma _{c}$ where $0<\rho _{1}=c<1$ and $\rho _{2}$ $=$ $\sqrt{1-c^{2}}%
.$ Thus for $c=\sqrt{2}/2$ ($\rho _{1}$ $=$ $\rho _{2}$ $=$ $c$)$,$ $\Sigma
_{c}$ is a closed, connected, embedded p-minimal torus with no singular
points (observing that $T$ is tangent at every point of $\Sigma _{c}$) and
the union of Legendrian great circles defined by $\varphi _{1}+\varphi
_{2}=a $, $0\leq a<2\pi $ with $0$ identified with $2\pi .$ In any case, $%
\Sigma _{c}$ is a torus of constant p-mean curvature. Also note that the
p-minimal torus $\Sigma _{\sqrt{2}/2}$ is not stable. This can be seen from
the second variation formula $(6.8)$ in which $\alpha =0,$ $A_{11}=0,$ and $%
W $ $=$ $2.$

\medskip

We can generalize Corollary F to the situation that the ambient
pseudohermitian 3-manifold is spherical. A pseudohermitian 3-manifold is
called spherical if it is locally CR equivalent to $(S^{3},\hat{J})$.

\bigskip

\noindent $\mathbf{Proof\ of\ Theorem\ E:}$
Locally near a singular point of $\Sigma ,$ we
may assume that $\Sigma $ is a $C^{2}$ smooth graph over the $xy$%
-plane in $(R^{3},J_{0},\lambda ^{2}\Theta _{0})$ for ($C^{\infty }$ smooth)
$\lambda >0.$ Here $(J_{0},\Theta _{0})$ denotes the standard
pseudohermitian structure of $H_{1}$ (see the Appendix). By Lemma 7.2, the
p-mean curvature of $\Sigma $ with respect to $(J_{0},\Theta _{0})$ equals $%
3\lambda ^{-1}e_{2}(\lambda )+\lambda$(original p-mean curvature)
which is a bounded function by assumption and the boundedness of $e_{2}(\lambda )$
(for the same reason as in the proof of Lemma 7.3). So the singular set $%
S_{\Sigma }$ (depending only on $\Sigma $ and the contact
structure) consists of finitely many isolated points and $C^{1}$
smooth closed curves in view of Theorem B. Also the extension
theorems (Corollary 3.6 and Theorem 3.10) hold in this situation.
Now the configuration of characteristic foliation on $\Sigma $ is
clear. The associated line field (extended to include those
defined on points of singular curves) has only isolated singular
points of index 1 in view of Lemma 3.8. Therefore the total index
sum of this line field is nonnegative.
This index sum is equal to the Euler characteristic number of the surface $%
\Sigma $ according to the Hopf index theorem for a line field (e.g., [Sp]).
On the other hand, the Euler characteristic number of $\Sigma $ equals $%
2-2g(\Sigma )$ where $g(\Sigma )$ denotes the genus of $\Sigma $. It follows
that $g(\Sigma )\leq 1.$
\begin{flushright}
Q.E.D.
\end{flushright}

\bigskip

\noindent $\mathbf{Appendix.\ Basic\ facts\ in\ pseudohermitian\ geometry}$

\medskip

Let $M$ be a smooth (paracompact) 3-manifold. A contact
structure or bundle $\xi $ on $M$ is a completely nonintegrable plane
distribution. A contact form is a 1-form annihilating $\xi $.
Let $(M,\xi )$ be a contact 3-manifold with an oriented contact structure $%
\xi $. We say a contact form $\Theta $ is oriented if $d\Theta (u,v)>0$ for $%
(u,v)$ being an oriented basis of $\xi $. There always exists a global
oriented contact form $\Theta $, obtained by patching together local ones
with a partition of unity. The Reeb vector field of $\Theta $
is the unique vector field $T$ such that $\Theta (T)=1$ and $\mathcal{L}%
_{T}\Theta =0$ or $d\Theta (T,{\cdot })=0$. A $CR$-structure compatible with
$\xi $ is a smooth endomorphism $J:{\xi}$ ${\rightarrow}$ ${\xi}$ such that $%
J^{2}=-Identity$. We say $J$ is oriented if $(X,JX)$ is an oriented basis of
$\xi $ for any nonzero $X$ ${\in }$ ${\xi }$. A pseudohermitian structure
compatible with $\xi $ is a $CR$-structure $J$ compatible with $\xi $
together with a global contact form $\Theta $.

Given a pseudohermitian structure $(J,\Theta )$, we can
choose a complex vector field $Z_{1}$, an eigenvector of $J$ with eigenvalue
$i$, and a complex 1-form ${\theta }^{1}$ such that $\{\Theta ,{\theta ^{1}},%
{\theta ^{\bar{1}}}\}$ is dual to $\{T,Z_{1},Z_{\bar{1}}\}$ (${\theta ^{\bar{%
1}}}={\bar{({\theta ^{1}})}}$,$Z_{\bar{1}}={\bar{({Z_{1}})}}$). It follows
that $d\Theta =ih_{1{\bar{1}}}{\theta ^{1}}{\wedge }{\theta ^{\bar{1}}}$ for
some nonzero real function $h_{1{\bar{1}}}$. If both $J$ and $\Theta $%
are oriented, then $h_{1{\bar{1}}}$ is positive. In this case
we call such a pseudohermitian structure $(J,\Theta )$ oriented, and we can
choose a $Z_{1}$ (hence $\theta ^{1}$) such that $h_{1{\bar{1}}}=1$. That is
to say

\begin{align}
    \tag{A.1}  d\Theta =i{\theta ^{1}}{%
\wedge }{\theta ^{\bar{1}}}.
\end{align}

We will always assume our pseudohermitian structure is oriented and
$h_{1{\bar{1}}}=1$. The pseudohermitian
connection of $(J,\Theta )$ is the connection $\nabla ^{p.h.}$ on $TM{%
\otimes }C$ (and extended to tensors) given by

\begin{align}
    \tag{A.2} {\nabla}^{p.h.}Z_{1}={\omega
_{1}}^{1}{\otimes }Z_{1},{\nabla}^{p.h.}Z_{\bar{1}}={\omega _{%
\bar{1}}}^{\bar{1}}{\otimes }Z_{\bar{1}},{\nabla}^{p.h.}T=0
\end{align}

\noindent in which the 1-form ${\omega _{1}}^{1}$ is uniquely determined
by the following equation with a normalization condition ([Ta], [We]):

\begin{align}
    \tag{A.3} &&d{\theta ^{1}}={\theta ^{1}}{\wedge }{%
\omega _{1}}^{1}+{A^{1}}_{\bar{1}}\Theta {\wedge }{\theta ^{\bar{1}},}\\
\tag{A.4} &&{\omega _{1}}^{1}+{\omega _{\bar{1}}}^{%
\bar{1}}=0.
\end{align}

\noindent The coefficient ${A^{1}}_{\bar{1}}$ in $(A.3)$
is called the (pseudohermitian) torsion. Since $h_{1{\bar{1}}}=1$, $A_{{\bar{%
1}}{\bar{1}}}=h_{1{\bar{1}}}{A^{1}}_{\bar{1}}={A^{1}}_{\bar{1}}$. And $%
A_{11} $ is just the complex conjugate of $A_{{\bar{1}}{\bar{1}}}$.
Differentiating ${\omega _{1}}^{1}$ gives

\begin{align}
    \tag{A.5} d{\omega _{1}}^{1}=W{\theta
^{1}}{\wedge }{\theta ^{\bar{1} }} + 2iIm(A_{11,{\bar{1}}}{\theta ^{1}}{\wedge }%
\Theta )
\end{align}

\noindent where $W$ is the Tanaka-Webster curvature. Write ${\omega _{1}}%
^{1}=i{\omega }$ for some real 1-form $\omega $ by $(A.4)$.
This $\omega $ is just the one used in previous sections. Write $Z_{1}=\frac{%
1}{2}(e_{1}-ie_{2})$ for real vectors $e_{1},e_{2}$. It follows that $%
e_{2}=Je_{1}.$ Let $e^{1}=Re({\theta ^{1}}),e^{2}=Im({\theta ^{1}%
})$. Then $\{e^{0}={\Theta },e^{1},e^{2}\}$ is dual to $%
\{e_{0}=T,e_{1},e_{2}\}$. Now in view of $(A.1),(A.2)$ and $%
(A.3),$ we have the following real version of structure
equations:

\begin{align}
    \tag{A.1r} &d\Theta =2e^{1}\wedge e^{2},&\\
\tag{A.2r} &{\nabla }^{p.h.}e_{1}=\omega
\otimes e_{2},{\nabla }^{p.h.}e_{2}=-\omega \otimes
e_{1},&\\
\tag{A.3r} &de^{1}=-e^{2}\wedge \omega \mbox{ mod }\Theta ;\; \;
de^{2}=e^{1}\wedge \omega \mbox{ mod }\Theta .&
\end{align}

Similarly, from $(A.5)$, we have the following equation for $W$:

\begin{align}
    \tag{A.5r} d{\omega }(e_{1},e_{2})=-2W.
\end{align}

\noindent Also by $(A.1),(A.3)$ we can deduce

\begin{align}
    \tag{A.6} [Z_{\bar{1}},Z_{1}]&=iT+{\omega _{1}}%
^{1}(Z_{\bar{1}})Z_{1}-{\omega _{\bar{1}}}^{\bar{1}}(Z_{1})Z_{\bar{1}},\\
\tag{A.7} [Z_{\bar{1}},T]&={A^{1}}_{\bar{1}}Z_{1}-{%
\omega _{\bar{1}}}^{\bar{1}}(T)Z_{\bar{1}}.
\end{align}

\noindent The real version of $(A.6),(A.7)$ reads

\begin{align}
    \tag{A.6r} [e_{1},e_{2}]&=-2T-\omega
(e_{1})e_{1}-\omega (e_{2})e_{2},\\
\tag{A.7r} [e_{1},T]&=(\text{Re}A_{11})e_{1}-((\text{%
Im}A_{11})+\omega (T))e_{2},\\  \nonumber
[e_{2},T]&=-((\text{Im}A_{\bar{%
1}\bar{1}})+\omega (T))e_{1}+(\text{Re}A_{\bar{1}\bar{1}})e_{2}.
\end{align}

\noindent Note that $A_{\bar{1}\bar{1}}={A^{1}}_{\bar{1}}$ since $h_{1{%
\bar{1}}}=1.$ We define the subgradient operator $\nabla _{b}$ acting on a
smooth function $f$ by

\begin{align}
    \tag{A.8} \nabla _{b}f=2\{(Z_{\bar{1}%
}f)Z_{1}+(Z_{1}f)Z_{\bar{1}}\}.
\end{align}

\noindent It is easy to see that the definition of $\nabla _{b}$ is independent of the
choice of unitary ($h_{1{\bar{1}}}=1$) frame $Z_{1}$. The real
version of $(A.8)$ reads

\begin{align}
    \tag{A.8r} \nabla
_{b}f=(e_{1}f)e_{1}+(e_{2}f)e_{2}.
\end{align}

Next we will introduce the 3-dimensional Heisenberg group $H_{1}$
(see, e.g., [FS])$.$ For two points $(x,y,z),$ $(x^{\prime },y^{\prime
},z^{\prime })$ $\in $ $R^{3},$ we define the multiplication as follows: $%
(x,y,z)\circ (x^{\prime },y^{\prime },z^{\prime })$ $=$ $(x+x^{\prime },$ $%
y+y^{\prime },$ $z+z^{\prime }+yx^{\prime }-xy^{\prime }).$ $R^{3}$ endowed
with this multiplication ''$\circ "$ forms a Lie group, called the
(3-dimensional) Heisenberg group and denoted as $H_{1}.$ It is a simple
exercise to verify that:

$$\hat{e}_{1}=\frac{\partial }{\partial x}+y\frac{\partial }{\partial z}%
,\hat{e}_{2}=\frac{\partial }{\partial y}-x\frac{\partial }{\partial z}%
,T_{0}=\frac{\partial }{\partial z}$$

\noindent form a basis for the left-invariant vector fields on $H_{1}.$ We can endow $%
H_{1}$ with a standard pseudohermitian structure. The plane distribution
spanned by $\hat{e}_{1},\hat{e}_{2}$ forms a contact structure $\xi _{0}$
(so that $\hat{e}_{1},\hat{e}_{2}$ are Legendrian, i.e., sitting in the
contact plane). The $CR$ structure $J_{0}$ compatible with $\xi _{0}$ is
defined by $J_{0}(\hat{e}_{1})$ $=$ $\hat{e}_{2},$ $J_{0}(\hat{e}_{2})$ $=$ $%
-\hat{e}_{1}.$ The contact form $\Theta _{0}$ $=$ $dz+xdy-ydx$ gives rise to
the Reeb vector field $T_{0}=\frac{\partial }{\partial z}.$ Observe that $%
\{dx,dy,\Theta _{0}\}$ is dual to $\{\hat{e}_{1},\hat{e}_{2},T_{0}\}$ and $%
d(dx+idy)=0.$ It follows from the structural equations (A.3), (A.4) and
(A.5) that the connection form associated to the coframe $dx+idy$, the
torsion and the Tanaka-Webster curvature are all zero.

\bigskip




\end{document}